# The Complex Hyperbolic Geometry of the Moduli Space of Cubic Surfaces

Daniel Allcock, James A. Carlson, and Domingo Toledo

*To Herb Clemens on his 60th birthday*

**Contents**

1. **Introduction**

A classical theorem of great beauty describes the connection between cubic curves and hyperbolic geometry: the moduli space of the former is a quotient of the complex hyperbolic line (or real hyperbolic plane). The purpose of this paper is to exhibit a similar connection for cubic surfaces: their space of moduli is a quotient of complex hyperbolic four-space. We will make a precise statement below.

The theorem for cubic curves can be established using periods of integrals, or, in modern language, Hodge structures. Indeed, the classifying space for the Hodge structure on the first cohomology of a cubic curve is uniformized by the upper half plane, isomorphic to the complex hyperbolic line; the monodromy group of the universal family of cubic curves is isomorphic to $SL(2, \mathbb{Z})$; and the *period map* that assigns to a cubic curve its Hodge structure defines an isomorphism between the moduli space of smooth cubic curves and the quotient of the upper half plane by $PSL(2, \mathbb{Z})$. This isomorphism holds in several categories: analytic spaces, orbifolds, and quasi-projective algebraic varieties.

For cubic surfaces the second cohomology is generated by the classes of algebraic cycles, and so the natural period map is constant. Nonetheless, we can still use Hodge theory to study their moduli by considering instead the cohomology of a suitable cyclic branched cover of projective

First author partially supported by an NSF postdoctoral fellowship. Second and third authors partially supported by NSF grants DMS 9625463 and DMS 9900543. Third author partially supported by the IHES.

space. Constructions of this kind go back at least to Picard [42]. To describe more fully the one used here, consider a cubic surface $S$ in $\mathbb{P}^3$, and let $T$ be the triple cover of $\mathbb{P}^3$ branched along $S$. It is a cyclic cubic threefold, that is, a hypersurface in $\mathbb{P}^4$ of degree three invariant under an action of a certain cyclic group of order 3. We then establish the following. First, the map that assigns to a smooth cyclic cubic threefold $T$ its Hodge structure may be regarded as a map from the moduli space to a quotient of the complex hyperbolic 4-space $\mathbb{C}H^4$ by the action a certain discrete group $\Gamma$. Second, this group, which is the monodromy group for the middle cohomology of the threefold, is generated by complex reflections. Moreover, it admits a natural surjection to the Weyl group $W(E_6)$, the famous group of the 27 lines on a cubic surface, and its projectivization, written $P\Gamma$, acts on $\mathbb{C}H^4$ in the same way as does $PU(4,1,\mathcal{E})$. The latter is the projective unitary group of the standard unimodular Hermitian form of signature $(4,1)$ over the ring $\mathcal{E}$ of Eisenstein integers (the integers in $\mathbb{Q}(\sqrt{-3})$). Third, the period map that associates to a cubic surface $S$ the Hodge structure of $H^3(T)$ maps the moduli space of smooth cubic surfaces to the complement of an explicit closed subvariety of $P\Gamma\backslash\mathbb{C}H^4$, namely, the projection of the set of mirrors of the generating complex reflections of $\Gamma$. This map is an isomorphism of analytic spaces. Just as in the case of cubic curves, it is also an isomorphism of orbifolds. See §2 for the detailed statements, and §7 for the identification of $\Gamma$.

For a still more precise statement and a still closer analogy to the case of cubic curves, we consider in §3 cubic surfaces which are *stable* in the sense of geometric invariant theory. As is well known, these are either smooth or have nodes (ordinary double points) as singularities. Then our main theorem (3.17) establishes an analytic isomorphism between the moduli space of stable cubic surfaces and $P\Gamma\backslash\mathbb{C}H^4$. Since a cubic curve is stable if and only if it is non-singular, our result is strictly analogous to the classical one for cubic curves (compare $P\Gamma$ with $PSL(2,\mathbb{Z})$). The identification asserted in the category of analytic spaces can also be stated for schemes and orbifolds. For the latter one must use a structure on $P\Gamma\backslash\mathbb{C}H^4$ different from that which comes naturally from the underlying analytic space. See (3.18)–(3.20) for a discussion. One can also consider the geometric invariant theory compactification of the moduli space. The period map extends to this space and defines an analytic isomorphism with the Satake compactification of the ball quotient.

Although a holomorphic map from the moduli space of cubic surfaces to $PU(4,1,\mathcal{E})\backslash\mathbb{C}H^4$ which is generically of maximum rank appears already in [13], the monodromy group $\Gamma$ was not determined there. Its determination in §7 relies on results of Libgober [29] and techniques developed in [1] to study a class of groups that includes $\Gamma$. The isomorphism of the moduli space of smooth cubic surfaces with the quotient of the complement of the mirrors in $\mathbb{C}H^4$ by $\Gamma$ then follows from a combination of fairly standard techniques in complex geometry and the Torelli theorem of Clemens and Griffiths [15] for cubic threefolds as described in our announcement [3].

The proof of the more refined theorem on moduli of stable cubic surfaces follows the general lines of the argument in [3], except that new technical subtleties arise in proving that the period map is an isomorphism on the divisor corresponding to nodal surfaces. We found it best to establish the theorem by using the identification of the latter divisor with the moduli space of six points in $\mathbb{P}^1$ in combination with the theorem of Deligne and Mostow which gives a complex hyperbolic structure to this space [20]. The resulting complex hyperbolic manifold is easily identified with the divisor of mirrors in $P\Gamma\backslash\mathbb{C}H^4$. See §9 for details on this point, and see (3.22) for a sketch of the proof of the main theorem. For a more elementary approach to establishing the isomorphism on the nodal divisor, one could further develop the calculus of fractional differentials of §6 so as to include a proof of Lemma 9.1 on all equisingular strata. We plan to use this approach in future papers devoted to the complex hyperbolic geometry of the moduli spaces of del Pezzo surfaces and cubic threefolds.

By Theorem 8.4 of [1], the monodromy group $\Gamma$ is not one of the discrete groups considered by Mostow [37] or Deligne and Mostow [20]; however, it contains some of their groups in a natural way. We have already appealed to the fact that divisor of nodal surfaces is uniformized by the first group of the list for $N = 6$ in §14.4 of [20], which is the same as the first group of Thurston's list [47]. Similarly, the codimension two subspace of surfaces with two nodes is uniformized by the first group of their list for $N = 5$, which is the second group on Thurston's list and is also the group studied by Picard in [42].

In addition to the proof of the main theorem we have included a few other results. In §7 we study the configuration of hyperplanes in $\mathbb{C}H^4$ consisting of the mirrors of $\Gamma$. It has the remarkable property that whenever two hyperplanes intersect, they do so at right angles. This allows one to prove that the moduli space of smooth cubic surfaces has contractible universal cover (see [2]). We also show that the set of mirrors in $\mathbb{C}H^4$ splits into 36 families indexed in a natural way by the 36 reflections of the Weyl group. Any two mirrors in one of these families are disjoint. See (7.28)–(7.30) for more details. From these facts it is easy to identify the toroidal compactification of $P\Gamma_\theta\backslash\mathbb{C}H^4$ with Naruki's cross-ratio variety [40], where $\Gamma_\theta$ is the congruence subgroup of $\Gamma$ introduced in (3.12). We do not, however, pursue these details.

In Theorem 10.2 we establish a simple but remarkable fact: the universal smooth cubic surface lies naturally in the projective tangent bundle to the moduli space of smooth cubic surfaces. As a consequence we are able in §11 to find the points of $\mathbb{C}H^4$ corresponding to the most symmetric surfaces. This suggests the most natural way to prove the theorem: invert the period map by finding automorphic forms on $\mathbb{C}H^4$ realizing the universal surface in the tangent bundle to the corresponding quotient of $\mathbb{C}H^4$. See [5] for a different family of automorphic forms, yielding an embedding of the space of marked cubic surfaces into $\mathbb{C}P^9$.

The results presented here grew out of our interest in the fundamental group of the space of smooth cubic surfaces. It had been studied earlier by Libgober [29], who discovered that it was a quotient of the Artin group of $E_6$, and much more recently by Looijenga [31], who found a presentation for it. Looijenga's work, which is a major step forward, solves in principle the problem of determining the kernel of Libgober's surjection from the Artin group of $E_6$. Our methods give a different insight into the structure of the fundamental group of the space of smooth cubic surfaces, one based on the fact that it is commensurable with the fundamental group of the complement of a totally geodesic divisor in a locally symmetric variety. These results therefore follow the pattern established in [34] and [48]. In [4] we will use this fact to prove that this fundamental group is not a lattice in any Lie group. In particular, the moduli space of smooth cubic surfaces is not the quotient of any bounded symmetric domain by any lattice. In other words, removal of a divisor from the locally symmetric variety is essential to any description of the moduli space of smooth surfaces.

We remark that moduli space of smooth cubic surfaces is the complement of a totally geodesic divisor in a locally symmetric variety in more that one way. Using the period map for configurations of six lines in $\mathbb{P}^2$ [33], one sees that our moduli space is also the complement of a totally geodesic divisor in a locally symmetric variety for $SO(2,4)$. The present situation is thus analogous to that of the moduli space for genus 2 curves, which is the complement of totally geodesic divisor in two different locally symmetric varieties, one for $SU(3,1)$, and another for $Sp(4,\mathbb{R})$; see [28].

An alternative approach to our theorem, developed by Hunt and van Geemen, is discussed in [27]. The differential equations satisfied by our period map (analogous to the hypergeometric equation) have been obtained by Sasaki and Yoshida [43].

We would like to explicitly record here our indebetedness to two papers that have been essential to this work, both in terms of their technical content and in terms of foundation and inspiration that they provide for further developments: the papers of Clemens and Griffiths [15] and of Deligne and Mostow [20]. We also thank H. Clemens, M. Kapovich, J. Kollár, Paul Roberts and particularly E. Looijenga for very useful discussions.

## 2. Moduli of smooth cubic surfaces

The aim of this section is to construct the period map for the moduli space of smooth cubic surfaces and to state our main results concerning this space. In the next section we will give more complete results which also treat the structure of the moduli space of stable cubic surfaces. We have gathered together the proofs of various lemmas in section 4.

(2.1) The period map is defined by associating to a cubic surface $S \subset \mathbb{P}^3$ the Hodge structure of a cyclic triple cover $T$ of $\mathbb{P}^3$ which is branched along $S$. For technical reasons we need to consider not only the surface $S$ but also a cubic form $F$ defining it. Let $(X_0, X_1, X_2, X_3)$ denote homogeneous coordinates for $\mathbb{P}^3$ and let $(X_0, X_1, X_2, X_3, Y)$ denote homogeneous coordinates for $\mathbb{P}^4$. We write $\mathcal{C}$ for the set of all nonzero cubic forms in $X_0, \ldots, X_3$. Each $F \in \mathcal{C}$ defines a cubic surface $S$ in $\mathbb{P}^3$ and also a threefold

$$T = \{(X_0, \ldots, X_3, Y) \in \mathbb{P}^4 : Y^3 - F(X_0, X_1, X_2, X_3) = 0\} \tag{2.1.1}$$

in $\mathbb{P}^4$. This is the triple cover of $\mathbb{P}^3$ branched over $S$, and we denote by $p : T \longrightarrow \mathbb{P}^3$ the covering map. *Throughout the paper, whenever we have a cubic form $F$ in mind we implicitly define the surface $S$ and threefold $T$ in this manner.* We will also use the analogous implicit definitions for $S'$ and $T'$ in terms of a form $F' \in \mathcal{C}$, etc. Fix a primitive cube root of unity $\omega$, and let $\sigma$ be the automorphism of $\mathbb{P}^4$ defined by

$$\sigma(X_0, X_1, X_2, X_3, Y) = (X_0, X_1, X_2, X_3, \omega Y) \, .$$

It leaves $T$ invariant and generates its group of branched covering transformations over $\mathbb{P}^3$. These choices of $\mathbb{P}^3$, $\mathbb{P}^4$ and $\sigma$ will be fixed throughout the paper, and $T$ will be called the cyclic cubic threefold associated to $F$. Let $\mathcal{C}_0$ denote the subspace of $\mathcal{C}$ consisting of forms which define surfaces $S$ in $\mathbb{P}^3$ which are smooth (as schemes). Thus $\mathcal{C}_0$ is the complement of the discriminant hypersurface $\Delta$. In order to define the period map we need to study the cohomology of $T$ for $F \in \mathcal{C}_0$, the corresponding monodromy over $\mathcal{C}_0$, and the way the Hodge structure of $T$ varies with $F$.

(2.2) We begin with the study of cohomology. First, $H^3(T, \mathbb{Z})$ is a free $\mathbb{Z}$-module of rank 10. This is standard and can be checked by computing Euler characteristics. See, for example, formula 8.2 of [13]. Since $H^3(\mathbb{P}^3, \mathbb{Z}) = 0$, the induced transformation $\sigma^*$ on $H^3(T, \mathbb{Z})$ fixes no vector except 0. (We will usually suppress the asterisk and write $\sigma$ for $\sigma^*$.) Thus its minimal polynomial is $t^2 + t + 1$, which is also the minimal polynomial of $\omega$. Taking $\omega$ to act as $\sigma$, and setting $\mathcal{E} = \mathbb{Z}[\omega]$, the ring of Eisenstein integers, we see that $H^3(T, \mathbb{Z})$ carries the structure of an $\mathcal{E}$-module. This $\mathcal{E}$-module, which we denote by $\Lambda(T)$, is free of rank five. Freeness follows from the fact that the ring of Eisenstein integers is a principal ideal domain.

(2.3) The symplectic form $\Omega$ on $H^3(T, \mathbb{Z})$, given by the cup product and evaluation on the fundamental cycle, defines a $\mathbb{Z}$-bilinear form $h$ on $\Lambda(T)$ by the formula

$$h(x, y) = -\frac{\Omega(\theta x, y) + \theta \Omega(x, y)}{2} \ . \tag{2.3.1}$$

This formula is interpreted as follows: the first $\theta$ is the endomorphism $\sigma - \sigma^{-1}$ and the second $\theta$ is the important Eisenstein integer $\omega - \omega^{-1} = \sqrt{-3}$. Theorem 4.1 shows that $h$ is actually a unimodular hermitian form over $\mathcal{E}$, $\mathcal{E}$-linear in its first argument and antilinear in its second. If $v \in \Lambda(T)$ then we define the norm of $v$ as $h(v, v)$ .

(2.4) To compute the signature of $h$, we consider the hermitian form on $H^3(T, \mathbb{C})$ defined by

$$h'(\alpha, \beta) = \theta \int_T \alpha \wedge \bar{\beta}.$$

With respect to it we have an orthogonal decomposition

$$H^3(T, \mathbb{C}) = H^3_\omega(T) \oplus H^3_{\bar\omega}(T) \tag{2.4.1}$$

into the eigenspaces of $\sigma$ corresponding to its eigenvalues $\omega$ and $\bar\omega$ repectively. Since $\sigma$ is defined on the real vector space $H^3(T, \mathbb{R})$, the two eigenspaces are exchanged by complex conjugation and so both are five-dimensional. Now consider the map

$$Z = \pi_{\bar\omega} \circ \iota : \Lambda(T) \longrightarrow H^3_{\bar\omega}(T) \tag{2.4.2}$$

where $\iota : \Lambda(T) \to \Lambda(T) \otimes_{\mathbb{Z}} \mathbb{C} \cong H^3(T, \mathbb{C})$ is the natural inclusion and $\pi_{\bar\omega}$ is the projection $H^3(T, \mathbb{C}) \to H^3_{\bar\omega}(T)$.

**(2.5) Lemma.** *$Z$ is an isometric embedding.*

The proof appears in (4.3). Hereafter we will identify $\Lambda(T)$ with the lattice $Z(\Lambda(T)) \subset H^3_{\bar\omega}(T)$. We also write $Z$ for the extension of this map $\Lambda(T) \otimes_{\mathcal{E}} \mathbb{C} \to H^3_{\bar\omega}$ and observe that another way to state Lemma 2.5 to say that

$$Z : (\Lambda(T) \otimes_{\mathcal{E}} \mathbb{C}, \sigma, h) \to (H^3_{\bar\omega}(T), \bar\omega, h') \qquad (2.5.1)$$

is an isomorphism of all the displayed structures. It follows that the signature of $h$ on $\Lambda(T)$ is the same as signature of $h'$ on $H^3_{\bar\omega}$. The signature of $h'$ on $H^3_{\bar\omega}$ can be determined by Hodge-theoretic calculations, carried out in (4.4) with the following result.

**(2.6) Lemma.** $H^3_{\bar\omega} = H^{2,1}_{\bar\omega} \oplus H^{1,2}_{\bar\omega}$; the first summand is one-dimensional, and $h'$ is negative-definite there; the second sumand is four-dimensional, and $h'$ is positive-definite there. In particular, the signatures of $h$ on $\Lambda(T)$ and $h'$ on $H^3_{\bar\omega}$ are $(4, 1)$.

(2.7) According to Theorem 7.1 of [1], there is a unique isometry class of five-dimensional free $\mathcal{E}$-modules endowed with a unimodular hermitian form of signature $(4, 1)$. Let $\Lambda$ be the standard such lattice $\mathcal{E}^{4,1}$, namely the free module $\mathcal{E}^5$ with the hermitian form $h$ defined by

$$h(x, y) = -x_0 \bar{y}_0 + x_1 \bar{y}_1 + \cdots + x_4 \bar{y}_4. \qquad (2.7.1)$$

Then $\Lambda(T)$ is isometric to $\Lambda$. Moreover, $\Lambda \otimes_{\mathcal{E}} \mathbb{C}$ is $\mathbb{C}^{4,1}$, the complex five-space endowed with the hermitian form given by the same formula. It follows that the pair $(H^3_{\bar\omega}(T), \Lambda(T))$ is isometric to the pair $(\mathbb{C}^{4,1}, \Lambda)$.

(2.8) We are now in a position to study the monodromy of $\Lambda(T)$ over the space of cubic forms. Over $\mathcal{C}$ is defined a universal family $\mathcal{S} \subset \mathcal{C} \times \mathbb{P}^3$ of surfaces

$$\mathcal{S} = \{(F, (X_0, X_1, X_2, X_3)) \in \mathcal{C} \times \mathbb{P}^3 : F(X_0, X_1, X_2, X_3) = 0\} \qquad (2.8.1)$$

and a universal family of cyclic cubic threefolds

$$\mathcal{T} = \{(F, (X_0, X_1, X_2, X_3, Y)) \in \mathcal{C} \times \mathbb{P}^4 : Y^3 - F(X_0, X_1, X_2, X_3) = 0\}. \qquad (2.8.2)$$

We write $\pi$ for either of the projections $\mathcal{S} \to \mathcal{C}$ or $\mathcal{T} \to \mathcal{C}$, since context will identify the intended map. Note that (2.8.2) is the global version of (2.1.1). Observe also that the total spaces $\mathcal{S}$ and $\mathcal{T}$ are smooth varieties. This is easily checked by differentiating the defining equations of $\mathcal{S}$ and $\mathcal{T}$ with respect to the coefficients of $F$ and observing that the simultaneous vanishing of all partial derivatives with respect to these coefficients implies that $X_0 = \cdots = X_3 = 0$. We write $\mathcal{S}_0$ and $\mathcal{T}_0$ for the locally trivial fibrations which are the restrictions of $\mathcal{S}$ and $\mathcal{T}$ to $\mathcal{C}_0$.

(2.9) The family $\mathcal{T}_0$ gives rise to the sheaf $R^3\pi_*(\mathbb{Z})$ over $\mathcal{C}_0$. Recall that it is the sheaf associated to the presheaf $U \to H^3(\pi^{-1}(U), \mathbb{Z})$. The automorphism $\sigma$ acts on it, and we claim that 1 is not an eigenvalue for any $U$. Indeed, if $U$ is contractible then $H^3(U \times \mathbb{P}^3, \mathbb{Z}) = 0$, and so $\sigma$ fixes no elements of $H^3(\pi^{-1}(U), \mathbb{Z})$ other than zero. Reasoning as in (2.2), we see that $R^3\pi_*(\mathbb{Z})$ is a sheaf of $\mathcal{E}$-modules. Defining a hermitian form $h$ as in (2.3.1), we see that $R^3\pi_*(\mathbb{Z})$ is a sheaf over $\mathcal{C}_0$ of unimodular hermitian $\mathcal{E}$-modules. We denote it by $\Lambda(\mathcal{T}_0)$. Let $H^3(\mathcal{T})$ denote the sheaf $R^3\pi_*(\mathbb{C})$ over $\mathcal{C}$, and denote by $H^3_{\bar\omega}(\mathcal{T})$ the subsheaf of $H^3(\mathcal{T})$ consisting of eigenvectors with eigenvalue $\bar\omega$. Over $\mathcal{C}_0$ the sheaf $H^3_{\bar\omega}(\mathcal{T}_0)$ is a local system vector spaces, hermitian of signature $(4, 1)$.

(2.10) Fix a basepoint $F_0 \in \mathcal{C}_0$. Its associated threefold $T_0$ is the fiber of $\mathcal{T}$ over $F_0$. Then we have the monodromy representation

$$\rho_0 : \pi_1(\mathcal{C}_0, F_0) \to \operatorname{Aut}(\Lambda(T_0)).$$

If $P$ denotes projectivization (passage to the associated projective space or projective linear group), we obtain a representation

$$P\rho_0 : \pi_1(\mathcal{C}_0, F_0) \to P\operatorname{Aut}(\Lambda(T_0)).$$

Let $\Gamma_0$ and $P\Gamma_0$ denote the images of $\rho_0$ and $P\rho_0$, and let $\tilde{\mathcal{C}}_0$ denote the covering space of $\mathcal{C}_0$ corresponding to $\ker(P\rho_0)$. Then $\tilde{\mathcal{C}}_0$ is a Galois cover of $\mathcal{C}_0$ with Galois group $P\Gamma_0$. In (3.9) we will give a convenient concrete model of this covering space. Over $\tilde{\mathcal{C}}_0$ the flat bundle of projective spaces obtained by pulling back $P(H^3_\omega)(\mathcal{T})$ is trivial: it is canonically isomorphic (by parallel translation) to the base $\tilde{\mathcal{C}}_0$ times the fiber at one point.

(2.11) It will be convenient to have a standard model for the monodromy representation and for the period map which will be introduced below. To this end, let us fix once and for all an isometry $\phi : \Lambda(T_0) \to \Lambda$. This choice induces an isomorphism of $\Gamma_0$ with the subgroup $\Gamma = \phi \Gamma_0 \phi^{-1}$ of $\operatorname{Aut}(\Lambda)$, which we call *the monodromy group*, and it induces a representation $\rho : \pi_1(\mathcal{C}_0, F_0) \to \operatorname{Aut}(\Lambda)$ with image $\Gamma$ which we call *the monodromy representation*. It is clear that $\rho$ depends (by conjugacy) on the choice of $\phi$; however, we will see in (2.14) that $\Gamma$ is independent of this choice.

(2.12) To describe $\Gamma$ precisely we will make use of the finite vector space obtained by reducing $\Lambda$ modulo $\theta = \sqrt{-3}$. It is easy to see that $\mathcal{E}/\theta\mathcal{E}$ is the field $\mathbb{F}_3$ of three elements. Thus $\Lambda/\theta\Lambda$ is a five-dimensional vector space over $\mathbb{F}_3$ which we denote by $V$. Since $\mathbb{F}_3$ has no non-trivial automorphisms, the hermitian form $h$ reduces to a quadratic form $q$ on $V$. Since $h$ is unimodular, $q$ is non-degenerate. We will never consider automorphisms of $V$ that do not preserve $q$, so we sometimes write $\operatorname{Aut}(V)$ for $\operatorname{Aut}(V, q)$. This is an orthogonal group, and the spinor norm homomorphism $\nu$ defined on it takes values in $\mathbb{F}_3^*/(\mathbb{F}_3^*)^2 = \{\pm 1\} \cong \mathbb{Z}/2\mathbb{Z}$. See §55 of [41] for background. Explictly, if $g \in \operatorname{Aut} V$, then $g$ may be expressed as a product of reflections in vectors $v_1, \cdots, v_k \in V$, and its spinor norm $\nu(g)$ is defined to be the square class of $q(v_1) \cdots q(v_k)$. We therefore have homomorphisms

$$\operatorname{Aut}(\Lambda, h) \to \operatorname{Aut}(V, q) \xrightarrow{\nu} \mathbb{Z}/2\mathbb{Z}. \tag{2.12.1}$$

By Lemma 4.5, both of these maps are surjective. We define $\operatorname{Aut}^+(V)$ to be the kernel of $\nu$ and $\operatorname{Aut}^+(\Lambda)$ to be the kernel of the composition. By the lemma these subgroups have index two in $\operatorname{Aut}(V)$ and $\operatorname{Aut}(\Lambda)$.

(2.13) We can now state our first main result: the identification of $\Gamma$. The proof occupies most of section 7. If $\zeta$ is a root of unity then the $\zeta$-reflection in a vector $v$ with $h(v, v) \neq 0$ is the linear isometry which multiplies $v$ by $\zeta$ and fixes $v^\perp$ pointwise. An explicit formula for this map is

$$x \mapsto x - (1-\zeta)\frac{h(x, v)}{h(v, v)} .$$

These transformations are called complex reflections, and complex reflections of orders 2, 3 and 6 are sometimes called biflections, triflections, and hexflections. A vector of $\Lambda$ with norm 1 (resp. 2) is called a short (resp. long) root, and it is easy to see that hexflections in short roots are isometries of $\Lambda$, as are biflections in long roots.

**(2.14) Theorem.** *The monodromy group $\Gamma$ is the group $\mathrm{Aut}^+(\Lambda)$. It coincides with the subgroup of $\mathrm{Aut}(\Lambda)$ generated by the hexflections in the short roots of $\Lambda$. Moreover, the inclusion of $\Gamma$ in $\mathrm{Aut}^+(\Lambda)$ induces an isomorphism $P\Gamma \cong P\mathrm{Aut}(\Lambda)$.*

The normality of $\mathrm{Aut}^+(\Lambda)$ in $\mathrm{Aut}(\Lambda)$ justifies our claim in (2.11) that $\Gamma$ is independent of $\phi$.

(2.15) Next we will define the period map from $\tilde{\mathcal{C}}_0$ to complex hyperbolic 4-space. If $W$ is a complex vector space of dimension $n+1$ with a hermitian form $h$ of signature $(n,1)$, we let $\mathbb{C}H(W)$, the *complex hyperbolic space* of $W$, denote the space of negative lines: the space of one-dimensional subspaces of $W$ on which $h$ is negative definite. The space $\mathbb{C}H(W)$ is an open subset of $P(W)$ and is biholomorphic to the unit ball in $\mathbb{C}^n$. If $W = \mathbb{C}^{n,1}$ we write $\mathbb{C}H^n$ for $\mathbb{C}H(W)$.

(2.16) For each $F \in \mathcal{C}_0$, we have seen in Lemma 2.6 that $H_{\bar\omega}^{2,1}(T) \subset H_{\bar\omega}^3(T)$ is a negative line, hence an element of $\mathbb{C}H(H_{\bar\omega}^3(T))$. By (2.10), the flat bundle $\mathbb{C}H(H_{\bar\omega}^3(T))$ trivializes over $\tilde{\mathcal{C}}_0$, so by parallel translation we can regard the collection of all the $H_{\bar\omega}^{2,1}(T)$ as points in $\mathbb{C}H(H_{\bar\omega}^3(T_0))$. The assignment

$$\left(\tilde{F} \in \tilde{\mathcal{C}}_0\right) \mapsto \text{ parallel translate of } H_{\bar\omega}^{2,1}(T) \in \mathbb{C}H(H_{\bar\omega}^3(T_0))$$

defines a map

$$g_0 : \tilde{\mathcal{C}}_0 \to \mathbb{C}H(H_{\bar\omega}^3(T_0)) \ .$$

It is $P\rho_0$-equivariant, and it is holomorphic since the Hodge filtration varies holomorphically (see also (6.7)). Our choice of $\phi : \Lambda(T_0) \to \Lambda$ in (2.11) induces an isometry $\mathbb{C}H(H_{\bar\omega}^3(T_0)) \to \mathbb{C}H^4$. Composing with $g_0$ yields a holomorphic map

$$g : \tilde{\mathcal{C}}_0 \to \mathbb{C}H^4$$

which we call the period map. Formulas for this map will be given in (3.11) and (6.7).

(2.17) Finally, we turn our attention to the moduli space of cubic surfaces. The general linear group $GL(4,\mathbb{C})$ operates on the left on $\mathcal{C}$ in the standard way: $(gF)(X) = F(g^{-1}X)$. This induces standard actions of $GL(4,\mathbb{C})$ on $\mathcal{S}$ and $\mathcal{T}$ by $g(F,X) = (gF, gX)$ and $g(F,(X,Y)) = (gF,(gX,Y))$, respectively. Let $D \subset GL(4,\mathbb{C})$ be the subgroup consisting of scalar matrices whose diagonal entries are a cube root of unity. Then $D$ is a central subgroup, cyclic of order three, which acts trivially on $\mathcal{C}$ and also on $\mathcal{S}$. The action of $D$ on $\mathcal{T}$ can be described by observing that each element of $D$ acts on the $\mathbb{P}^4$ containing the cubic threefolds in the same manner as some power of the branched covering transformation $\sigma$. Let $G = GL(4,\mathbb{C})/D$. Then $G$ acts effectively on $\mathcal{C}$ and on $\mathcal{S}$, preserving the subspace $\mathcal{C}_0$. We will see in (3.9) that the action on $\mathcal{C}_0$ lifts to an action of $G$ on $\tilde{\mathcal{C}}_0$.

(2.18) It is known (see (3.1)) that $G$ acts properly and with finite isotropy groups on $\mathcal{C}_0$. Thus the quotient $G\backslash\mathcal{C}_0$ is an analytic space and indeed a complex analytic orbifold. It is also true (3.14) that $G$ acts freely and properly on $\tilde{\mathcal{C}}_0$, so that $G\backslash\tilde{\mathcal{C}}_0$ is a complex analytic manifold and the orbit map $\tilde{\mathcal{C}}_0 \to G\backslash\tilde{\mathcal{C}}_0$ is a principal $G$-bundle. Our moduli spaces are

$$M_0 = G\backslash\mathcal{C}_0 \quad \text{and} \quad M_0^f = G\backslash\tilde{\mathcal{C}}_0.$$

$M_0$ is the moduli space of smooth cubic surfaces, and $M_0^f$ will be discussed in (3.15) as the moduli space of "framed" smooth cubic surfaces. Since a bounded holomorphic function on $G$ is constant, the period map $g$ is constant on $G$-orbits. Consequently it descends to a holomorphic map

$$g : M_0^f = G\backslash\tilde{\mathcal{C}}_0 \to \mathbb{C}H^4.$$

By the equivariance of the $P\Gamma$-actions on domain and range, this descends further, to a holomorphic map

$$g : M_0 = G\backslash\mathcal{C}_0 \to P\Gamma\backslash\mathbb{C}H^4.$$

(2.19) We can now state our main theorem concerning moduli of smooth cubic surfaces. This is a special case of our more general results for stable cubic surfaces, which appear in the next section. Write $\mathcal{R}$ for the set of short roots of $\Lambda$, and observe that for each $r \in \mathcal{R}$ there is a natural inclusion $\mathbb{C}H(r^\perp) \subset \mathbb{C}H^4$ induced by the inclusion of $r^\perp$ into $\Lambda$. The short roots define the hyperplane arrangement

$$\mathcal{H} = \bigcup_{r \in \mathcal{R}} \mathbb{C}H(r^\perp) \subset \mathbb{C}H^4 \ ,$$

and it is obvious that $\mathcal{H}$ is preserved by $P\Gamma$. Of course $\mathbb{C}H^4 - \mathcal{H}$ is a complex manifold, and we regard $P\Gamma\backslash(\mathbb{C}H^4 - \mathcal{H})$ as equipped with the orbifold structure that arises naturally from its description as a quotient of a manifold by a discrete group.

(2.20) **Theorem.** *The period map $g$ induces isomorphisms $M_0^f \longrightarrow \mathbb{C}H^4 - \mathcal{H}$ and $M_0 \longrightarrow P\Gamma\backslash(\mathbb{C}H^4 - \mathcal{H})$. The first is an isomorphism of complex manifolds and the second is an isomorphism of complex analytic orbifolds.*

3. **Moduli of stable cubic surfaces**

To prove Theorem 2.20 and its extensions, Theorems 3.17 and 3.20, we consider the moduli of cubic surfaces equipped with extra structure, as well as suitable compactifications of these spaces. In this section we introduce these structures and state the main results of the paper, which concern the moduli space of stable cubic surfaces and its Geometric Invariant Theory (GIT) compactification. As in section 2 we have gathered the proofs of various lemmas in section 4.

(3.1) To compactify the moduli spaces of section 2, we will need the notions of stability and semi-stability from GIT, specialized to cubic forms in four variables. Background for the claims here can be found in [38], particularly p. 80, and in [39], particularly p. 51. A form is *stable* if its orbit under $SL(4,\mathbb{C})$ is closed and its isotropy group is finite. A cubic form is stable if and only if its zero scheme is either smooth or has only ordinary double points. (Ordinary double points, also called nodes or $A_1$ singularities, are locally analytically equivalent to the singularity defined by $x^2 + y^2 + z^2 = 0$.) We write $\mathcal{C}_s$ for the space of stable cubic forms and $\Delta_s$ for $\Delta \cap \mathcal{C}_s$. A form is called semi-stable if the closure of its $SL(4,\mathbb{C})$-orbit does not contain 0. A cubic form is semi-stable if and only if its zero scheme has no singularities besides nodes and cusps. (The latter are locally analytically equivalent to $x^2 + y^2 + z^3 = 0$ and are also called $A_2$ singularities.) For a proof that these are the semi-stable cubic forms in four variables see §19 of [26]. We write $\mathcal{C}_{ss}$ for the space of semi-stable cubic forms and $\Delta_{ss}$ for $\Delta \cap \mathcal{C}_{ss}$. The only closed $SL(4,\mathbb{C})$-orbit in $P(\Delta_{ss} - \Delta_s)$ is the orbit of $X_0^3 - X_1 X_2 X_3$; therefore the closure of the orbit of any semistable cubic surface with a cusp contains this particular tricuspidal surface. Although this last fact is known (see for example the introduction to [40]), we have not found a proof in the literature and therefore provide one in (4.6).

(3.2) We now develop the notion of a marking of a cubic surface in terms of cohomology. It is equivalent to the one used by Naruki in [40]. Let $L$ denote the lattice $\mathbb{Z}^{1,6}$, namely $\mathbb{Z}^7$ with the bilinear form $x_0 y_0 - x_1 y_1 - \cdots - x_6 y_6$. Let $\eta \in L$ denote the norm 3 vector $(-3, 1, 1, 1, 1, 1, 1)$. It is known [32] that if $S$ is a smooth cubic surface, then the lattice $L(S) = H^2(S, \mathbb{Z})$ is isometric to $L$ by an isometry that takes the hyperplane class $\eta(S)$ to $\eta$. Define a marking of a cubic form $F \in \mathcal{C}_0$ to be an isometry $m : L(S) \to L$ which takes $\eta(S)$ to $\eta$. We can now define the space $\mathcal{M}_0$ of marked smooth cubic forms. As a set, it is the collection of all markings of all smooth cubic forms. To give it the structure of a complex manifold, let $L(\mathcal{S}_0)$ be the sheaf $R^2 \pi_*(\mathbb{Z})$ over $\mathcal{C}_0$, where $\pi : \mathcal{S}_0 \to \mathcal{C}_0$ is the universal smooth cubic surface of (2.8). Observe that the stalk of $L(\mathcal{S}_0)$ over $F$ is canonically isomorphic to $L(S)$. The set $\mathcal{M}_0$ is therefore in natural one-to-one correspondence with the subsheaf of $\text{Hom}(L(\mathcal{S}_0), \mathcal{C}_0 \times L)$ consisting of homomorphisms which carry each stalk isometrically to $L$ in such a way as to identify the hyperplane class with $\eta$. This one-to-one correspondence defines the structure of a complex manifold on $\mathcal{M}_0$. We write $p : \mathcal{M}_0 \longrightarrow \mathcal{C}_0$ for the natural projection. Now observe that the action of $G$ on forms defined in (2.17) extends naturally to $\mathcal{M}_0$: if $g \in G$ and $m \in \mathcal{M}_0$ is a marking of $F$, then $m \circ g^*$ is a marking of $gF$. Observe also that the group $\text{Aut}(L, \eta)$, which is isomorphic to the Weyl group $W(E_6)$, acts on $\mathcal{M}_0$: if $w \in \text{Aut}(L, \eta)$ and $m \in \mathcal{M}_0$ is a marking of $F$, then $w \circ m$ is also a marking of $F$. These two group actions commute. It is clear that $\mathcal{M}_0$ is a covering space of $\mathcal{C}_0$ with Galois group $\text{Aut}(L, \eta) \cong W(E_6)$. It is classical that the monodromy of $\mathcal{M}_0$ on $L(\mathcal{S}_0)$ over some $F \in \mathcal{C}_0$ realizes every isometry of $(L(S), \eta(S))$. Therefore $\mathcal{M}_0$ coincides with the covering space of $\mathcal{C}_0$ corresponding to the subgroup of $\pi_1(\mathcal{C}_0)$ that is the kernel of the monodromy of $L(\mathcal{S}_0)$. Thus $\mathcal{M}_0$ is connected.

(3.3) We define the space $\mathcal{M}_s$ of marked stable cubic forms to be the Fox completion of $\mathcal{M}_0$ over $\mathcal{C}_s$. More precisely, in the terminology of [22], $p : \mathcal{M}_s \longrightarrow \mathcal{C}_s$ is the completion of the spread $p : \mathcal{M}_0 \longrightarrow \mathcal{C}_s$, where $p$ is as defined in (3.2). Formally, a point of $\mathcal{M}_s$ lying over $F \in \mathcal{C}_s$ is a function $m$ which assigns to each neighborhood $W$ of $F$ a connected component $m(W)$ of $p^{-1}(W - \Delta)$ in such a way that if $W' \subset W$ then $m(W') \subset m(W)$. It is clear that $m$ is determined by its values at connected neighborhoods of $F$. In our setting it is also true that $m$ is determined by its value on a fixed suitably chosen neighborhood of $F \in \mathcal{C}$. Indeed, let $W$ be a connected neighborhood of $F \in \mathcal{C}_s$ with the property that for any connected neighborhood $W'$ of $F$, $W' \subset W$, the inclusion $W' - \Delta \subset W - \Delta$ induces a surjection on fundamental groups. Then $m(W)$ is a component of the preimage of $W - \Delta$; if $W' \subset W$ is any connected neighborhood of $F$ then the preimage of $W'$ in $m(W)$ is connected, so that $m(W')$ is uniquely determined. The existence of such $W$ follows from the discussion below.

(3.4) In the same manner we define the space $\mathcal{M}_{ss}$ of marked semi-stable forms to be the Fox completion of $\mathcal{M}_0$ over $\mathcal{C}_{ss}$. This space naturally contains $\mathcal{M}_s$, and these completions carry natural topologies (see §2 of [22]) and natural complex manifold structures (see (3.6) and (3.8) below). We will sometimes refer to a point of $\mathcal{M}_{ss}$ lying over $F \in \mathcal{C}_{ss}$ as a marking of $F$. Because of the intrinsic nature of the Fox completion, the actions of $G$ and $\operatorname{Aut}(L, \eta) = W(E_6)$ on $\mathcal{M}_0$ extend to $\mathcal{M}_s$ and $\mathcal{M}_{ss}$.

(3.5) To understand the Fox completions we must understand the inclusions $\Delta_s \to \mathcal{C}_s$ and $\Delta_{ss} \to \mathcal{C}_{ss}$ as well as the local monodromy of $\mathcal{M}_0$ near a point of $\Delta_s$ or $\Delta_{ss}$. Let us begin by studying $\Delta_s$. A cubic surface can have at most four nodes, and the space $\Delta_s^k$ of cubic forms with exactly $k \geq 1$ nodes and no other singularities has codimension $k$ in $\mathcal{C}$. If $F_0 \in \Delta_s^k$, then $\Delta$ is a locally normal crossing divisor at $F_0$. That is, there are local analytic coordinates $z_1, \cdots, z_{20}$ for $\mathcal{C}$ centered at $F_0$ such that $\Delta$ has local equation $z_1 \cdots z_k = 0$ at $F_0$. In particular, $\Delta_s^1$ is smooth; it also accounts for all of the smooth points of $\Delta$. These facts can be checked by referring to the normal forms and slices described in [10]. It follows from this description that any $F_0 \in \Delta_s^k$ has a neighborhood $W$ with $\pi_1(W - \Delta) \cong \mathbb{Z}^k$, with natural generators for the local fundamental group defined by positively-oriented meridians linking the $k$ components of $\Delta$ at $F_0$. If $F \in W - \Delta$, then the monodromy of $\mathcal{M}_0$ over $W - \Delta$ acts on $L(S)$, and it is well-known that this group is generated by the reflections in the vanishing cycles associated to the nodes of $F$. Since distinct nodes have orthogonal vanishing cycles, the image of $\pi_1(W - \Delta)$ is a group of the form $(\mathbb{Z}/2)^k$.

(3.6) According to the foregoing discussion, the restriction of the covering map to each component $V$ of the preimage of $W - \Delta$ is equivalent to the map

$$(y_1, \ldots, y_{20}) \mapsto (z_1, \ldots, z_{20}) = (y_1^2, \ldots, y_k^2, y_{k+1}, \ldots, y_{20}) \tag{3.6.1}$$

on the complement of the set $y_1 \cdots y_k = 0$ in some neighborhood of 0 in $\mathbb{C}^{20}$. The Fox completion adjoins to $V$ a copy of the set $y_1 \cdots y_k = 0$ and extends the projection map in the obvious way. This makes it clear that $\mathcal{M}_s$ is a topological manifold and that the branched covering map $\mathcal{M}_s \to \mathcal{C}_s$ is a local homeomorphism on the preimage of $\Delta_s^k$ for each $k$. It is also clear that $\mathcal{M}_s$ has a unique complex manifold structure compatible with these structures; this is given by coordinate charts like the $y_i$ on the left side of (3.6.1).

(3.7) Now we will discuss the structure of $\Delta_{ss}$ in $\mathcal{C}_{ss}$. A cubic surface has at most three cusps; a tricuspidal surface has no other singularities; a cubic surface with two cusps can have at most one node; and one with a single cusp can have at most two nodes. The space $\Delta_{ss}^{a,b}$ of cubic forms with exactly $a$ nodes and $b$ cusps, and no other singularities, has codimension $a + 2b$ in $\mathcal{C}$. Note that $\Delta_{ss}^{k,0}$ coincides with the space $\Delta_s^k$ introduced in (3.5). If $F \in \Delta_{ss}^{0,1}$, then there are coordinates centered at $F$ in which $\Delta$ has local equation $z_1^3 - z_2^2 = 0$. Thus, if $W$ is a suitable small neighborhood of $F$, then $\pi_1(W - \Delta) \cong B_3$, where $B_3$ denotes the classical braid group on three strands. If $F \in \Delta_{ss}^{a,b}$, then in a neighborhood of $F$ one can describe $\Delta$ as the transverse intersection of $a$ smooth hypersurfaces of the form $\Delta_s^1$ and $b$ smooth hypersurfaces containing a codimension one stratum of cusps of the form $\Delta_{ss}^{0,1}$. Thus, if $W$ is a small neighborhood of $F$, then $\pi_1(W - \Delta) \cong \mathbb{Z}^a \times B_3^b$. The facts about $\Delta_{ss}$ used here may be checked by referring to [10].

(3.8) The monodromy action of this group $\mathbb{Z}^a \times B_3^b$ on the restriction of $L(\mathcal{S}_0)$ to $W - \Delta$ is generated by the reflections in the roots in the spaces of vanishing cycles of the singularities of $F$. These form a root system of type $A_1^a A_2^b$, so the local monodromy group is $(\mathbb{Z}/2)^a \times S_3^b$. In order to describe a neighborhood of a point of $\mathcal{M}_{ss}$ lying over $F$, we first recall a fact about the natural action of $S_3$ on $A = \{(u_1, u_2, u_3) \in \mathbb{C}^3 : u_1 + u_2 + u_3 = 0\}$. The quotient $B = A/S_3$ is isomorphic (as an analytic space) to $\mathbb{C}^2$. Furthermore, the three involutions of $S_3$ act by reflections across the hyperplanes $u_1 = u_2$, $u_2 = u_3$ and $u_3 = u_1$, and $S_3$ acts freely away from the union $D$ of these mirrors. Finally, there are coordinates $(v_1, v_2)$ on $B$ such that the image of the mirrors under the natural map $j : A \to B$ is the variety $E$ defined by $v_1^3 - v_2^2 = 0$. From this it follows that $j : A \to B$ is the Fox completion of the covering map $(A - D) \to (B - E)$. We can now determine the local structure of the Fox completion $\mathcal{M}_{ss}$. Since $F \in \Delta_{ss}^{a,b}$, we may choose a neighborhood $W$ of $F$ that is isomorphic to a neighborhood of 0 in

$$\underbrace{\mathbb{C} \times \cdots \times \mathbb{C}}_{a \text{ factors}} \times \underbrace{B \times \cdots \times B}_{b \text{ factors}} \times \mathbb{C}^{20-a-2b}$$

in such a way that $W - \Delta$ corresponds to the intersection of a neighborhood of 0 with

$$\bigl(\mathbb{C} - \{0\}\bigr) \times \cdots \times \bigl(\mathbb{C} - \{0\}\bigr) \times (B - E) \times \cdots \times (B - E) \times \mathbb{C}^{20-a-2b}\ .$$

Then a component $V$ of the preimage of $W - \Delta$ in $\mathcal{M}_0$ is a copy of the intersection of

$$\bigl(\mathbb{C} - \{0\}\bigr) \times \cdots \times \bigl(\mathbb{C} - \{0\}\bigr) \times (A - D) \times \cdots \times (A - D) \times \mathbb{C}^{20-a-2b}$$

with a neighborhood of the origin in $\mathbb{C}^{20}$. Moreover, the projection to $\mathcal{C}_0$ is given by squaring each of the first $a$ coordinates, applying $j$ to each factor $A$, and acting by the identity of $\mathbb{C}^{20-a-2b}$. Passing to the Fox completion of $V \to W - \Delta$ over $W$ enlarges $V$ to a neighborhood of the origin in

$$\mathbb{C} \times \cdots \times \mathbb{C} \times A \times \cdots \times A \times \mathbb{C}^{20-a-2b}$$

in the obvious way, and we see that each component of the preimage of $W$ in $\mathcal{M}_{ss}$ is homeomorphic to a neighborhood of 0 in $\mathbb{C}^{20}$. As in the stable case, there is a unique complex manifold structure on $\mathcal{M}_{ss}$ compatible with the various structures already present. Finally, the preimage of $\Delta_{ss}$ in $\mathcal{M}_{ss}$ may be described in local analytic coordinates as the intersection of a suitable neighborhood of the origin with the hyperplanes associated to the root system $A_1^a \oplus A_2^b$.

(3.9) In a similar way we define a new structure called a *framing*. Recall from (2.2), (2.7) and (2.9) the definitions of the $\mathcal{E}$-lattices $\Lambda(T)$ and $\Lambda$ and the sheaf $\Lambda(\mathcal{T}_0)$ over $\mathcal{C}_0$. By a framing of a smooth form $F \in \mathcal{C}_0$ we mean a projective equivalence class $Pf$, where $f : \Lambda(T) \to \Lambda$ is an isometry; to form a projective equivalence class we identify isometries that differ by multiplication by a unit of $\mathcal{E}$. Now define the space $\mathcal{F}_0$ of framed smooth cubic forms as follows. As a set, $\mathcal{F}_0$ is the collection of all framings of all smooth cubic forms. Since the stalk of $\Lambda(\mathcal{T}_0)$ at $F \in \mathcal{C}_0$ is canononically isomorphic to $\Lambda(T)$, the set $\mathcal{F}_0$ is in natural one-to-one correspondence with the subsheaf of $PHom(\Lambda(\mathcal{T}_0), \mathcal{C}_0 \times \Lambda)$ consisting of projective equivalence classes of homomorphisms that are isometries on each stalk. We use this one-to-one correspondence to view $\mathcal{F}_0$ a complex manifold. In computations we will often work with isometries rather than equivalence classes and leave it to the reader to check that the results depend only on the equivalence class. We may also denote $Pf$ simply by $f$. The action of $g \in GL(4, \mathbb{C})$ on $\mathcal{C}_0$ defined in (2.17) lifts to an action of $GL(4, \mathbb{C})$ on $\mathcal{T}_0$, and this induces an action of $G$ on $\mathcal{F}_0$. Explicitly, if $g \in GL(4, \mathbb{C})$ and $Pf$ is a framing of $F \in \mathcal{C}_0$, then $P(f \circ g^*)$ is a framing of $gF$ that depends only on the image of $g$ in $G$. Here we have used the fact that an element of $D$ acts on $\mathcal{T}_0$ by a power of the branched covering transformation $\sigma$ which induces the scalar $\omega$ on $\Lambda(\mathcal{T}_0)$. Similarly, if $w \in P\operatorname{Aut}(\Lambda)$ then $w \circ Pf$ is a framing of $F$. Thus commuting actions of $G$ and $P\operatorname{Aut}(\Lambda)$ on $\mathcal{F}_0$ are defined. Moreover, it is clear that $\mathcal{F}_0$ is a covering space of $\mathcal{C}_0$ with Galois group $P\operatorname{Aut}(\Lambda)$. Indeed, since $P\Gamma = P\operatorname{Aut}\Lambda$ by Theorem 2.14, this space coincides with the covering $\tilde{\mathcal{C}}_0$ introduced in (2.10), and so $\mathcal{F}_0$ is connected.

(3.10) We define the space $\mathcal{F}_s$ of framed stable forms to be the Fox completion of $\mathcal{F}_0$ over $\mathcal{C}_s$. We call a preimage of $F \in \mathcal{C}_s$ in $\mathcal{F}_s$ a framing of $F$, and we call an element of $\mathcal{F}_s$ a framed stable form. As before, the Fox completion can be understood in terms of the local structure of $\Delta_s$ and the local monodromy groups. Imitating the case of marked surfaces, we take for each $F_0$ in $\Delta_s^k$ a neighborhood $W$ of $F_0$ in $\mathcal{C}_0$ such that $\pi_1(W - \Delta) \cong \mathbb{Z}^k$. By (5.4) and (5.7), the image of the local fundamental group in $P\Gamma$ is a group $(\mathbb{Z}/6)^k$, so that the restriction of the covering map $\mathcal{F}_0 \to \mathcal{C}_0$ is just like (3.6.1) but with $y_1^2, \ldots, y_k^2$ replaced by $y_1^6, \ldots, y_k^6$. Passing to the Fox completion adjoins the set $y_1 \cdots y_k = 0$ in the same way as before, and $\mathcal{F}_s$ is a complex manifold in a natural way. In particular, $\mathcal{F}_s$ is a cover of $\mathcal{C}_s$ with sixfold branching over each component of the preimage of $\Delta_s^1$, whereas $\mathcal{M}_s$ has only twofold branching there. Since the Fox completion of $\mathcal{F}_0$ over $\mathcal{C}_{ss}$ is not an analytic space and not even locally compact, there is no reasonable notion of a framed semistable form.

(3.11) The introduction of $\mathcal{F}_0$ allows us to give an explicit expression for the period map $\tilde{\mathcal{C}}_0 \to \mathbb{C}H^4$ introduced in (2.16). To this end, observe that the choice of $F_0$ and $\phi$ in (2.10) and (2.11) defines a basepoint for $\mathcal{F}_0$ and hence defines an isomorphism between $\tilde{\mathcal{C}}_0$ and $\mathcal{F}_0$. If $Pf \in \mathcal{F}_0$ lies over $F \in \mathcal{C}_0$ and is represented by an isometry $f : \Lambda(T) \longrightarrow \Lambda$, then $(f \otimes 1) \circ Z^{-1} : H^3_\omega(T) \longrightarrow \Lambda \otimes_\mathcal{E} \mathbb{C} = \mathbb{C}^{4,1}$ is an isometry. Here $Z$ is as in (2.5.1). Let

$$f_* : P(H^3_\omega) \xrightarrow{\cong} P(\mathbb{C}^{4,1})$$

be the resulting identification of projective spaces. Then the value of the period map at $Pf$ is

given by the formula
$$g(Pf) = f_*(H_{\bar{\omega}}^{2,1}) \in \mathbb{C}H^4 \subset P(\mathbb{C}^{4,1}). \tag{3.11.1}$$

We will give an even more explicit formula in terms of integrals in (6.7).

(3.12) Next we explain the relation between $\mathcal{F}_s$ and $\mathcal{M}_s$. Recall the definition in (2.12) of the $\mathbb{F}_3$-vector space $V$ and its quadratic form, and let

$$\Gamma_\theta = \ker\bigl(\mathrm{Aut}(\Lambda) \longrightarrow \mathrm{Aut}(V)\bigr) = \ker\bigl(\Gamma \longrightarrow \mathrm{Aut}^+(V)\bigr).$$

The subscript $\theta$ reflects the fact that $\Gamma_\theta$ is the level $\theta$ congruence subgroup of $\Gamma$. By (4.5) and (2.14), $\Gamma$ maps surjectively to $\mathrm{Aut}^+(V)$, so there is an exact sequence

$$1 \longrightarrow \Gamma_\theta \longrightarrow \Gamma \longrightarrow \mathrm{Aut}^+(V) \longrightarrow 1. \tag{3.12.1}$$

The center of $\mathrm{Aut}^+(V)$ is trivial because the central involution of $V$ has spinor norm $-1$. The center of $\Gamma = \mathrm{Aut}^+(\Lambda)$ consists of the cube roots of unity acting as scalars, all of which lie in $\Gamma_\theta$. Furthermore, $P\mathrm{Aut}^+(V) = P\mathrm{Aut}(V)$. Therefore by projectivizing (3.12.1) we obtain another exact sequence

$$1 \longrightarrow P\Gamma_\theta \longrightarrow P\Gamma \longrightarrow P\mathrm{Aut}(V) \longrightarrow 1.$$

According to [17], p. 26, $P\mathrm{Aut}(V) = PGO_5(3)$ is isomorphic to the Weyl group $W(E_6)$. Therefore $P\Gamma_\theta \backslash \mathcal{F}_0$ and $\mathcal{M}_0$ are covering spaces of $\mathcal{C}_0$ with isomorphic Galois groups. This suggests the next result, whose proof appears in (4.8)–(4.9) (but see also (4.10)).

**(3.13) Lemma.** *The spaces $\Gamma_\theta \backslash \mathcal{F}_0$ and $\mathcal{M}_0$ are isomorphic as covering spaces of $\mathcal{C}_0$. The spaces $\Gamma_\theta \backslash \mathcal{F}_s$ and $\mathcal{M}_s$ are isomorphic as branched covering spaces of $\mathcal{C}_s$. Both of these isomorphisms are $G$-equivariant.*

**(3.14) Lemma.** *$G$ acts freely and properly on $\mathcal{M}_s$, $\mathcal{F}_s$, and, in particular, on $\tilde{\mathcal{C}}_0 = \mathcal{F}_0 \subset \mathcal{F}_s$.*

The proof appears in (4.11)–(4.14).

(3.15) Now we gather together all of the moduli spaces considered in this paper:

$$\begin{aligned}
M_0^f &= G\backslash \mathcal{F}_0, & M_s^f &= G\backslash \mathcal{F}_s, & & \\
M_0^m &= G\backslash \mathcal{M}_0, & M_s^m &= G\backslash \mathcal{M}_s, & M^m &= PG\backslash\backslash P\mathcal{M}_{ss}, \\
M_0 &= G\backslash \mathcal{C}_0, & M_s &= G\backslash \mathcal{C}_s, & M &= PG\backslash\backslash P\mathcal{C}_{ss}.
\end{aligned} \tag{3.15.1}$$

These are the moduli spaces of framed smooth surfaces, framed stable surfaces, marked smooth surfaces, etc. We introduced $M_0$ and $M_0^f$ in (2.18). By Lemma 3.14, $G$ acts freely on $\mathcal{M}_0$, $\mathcal{M}_s$, $\mathcal{F}_0$ and $\mathcal{F}_s$, so that $M_0^m$, $M_s^m$, $M_0^f$ and $M_s^f$ are complex manifolds. From geometric invariant theory, we know that $G$ acts properly on $\mathcal{C}_s$. Therefore $M_s$ and its subspace $M_0$ are analytic spaces and even complex analytic orbifolds. The space $M = PG\backslash\backslash P\mathcal{C}_{ss}$ is defined to be $SL(4,\mathbb{C})\backslash\backslash P\mathcal{C}$, the geometric invariant theory quotient of $P\mathcal{C}_{ss}$ by $SL(4,\mathbb{C})$. By definition this is $P\mathcal{C}_{ss}$ modulo the equivalence relation closure (in $P\mathcal{C}_{ss}$) of $PSL(4,\mathbb{C})$-orbits, or equivalently of $PG$-orbits. GIT implies that the quotient by this equivalence relation is a projective variety whose graded coordinate ring is the ring of $G$-invariant polynomials on $\mathcal{C}$. GIT also implies that points of $M$ correspond to the closed orbits of $PG$ in $P\mathcal{C}_{ss}$. By Lemma 4.6, $M$ is the union of $M_s$ with the equivalence class of $X_0^3 - X_1X_2X_3$. The space $M^m$ is the following GIT quotient. Let $p : P\mathcal{M} \longrightarrow P\mathcal{C}$ be the natural projection. It is easy to see that $p^*\mathcal{O}(1)$ is an ample line bundle on $P\mathcal{M}$ and that there are commuting actions of $W(E_6)$ and $SL(4,\mathbb{C})$ on $p^*\mathcal{O}(1)$. It is also easy to see that $m \in \mathcal{M}$ is stable (respectively semistable) with respect to the action of $SL(4,\mathbb{C})$ on $\mathcal{M}$ and the linearized line bundle $p^*\mathcal{O}(1)$ if and only if $F = p(m) \in \mathcal{C}$ is stable (respectively semi-stable) with respect to the action of $SL(4,\mathbb{C})$ on $\mathcal{C}$ and the linearized line bundle $\mathcal{O}(1)$. Thus the previously defined spaces $\mathcal{M}_s$ and $\mathcal{M}_{ss}$ coincide with the spaces of stable and semi-stable elements of $\mathcal{M}$ for this particular linearization of the $SL(4,\mathbb{C})$-action on $\mathcal{M}$. We therefore define $M^m = PG \backslash\backslash P\mathcal{M}_{ss} = SL(4,\mathbb{C}) \backslash\backslash P\mathcal{M}_{ss}$ using this linearization of the $SL(4,\mathbb{C})$-action. Because the actions of $P\Gamma$ and $G$ commute, the spaces of (3.15.1) fit into the commutative diagram

$$\begin{array}{ccccc} M_0^f & \longrightarrow & M_s^f & & \\ \downarrow & & \downarrow & & \\ M_0^m & \longrightarrow & M_s^m & \longrightarrow & M^m \\ \downarrow & & \downarrow & & \downarrow \\ M_0 & \longrightarrow & M_s & \longrightarrow & M, \end{array}$$

where the horizontal maps are inclusions and the vertical maps are coverings or branched coverings. There is no space $M^f$ because, as explained in (3.10), there is no reasonable space $\mathcal{F}_{ss}$.

(3.16) Recall from (2.16) the definition of the period map $g : \mathcal{F}_0 \longrightarrow \mathbb{C}H^4$. Using the Riemann extension theorem and the fact that $\mathbb{C}H^4$ is a bounded domain in $\mathbb{C}^4$, we see that $g$ extends to a map $\mathcal{F}_s \longrightarrow \mathbb{C}H^4$. Since $g$ is constant on $G$-orbits in $\mathcal{F}_0$, so is the extension. Consequently there is a holomorphic quotient map $M_s^f \longrightarrow \mathbb{C}H^4$ which we denote by the same symbol $g$ and refer to as a period map. Dividing by the $\Gamma_\theta$ and $\Gamma$-actions and using Lemma 3.13, we obtain period maps corresponding to the horizontal arrows of the following commutative diagram.

$$\begin{array}{ccc} M_s^f & \longrightarrow & \mathbb{C}H^4 \\ \downarrow & & \downarrow \\ M_s^m & \longrightarrow & P\Gamma_\theta \backslash \mathbb{C}H^4 \\ \downarrow & & \downarrow \\ M_s & \longrightarrow & P\Gamma \backslash \mathbb{C}H^4 \end{array} \quad (3.16.1)$$

(The vertical maps are branched covers.) We can now state our main theorem. In it we write $\mathcal{H}^k$ for the set of points of $\mathbb{C}H^4$ that lie on exactly $k$ components of the hyperplane arrangement $\mathcal{H}$ introduced in (2.19). In particular, $\mathcal{H}^0 = \mathbb{C}H^4 - \mathcal{H}$.

**(3.17) Theorem.** *The top map of (3.16.1) is an isomorphism of complex manifolds, and the other horizontal maps are isomorphisms of analytic spaces. Furthermore, for each $k = 0, \ldots, 4$, the period map carries the framed cubic forms with $k$ nodes and no other singularities onto $\mathcal{H}^k$. Finally, the maps of the bottom two rows of (3.16.1) extend to isomorphisms of compact analytic spaces as in the diagram below.*

$$\begin{array}{ccc} M^m & \longrightarrow & \overline{P\Gamma_\theta \backslash \mathbb{C}H^4} \\ \downarrow & & \downarrow \\ M & \longrightarrow & \overline{P\Gamma \backslash \mathbb{C}H^4} \end{array} \quad (3.17.1)$$

*The bars indicate the Satake compactifications of the two ball quotients.*

(3.18) We have shown that $M_s^f$, with its natural complex manifold structure, is isomorphic to $\mathbb{C}H^4$. As mentioned in (3.15), $M_s^m$ is a complex manifold and $M_s$ is a complex orbifold. Theorem 3.17 shows that $M_s^m$ and $M_s$ are isomorphic as analytic spaces to the ball quotients $P\Gamma_\theta\backslash\mathbb{C}H^4$ and $P\Gamma\backslash\mathbb{C}H^4$. As quotients of $\mathbb{C}H^4$ by discrete groups, the latter spaces are also orbifolds. However, the analytic space isomorphisms are *not* orbifold isomorphisms. To see this for $M_s^m$, one need only observe that (as an orbifold) $P\Gamma_\theta\backslash\mathbb{C}H^4$ is not a manifold. Indeed, $P\Gamma_\theta$ contains elements with fixed points in $\mathbb{C}H^4$ (see for example (7.25)). To see the non-isomorphism for $M_s$, one uses the fact that there exists a nodal cubic surface $S$ with no symmetries. Since $S$ has no symmetries, the corresponding point of $M_s$ is a manifold point. But because $S$ has a node, the corresponding point of $P\Gamma\backslash\mathbb{C}H^4$ lies in the image of a hyperplane of $\mathcal{H}$. Such a point is not a manifold point of the orbifold $P\Gamma\backslash\mathbb{C}H^4$ because $P\Gamma$ contains the hexflections which fix the hyperplane pointwise.

(3.19) Despite the foregoing discussion, we can define orbifold stuctures on the analytic spaces $P\Gamma_\theta\backslash\mathbb{C}H^4$ and $P\Gamma\backslash\mathbb{C}H^4$ so that the two lower horizontal maps in (3.16.1) are orbifold isomorphisms. The idea, say for $P\Gamma\backslash\mathbb{C}H^4$, is to start with the complex hyperbolic orbifold $P\Gamma\backslash\mathbb{C}H^4$ and replace the generic points of the image of $\mathcal{H}$, which are orbifold points with local group $\mathbb{Z}/6$, by ordinary manifold points. To make this precise, and to deal with the non-generic points of the image of $\mathcal{H}$, consider $x \in P\Gamma\backslash\mathbb{C}H^4$ and let $y$ be a point in $\mathbb{C}H^4$ lying over the point of $P\Gamma\backslash\mathbb{C}H^4$ that corresponds to $x$. Write $H$ for the stabilizer of $y$ in $P\Gamma$ and $N$ for the normal subgroup of $H$ generated by the reflections in the short roots of $\Lambda$ whose mirrors pass through $y$. Choose a small $H$-invariant neighborhood $U$ of $y$ and observe that the normality of $N$ implies that the map $U \to U/H$ factors as the composition

$$U \to U/N \to (U/N)/(H/N) = U/H \ .$$

Observe that $U/N$ may be viewed as a complex manifold, since $N$ is generated by reflections in orthogonal hyperplanes. Ignore the orbifold structure on $U/N$, regarding it instead as neighborhood $V$ of $0$ in $\mathbb{C}^4$. Then take the natural map from $V$ to $V/(H/N)$ to a neighborhood of $x$ in $P\Gamma\backslash\mathbb{C}H^4$ as an orbifold chart for $P\Gamma\backslash\mathbb{C}H^4$. It is easy to see that these charts are compatible with each other. Consequently they equip $P\Gamma\backslash\mathbb{C}H^4$ with a complex analytic orbifold structure. For $P\Gamma_\theta\backslash\mathbb{C}H^4$ we apply exactly the same construction. The only difference is that $N$ will be a product of $(\mathbb{Z}/3)$'s rather than $(\mathbb{Z}/6)$'s. It is even true that under this new orbifold structure, $P\Gamma_\theta\backslash\mathbb{C}H^4$ is a manifold; this follows from Lemma 7.27. Finally, we note that away from the image of $\mathcal{H}$, the new orbifold structure coincides with the old one.

**(3.20) Theorem.** *Let $P\Gamma\backslash\mathbb{C}H^4$ and $P\Gamma_\theta\backslash\mathbb{C}H^4$ be given the orbifold structures just defined, and let $M_s$ and $M_s^m$ be given their natural orbifold structures. Then the two lower maps of (3.16.1) are orbifold isomorphisms.*

(3.21) We note that although $M_s$ is an orbifold and not a manifold, its orbifold fundamental group $\pi_1(M_s)$ is trivial. The description of $M_s$ in terms of $P\Gamma\backslash\mathbb{C}H^4$ shows that $\pi_1(M_s)$ equals the orbifold fundamental group $\pi_1(P\Gamma\backslash\mathbb{C}H^4) = P\Gamma$, modulo the relations "reflections in the short roots of $\Lambda$ are trivial." By (7.1) and (7.21), these reflections generate $P\Gamma$. Therefore $\pi_1(M_s)$ is trivial. In the same way, it follows from (7.25) that the orbifold (indeed manifold) $M_s^m$ is also simply connected.

(3.22) Here now is the strategy for the proofs of Theorems 3.17 and 3.20 to be given in sections 4–9. First, we have to establish the various claims made in this section and the previous one that are required for our constructions to make sense. The material of section 2 up to (2.14) is verified in order in (4.1)–(4.5). This allows one to establish Lemma 5.4, which gives the local monodromy of $\Lambda(\mathcal{T}_0)$ near a nodal surface. Then all of section 7 follows, justifying (2.14). The remaining claims of section 2 are subsumed in more precise statements in section 3. The foundational results of section 3 needed to verify that these statements make sense are proved in (4.6)–(4.14), leaving only the proofs of Theorems 3.17 and 3.20. These occupy section 9, which is divided into four subsections. The first subsection proves that $M_0 \to P\Gamma\backslash(\mathbb{C}H^4 - \mathcal{H})$ is an isomorphism, and consists of three steps: (1) a local computation, Lemma 9.1, which relies on some material on fractional differentials developed in section 6; (2) the Torelli theorem of Clemens and Griffiths for cubic threefolds; (3) the existence of the extension of $g$ to $M \to \overline{P\Gamma\backslash\mathbb{C}H^4}$, which is obtained in section 8 by studying the period map near cuspidal surfaces. The second subsection proves that $g : (M_s - M_0) \longrightarrow P\Gamma\backslash\mathcal{H}$ is a homeomorphism by appealing to a theorem of Deligne and Mostow and some facts proved in sections 5–7. The third subsection gives the remaining simple details to complete the proof of Theorem 3.17. The fourth subsection proves Theorem 3.20.

## 4. Proofs of lemmas

We now give the the statements and/or proofs of various lemmas used in the two preceding sections.

**(4.1) Lemma.** *Let $M$ be a free $\mathcal{E}$-module whose underlying $\mathbb{Z}$-module is equipped with a unimodular symplectic form $\Omega$. Then the function*

$$h(x, y) = -\frac{\Omega(\theta x, y) + \theta \Omega(x, y)}{2} \qquad (4.1.1)$$

*is an $\mathcal{E}$-valued unimodular Hermitian form on $M$ (linear in its first argument and antilinear in its second).*

**Proof:** We first we show that $h(x, y) \in \mathcal{E}$ for all $x, y \in M$. Because $\Omega$ is $\mathbb{Z}$-valued, $h(x, y) = a/2 + \theta b/2$ for some $a, b \in \mathbb{Z}$. Since $\mathcal{E}$ consists of the numbers of this form with $a \equiv b \pmod 2$, it suffices to show that $\Omega(\theta x, y) \equiv \Omega(x, y) \pmod 2$. This amounts to showing that $\Omega(\theta x - x, y)$ is even, which is true because $\theta - 1 = 2\omega$.

(4.2) Next we show that $h(x, y) = \overline{h(y, x)}$. To do this, consider $\Omega$ to be defined by $\mathbb{R}$-linear extension on all of $M \otimes \mathbb{R}$, and $h$ to be defined on all of $M \otimes \mathbb{R}$ by the formula (4.1.1). The imaginary part of $h(x, y)$ is $-\theta\, \Omega(x, y)/2$, which obviously changes sign under exchange of $x$ and $y$. The real part is $-\Omega(\theta x, y)/2$, so it suffices to show $\Omega(\theta x, y) = \Omega(\theta y, x)$:

$$\Omega(\theta x, y) = \Omega(\sigma x - \sigma^{-1} x, y) = \Omega(\sigma x, y) - \Omega(\sigma^{-1} x, y)$$
$$= \Omega(x, \sigma^{-1} y) - \Omega(x, \sigma y) = \Omega(x, \sigma^{-1} y - \sigma y)$$
$$= \Omega(x, -\theta y) = \Omega(\theta y, x).$$

Now we prove $\mathbb{C}$-linearity in the first variable. Since $\mathbb{R}$-linearity follows from that of $\Omega$, it suffices to prove that $h(\theta x, y) = \theta h(x, y)$:

$$h(\theta x, y) = -\frac{1}{2}\Big(\Omega(-3x, y) + \theta \Omega(\theta x, y)\Big) = -\frac{1}{2}\Big(-3\Omega(x, y) + \theta \Omega(\theta x, y)\Big)$$
$$= -\frac{\theta}{2}\Big(\theta \Omega(x, y) + \Omega(\theta x, y)\Big) = \theta h(x, y).$$

Finally, we show that $h$ is unimodular. Suppose $x \in M \otimes \mathbb{R}$ satisfies $h(x, y) \in \mathcal{E}$ for all $y \in M$; we must prove $x \in M$. For each $y \in M$, $h(x, y) \in \mathcal{E}$ implies that the imaginary part of $h(x, y)$ is an integral multiple of $\theta/2$, so that $\Omega(x, y) \in \mathbb{Z}$. Since $\Omega$ is unimodular, $x \in M$.

**Proof of Lemma 2.5**

(4.3) Let $\xi \in H^3(T, \mathbb{Z})$ and introduce the temporary notation $\Theta = \sigma - \sigma^{-1}$, viewed as an operator on cohomology. On $H^3_\omega(T, \mathbb{C})$, $\sigma$ acts as scalar multiplication by $\omega$, and on $H^3_{\bar\omega}(T, \mathbb{C})$ it acts by $\bar\omega$. It follows that $\Theta(\xi_\omega) = \theta \xi_\omega$ and $\Theta(\xi_{\bar\omega}) = -\theta \xi_{\bar\omega}$. By definition, $h(\xi, \xi) = -\frac{1}{2}\Omega(\Theta \xi, \xi) - \frac{\theta}{2}\Omega(\xi, \xi)$. Using the antisymmetry of $\Omega$ and the properties of $\Theta$ given above, we find

$$h(\xi, \xi) = -\Omega(\theta \xi_\omega - \theta \xi_{\bar\omega}, \xi_\omega + \xi_{\bar\omega})/2$$
$$= -(\theta/2)[\Omega(\xi_\omega, \xi_{\bar\omega}) - \Omega(\xi_{\bar\omega}, \xi_w)]$$
$$= -\theta \Omega(\xi_\omega, \xi_{\bar\omega}).$$

Next we compute

$$h'(Z(\xi), Z(\xi)) = \theta \int_T \xi_{\bar\omega} \wedge \overline{\xi_{\bar\omega}} = \theta \int_T \xi_{\bar\omega} \wedge \xi_\omega = \theta \Omega(\xi_{\bar\omega}, \xi_\omega) = -\theta \Omega(\xi_\omega, \xi_{\bar\omega}).$$

This completes the proof.

**Proof of Lemma 2.6**

(4.4) Suppose $F \in \mathcal{C}_0$. Since $\sigma$ is holomorphic, the Hodge decomposition refines the eigenspace decomposition. It is well-known that $H^{3,0}(T)$ and $H^{0,3}(T)$ vanish. Because $\sigma$ is defined on $H^3(T, \mathbb{R})$, $H^3_\omega$ and $H^3_{\bar\omega}$ have the same dimension, so each has dimension 5. In §5 of [13] it is proved that $H^{2,1}_{\bar\omega}$ is one dimensional. In Theorem 6.5 we will prove this fact again using suitable fractional differentials. From this the Hodge numbers of $H^3_{\bar\omega}$ are determined:

$$\dim H^{2,1}_{\bar\omega} = 1, \qquad \dim H^{1,2}_{\bar\omega} = 4.$$

From the Hodge numbers we obtain the signature of $h'$. Let $\phi$ be a nonzero vector in $H^{2,1}_{\bar\omega}(T)$. Then the inner product $h'(\phi, \phi)$ is a sum of integrals which in local coordinates have the form

$$\theta \int_T |f|^2 dx \wedge dy \wedge d\bar{z} \wedge d\bar{x} \wedge d\bar{y} \wedge dz = \theta \int_T |f|^2 dx \wedge d\bar{x} \wedge dy \wedge d\bar{y} \wedge dz \wedge d\bar{z}.$$

The last expression is negative, so $h'$ is negative-definite on $H^{2,1}_{\bar\omega}$; a similar computation shows that it is positive-definite on $H^{1,2}_{\bar\omega}$. Consequently the signature of $h'$ on $H^3_{\bar\omega}$ is $(4, 1)$.

**(4.5) Lemma.** *Both maps of (2.12.1) are surjective.*

**Proof:** It is trivial to check that $\nu$ takes both possible values. Since $\mathrm{Aut}(V)$ is generated by reflections in its elements of nonzero norm, it suffices to show that such reflections are induced by automorphisms of $\Lambda$. One verifies directly that each nonisotropic vector of $V$ is the image of a lattice vector of norm $\pm 1$ or $\pm 2$. Biflections (order 2 reflections) in such vectors preserve $\Lambda$ and reduce to the desired automorphisms of $V$.

**(4.6) Lemma.** *The cubic form $F = X_0^3 - X_1 X_2 X_3$ represents the unique closed $SL(4, \mathbb{C})$ orbit in $P(\mathcal{C}_{ss} - \mathcal{C}_s)$.*

(4.7) **Proof:** By p. 51 of [39] or §19 of [26], $F$ is semistable, since $S$ has three cusps and no other singularities. Now suppose $F' \in \mathcal{C}_{ss} - \mathcal{C}_s$; we will show that the closure $K$ of its orbit contains $F$. Since $F' \notin \mathcal{C}_s$, the surface $S'$ has an $A_2$ singularity at a point $p$. We may choose affine coordinates $x_0, x_2, x_3$ centered at $p$, such that $S$ has local equation $Q(x_0, x_2, x_3) + C(x_0, x_2, x_3)$, where $Q$ (resp. $C$) is a quadratic (resp. cubic) form. Since $p$ is an $A_2$ singularity, $Q$ is a quadratic form of rank 2. Therefore, by a linear change of coordinates we may suppose that $Q = -x_2 x_3$. Passing to homogeneous coordinates, we find that $F' = -X_1 X_2 X_3 + C(X_0, X_2, X_3)$. We let $t$ approach $\infty$ in the 1-parameter group $\mathrm{diag}[1, t^2, t^{-1}, t^{-1}]$. In the limit all terms of $C$ vanish except the $X_0^3$ term, say $bX_0^3$. Therefore $K$ contains the form $bX_0^3 - X_1 X_2 X_3$. If $b = 0$ then a further degeneration shows that $0 \in K$, which is impossible. Therefore $b \neq 0$. By rescaling $X_0$ we may take $b = 1$, which proves $F \in K$.

**Proof of Lemma 3.13**

(4.8) Choose a basepoint $F$ of $\mathcal{C}_0$. Our goal is to show that the kernel of the monodromy representation of $\pi_1(\mathcal{C}_0, F)$ on $L(S)$ coincides with the kernel of the monodromy representation on the five-dimensional vector space $V(T) = \Lambda(T)/\theta\Lambda(T)$ over $\mathbb{F}_3$. The first step is to build a 5-dimensional $\mathbb{F}_3$ vector space from $L(S)$; we shall call it $V(S)$. To this end let $L_0(S)$ be the primitive cohomology $\eta(S)^\perp$ in $L(S)$; it is well-known to be a copy of the root lattice $-E_6$ of determinant 3. Write $L_0'(S)$ for the lattice dual to $L_0(S)$ and recall from [8] VI §4, exercise 2, that $L_0(S)/3L_0'(S)$ is a nondegenerate quadratic space of dimension five over $\mathbb{F}_3$. This is our space $V(S)$. Inner products are computed by reducing the inner products of lattice vectors modulo 3. One observes that $\mathrm{Aut}(L(S), \eta(S))$ acts faithfully on both $L_0(S)$ and $V(S)$.

(4.9) Now we show that the monodromy representations of $\pi_1(\mathcal{C}_0, F)$ on $V(S)$ and $V(T)$ are equivalent. By the discussion of (7.1), the fundamental group is the image of the Artin group $\mathcal{A}(E_6)$ with the six standard generators mapping to meridians $\mu_1, \ldots, \mu_6$ of $\Delta$. The monodromy action on $L(T)$ is well-known: each $\mu_i$ maps to the reflection in some root $r_i$ (a norm $-2$ vector) of $L_0$. Furthermore, the $r_i$ may be taken to form a simple root system of type $-E_6$, so that $r_i$ and $r_j$ have inner product 1 or 0 according to whether the corresponding vertices of the Dynkin diagram a joined or not. The images $\bar{r}_i$ of the $r_i$ in $V(S)$ span $V(S)$ and have the same mutual inner products as the $r_i$, and all of them have norm 1. Now, by (7.5) the $\mu_i$ act on $\Lambda(T)$ by hexflections in lattice vectors $r_i'$ of norm 1, with pairwise inner products given by the same rule.

Therefore the $\mu_i$ act on $V(S)$ by the biflections in the images $\bar{r}'_i$ of the $r'_i$. Since the $\bar{r}'_i$ have the same inner product matrix as the $\bar{r}_i$, the map $\bar{r}_i \mapsto \bar{r}'_i$ defines an isometry from $V(S)$ to $V(T)$. The monodromy representations on $V(S)$ and $V(T)$ are now visibly equivalent.

**Alternate proof (sketch) for Lemma 3.13**

(4.10) A longer but more satisfying proof can be given by constructing a natural isometry from $V(S)$ to $V(T)$ that identifies the monodromy representations of $\pi_1(\mathcal{C}_0, F)$ on them. We begin by building a map $k : L_0(S) \to V(T)$. If $\alpha \in L_0(S)$, consider its Poincaré dual $\hat{\alpha} \in H_2(S, \mathbb{Z})$. Since $\hat{\alpha}$ is a primitive homology class of $S \subset T$ and $T$ has no primitive second homology, $\hat{\alpha}$ bounds a 3-chain $\hat{\beta}$ in $T$. Since $\hat{\beta}$, $\sigma_*^{-1}(\hat{\beta})$ and $\sigma_*(\hat{\beta})$ all have the same boundary, $\hat{\gamma} = \sigma_*^{-1}\hat{\beta} - \sigma_*\hat{\beta}$ represents an element of $H_3(T, \mathbb{Z})$. We let $\gamma \in \Lambda(T) = H^3(T, \mathbb{Z})$ be the Poincaré dual of $\hat{\gamma}$, and define $k(\alpha)$ as the reduction of $\gamma$ modulo $\theta$. Despite the ambiguity in the choice of $\hat{\beta}$, this map is well-defined, because if we replaced $\hat{\beta}$ by some other chain $\hat{\beta}'$ then the corresponding $\gamma'$ would differ from $\gamma$ by an element of $(\sigma^* - (\sigma^*)^{-1})\Lambda(T) = \theta\Lambda(T)$. One can show by an explicit geometric construction and computation of intersection numbers that $k$ vanishes on $3L'_0(S)$ and that the induced map $\bar{k} : V(S) = L_0(S)/3L'_0(S) \to V(T)$ is an isometry. In fact, once one shows that $k$ is non-zero, simple equivariance and irreducibility arguments force $k$ to descend to an isomorphism $\bar{k} : V(S) = L_0(S)/3L'_0(S) \to V(T)$ of $\pi_1(\mathcal{C}_0, F)$-modules. A comparison of the determinants of $V(S)$ and $V(T)$ then shows that $\bar{k}$ must be an isometry.

**Proof of Lemma 3.14**

By Lemma 3.13 the claims for $\mathcal{F}_s$ follow from those for $\mathcal{M}_s$. Since the properness of the $G$-action on $\mathcal{M}_s$ follows from that of the action on $\mathcal{C}_s$, it suffices to show that $G$ acts freely on $\mathcal{M}_s$. Freeness follows from the next two lemmas.

(4.11) **Lemma.** *Let $m \in \mathcal{M}_s$ lie over $F \in \mathcal{C}_s$. Then $m$ determines:*

i. *A bijection $b_m$ between the set of nodes $\{p_1, \cdots, p_k\}$ (if any) of $S$ and a set $\{b_m(p_1), \cdots, b_m(p_k)\}$ of mutually orthogonal sublattices of $L$, each generated by a root $r_i$. We write $R_i = b_m(p_i)$ and $R = R_1 \oplus \cdots \oplus R_k$.*

ii. *An isometric embedding $i_m : (L(S), \eta(S)) \longrightarrow (L, \eta)$ with image the orthogonal complement of $R$.*

**Proof:** If $m \in \mathcal{M}_0$, then the first assertion is vacuous and the second assertion is the definition of a marked smooth form. If $m$ lies over $F \in \Delta_s^k$, $k > 0$, then $m$ may be regarded as a choice of component $C$ of the preimage in $\mathcal{M}_0$ of $W - \Delta$, where $W$ is as in (3.3). Picard-Lefschetz theory tells us first that $L(\mathcal{S})|_{W-\Delta}$ has mutually orthogonal rank-one local subsystems of vanishing cycles, $\mathcal{V}_1, \cdots, \mathcal{V}_k$, each in natural one-to-one correspondence with the nodes $p_1, \cdots, p_k$ of $S$, and second that there is a constant subsystem $\mathcal{I} = \mathcal{V}^\perp$ of rank $7 - k$ of invariant cohomology, where $\mathcal{V} = \mathcal{V}_1 \oplus \cdots \oplus \mathcal{V}_k$ is the vanishing cohomology. Moreover, there is a deformation retraction $r : \mathcal{S}|_W \longrightarrow S$ such that the constant system $\mathcal{I}$ is trivialized by the isomorphism $(W - \Delta) \times L(S) \longrightarrow \mathcal{I}|_{W-\Delta}$ defined by $(F_1, x) \longrightarrow i_{F_1}^* r^* x$, where $F_1 \in W - \Delta$ and $i_{F_1} : S' \longrightarrow \mathcal{S}$ is the inclusion of the fiber $S'$ over $F_1$. We choose $F_1 \in W - \Delta$, and a marking $m_1 : (L(S_1), \eta(S_1)) \to (L, \eta)$ of $F_1$ that lies in $C$. Then $C$ consists of the markings $m_2 : (L(S_2), \eta(S_2)) \to (L, \eta)$ of forms $F_2 \in W - \Delta$ that may be obtained by parallel transport of $m_1$ along paths in $W - \Delta$. Since parallel translation preserves every subsystem, it follows that $m_1^{-1} m_2(\mathcal{V}_i) = \mathcal{V}_i$ for $i = 1, \cdots, k$, and $m_1^{-1} m_2|_\mathcal{I} = id$. Therefore $m_1(\mathcal{V}_i) = m_2(\mathcal{V}_i)$ for each $i$. We denote this common value by $b_m(p_i)$ and observe that these give mutually orthogonal sublattices of $L$, each generated by the image of a vanishing cycle, hence generated by a root. This proves the first assertion of the lemma. The second assertion follows because $m_1$ and $m_2$ give the same isometric embedding of $L(S)$ in $L$. Thus we take $i_m$ to be be the common value of $m_1$ and $m_2$.

(4.12) **Remark.** The previous lemma suggests that a point $m \in \mathcal{M}_s$ lying over $F \in \Delta_s$ should be described by a pair $(b_m, i_m)$. This is almost the case. The only difficulty is that $b_m$ and $i_m$ cannot in general be prescribed independently. Since a complete discussion would take us too far afield, we mention only the fact that for $k = 1, 2, 3$, $i_m$ does determine $b_m$. For $k = 4$, this is not the case, and in addition $b_m$ must satisfy some relations imposed by $i_m$.

**(4.13) Lemma.** *Suppose $F \in \mathcal{C}_s$ and that $g \in G$ fixes $F$, acting trivially on $H^2(S, \mathbb{Z})$. If $F$ is singular, suppose in addition that $g$ fixes each node of $S$. Then $g$ is the identity.*

(4.14) **Proof:** It suffices to prove that $g$ induces the identity on $S$. To this end, we will use some easily verified facts about Poincaré duality on $S$ and the homology classes of the lines on $S$. First, since $S$ is a rational homology manifold, there is a Poincaré duality isomorphism $H^2(S, \mathbb{R}) \cong H_2(S, \mathbb{R})$. As in [49], this isomorphism is actually defined on the level of $\mathbb{Z}$-modules: it is the transport to $S$ via a resolution of singularities $\bar{S} \longrightarrow S$ of the standard Poincaré duality isomorphism on $\bar{S}$. More precisely, it is the composition $H^2(S) \longrightarrow H^2(\bar{S}) \longrightarrow H_2(\bar{S}) \longrightarrow H_2(S)$, where $\bar{S} \longrightarrow S$ is the blow-up of the nodes. This map is injective, and its image has index 1, 2, 4, 8 or 8 according to whether $S$ has 0, 1, 2, 3 or 4 nodes.

Let $l$ be a line on $S$, and let $[l]$ denote its homology class in $H_2(S, \mathbb{Z})$. Then there are three cases to consider:

1. The line $l$ does not contain a node of $S$. Then $l$ has self-intersection $-1$ in the usual sense: the Poincaré dual $\lambda$ of its homology class has square $-1$. It follows that the line is rigid in its homology class: it is the unique algebraic curve in $[l]$.

2. The line $l$ contains exactly one node of $S$. Then $l$ still has negative self-intersection, and is rigid in its homology class. More precisely, there is a cohomology class $\lambda$ with $\lambda^2 = -2$ that is Poincaré dual to $2[l]$. It is easy to verify that $l$ itself is not in the image of the Poincaré duality map $H^2(S, \mathbb{Z}) \to H_2(S, \mathbb{Z})$. See the related example on p. 181 of [49].

3. The line $l$ contains two nodes of $S$. Then $l$ has zero self-intersection. More precisely, there is a cohomology class $\lambda$ which is Poincaré dual to $2[l]$, and $\lambda^2 = 0$.

We now return to the proof that $g$ is the identity on $S$. The key point is that $g$ preserves all the lines of $S$: lines of type 1 or 2 are preserved because they are rigid in their homology class, and $g$ acts as the identity on $H^2(S, \mathbb{Z})$. Lines of type 3 are preserved because $g$ fixes each node.

To complete the proof, we examine cases according to singularities of $S$. First, if $S$ is smooth, and if $l_1, \cdots, l_6$ are six disjoint lines on $S$, then $S$ can be blown down to $\mathbb{P}^2$ with the $l_i$ mapping to points $p_i$ in general position. Since $S-(l_1\cup\cdots\cup l_6) \cong \mathbb{P}^2-\{p_1,\cdots,p_6\}$, $g$ induces an automorphism of $\mathbb{P}^2-\{p_1,\cdots,p_6\}$ which by Hartogs' theorem extends to an automorphism of $\mathbb{P}^2$. This extension fixes six points in general position, hence is the identity of $\mathbb{P}^2$, and so $g$ is the identity on $S$. If $S$ has one node $q$, then projection from $q$ defines a birational map of $S$ to $\mathbb{P}^2$ that blows $q$ up to a non-singular conic $C \subset \mathbb{P}^2$ and contracts the six lines $l_1, \cdots, l_6$ through $q$ to points $p_1, \cdots, p_6$ on $C$. The map $g$ of $\mathbb{P}^3$ preserving $S$ again induces a map of $\mathbb{P}^2 - \{p_1, \cdots, p_6\}$ preserving $C$. It extends as before to an automorphism of $\mathbb{P}^2$ preserving $C$ and fixing $p_1, \cdots, p_6$, hence is the identity. If $S$ has $k = 2, 3$ or $4$ nodes then we choose one, say $q$; projecting $S$ to $\mathbb{P}^2$ from $q$ as before blows $q$ up to a non-singular conic $C \subset \mathbb{P}^2$. Now there are only $7-k$ lines in $S$ through $q$. The map induced by $g$ extends to all of $\mathbb{P}^2$, preserves $C$, and fixes the corresponding $7-k$ points on $C$. Since $7-k \geq 3$, this map is the identity on $C$, hence on $\mathbb{P}^2$, hence $g$ is the identity on $S$.

(4.15) **Remark.** The only stable surface admitting a nontrivial automorphism acting trivially on cohomology is the 4-nodal surface. The nontrivial automorphisms which act trivially there are just those that permute the nodes in two orbits of size two.

## 5. Topology of nodal degenerations

(5.1) In order to extend the period map to the various spaces considered in (3.15.1) and to prove our Theorem 3.17, we must understand the behavior of the period map near cubic forms $F_0 \in \Delta_{ss}$. To this end we must understand the topology of the family $\mathcal{T} \xrightarrow{\pi} \mathcal{C}$ of (2.8) near $F_0$, as well as the local structure of the sheaf $\Lambda(\mathcal{T})$. In this section we will study these problems for $F_0 \in \Delta_s$ (the nodal case). The cuspidal case will be treated in §8.

(5.2) Consider first a surface with a single node defined by $F_0 \in \Delta_s^1$. Then $F_0$ has a fundamental system of neighborhoods $U$ in $\mathcal{C}$ of the form $U = U_0 \times D$ where $U_0$, biholomorphic to a ball, is a neighborhood of $F_0$ in $\Delta_s^1$ and $D$ is biholomorphic to the disk $\{z : |z| < 1\} \subset \mathbb{C}$. The families $\mathcal{S}|_U$ and $\mathcal{T}|_U$ are topologically equivalent as families over $U$ to $U_0 \times (\mathcal{S}|_D)$ and $U_0 \times (\mathcal{T}|_D)$, respectively. See §I1.11 of [23]. Because of this it will be enough to study families parametrized by a disk transverse to $\Delta_s^1$ at $F_0$. Let $p \in \mathbb{P}^3$ be the unique singular point of $F_0$. Then there is a closed ball $B$ centered at $p$ and an $\epsilon > 0$ so that all fibers $\mathcal{S}_t$ are transverse to $\partial B$ for $t \in \bar{D}_\epsilon$, where $D_\epsilon = \{|t| < \epsilon\} \subset D$. Moreover, there is a decomposition

$$\mathcal{S}|_{D_\epsilon} = \mathcal{S}' \cup \mathcal{S}'',$$

where $\mathcal{S}' = (\mathcal{S}|_{D_\epsilon}) \cap (D \times B)$ and $\mathcal{S}'' = (\mathcal{S}|_{D_\epsilon}) \cap (D \times \overline{\mathbb{P}^3 - B})$, such that $\mathcal{S}''$ is topologically trivial, and such that the interior of $\mathcal{S}'$ is biholomorphically equivalent to a neighborhood of the origin in the family

$$x^2 + y^2 + z^2 = t. \tag{5.2.1}$$

We take $F$ to be the cubic form associated to a point of $D_\epsilon^*$. Then the family $\mathcal{S}|_{D_\epsilon^*}$ over the punctured disk $D_\epsilon^*$ is homotopy equivalent to the family over the circle with fiber $S$ and with geometric monodromy $\phi : S \longrightarrow S$, where $S = S' \cup S''$, $S' = S \cap B$, $S'' = S \cap \overline{\mathbb{P}^3 - B}$, $\phi|_{S''} = id$, $S'$ is the Milnor fiber of $x^2 + y^2 + z^2 = 0$, and $\phi|_{S'}$ is the monodromy of the corresponding Milnor fibration. See §3 of [19] or §2 and §3 of [30] for more details.

(5.3) Given this discussion, we can understand the topology of $\mathcal{T}|_D$, where $D$ is a disk transverse to $\Delta_s^1$ as above. Pulling back along the projection $\pi : \mathcal{T}|_D \longrightarrow D \times \mathbb{P}^3$, we obtain a decomposition

$$\mathcal{T}|_{D_\epsilon} = \mathcal{T}' \cup \mathcal{T}'', \tag{5.3.1}$$

where $\mathcal{T}''$ is topologically trivial and the interior of $\mathcal{T}'$ is biholomorphic to a neighborhood of the origin in the family

$$x^2 + y^2 + z^2 + w^3 = t , \tag{5.3.2}$$

with $\sigma$ acting by multiplication of $w$ by $\omega$. Moreover, the family $\mathcal{T}|_{D_\epsilon^*}$ over the punctured disk is homotopy equivalent to the family over the circle with fiber $T$ and monodromy $\psi : T \longrightarrow T$ where $T = T' \cup T''$, $\psi|_{T''} = id$, $T'$ is diffeomorphic to the Milnor fiber of $x^2 + y^2 + z^2 + w^3 = 0$, and $\psi|_{T'}$ is the monodromy of the corresponding Milnor fibration. Recall that the *vanishing homology* of the family $\mathcal{T}|_D$ is the image of $H_3(T', \mathbb{Z})$ in $H_3(T, \mathbb{Z})$ and that the *vanishing cohomology* $V$ is the subspace of $H^3(T, \mathbb{Z})$ which is Poincaré dual to the vanishing homology. The space $T'$ itself is the space of (geometric) vanishing cycles, cf. §3 of [19]. Let $\Psi$ denote the monodromy transformation $\psi^*$ on cohomology. It is clear that $\Psi$ acts $\mathcal{E}$-linearly on $\Lambda(T)$, and that $\sigma$ preserves $V$, so that $V$ is an $\mathcal{E}$-submodule of $\Lambda(T)$.

(5.4) **Lemma.** *Let $F_0 \in \Delta_s^1$. Then the vanishing cohomology $V$ for the family (5.3.1) has rank one as an $\mathcal{E}$-submodule of $\Lambda(T)$, with a generator $\delta$ satisfying $h(\delta, \delta) = 1$. The monodromy tranformation $\Psi$ acts on $V$ by multiplication by $-\omega$ and on the subspace $V^\perp \subset \Lambda(T)$ by the identity. Therefore $\Psi$ is the complex $(-\omega)$-reflection in $\delta$.*

(5.5) **Proof:** The monodromy tranformation $\Psi$ on cohomology is $\psi^*$, where $\psi : T \longrightarrow T$ and $T = T' \cup T''$ and $\psi$ preserves the decomposition, as in (5.3). Now $T' \cap T'' = \partial T'$ is the boundary of the Milnor fiber, and for the singularity (5.3.2) it is well-known to be a homology sphere, and in fact is a sphere. See the first example in §9 of [35]. From the Meyer-Vietoris sequence it is then clear that

$$H^3(T, \mathbb{Z}) = H^3(T', \mathbb{Z}) \oplus H^3(T'', \mathbb{Z}) = V \oplus V^\perp$$

is a direct sum decomposition over $\mathbb{Z}$, orthogonal with respect to the symplectic pairing, and that $\Psi$ acts as the identity on the second summand. The action of $\Psi$ on the first summand is well known, and can be derived from the Sebastiani-Thom theorem [46]. Note first that the space $V(k)$ of vanishing cycles for the degeneration of zero-dimensional varieties $z^k = t$ is the $\mathbb{Z}$-module spanned by the differences of roots $a_i - a_{i+1}$. Taken with a counterclockwise ordering of the solution set $\{ a_i(t) : i = 0 \ldots k - 1 \}$ of $z^k = re^{i\theta}$ on the circle of radius $r^{1/k}$, the monodromy transformation on $\mathcal{V}(k)$ is given by cyclic shift: $\Psi(a_i - a_{i+1}) = a_{i+1} - a_{i+2}$ (indices viewed modulo $k$). For a sum of powers $z_1^{k_1} + \cdots + z_n^{k_n}$, the space of vanishing cycles is isomorphic to the tensor product $V(k_1) \otimes \cdots \otimes V(k_n)$ and the monodromy transformation is the tensor product of the monodromy

transformations:
$$\Psi = \Psi_1 \otimes \cdots \otimes \Psi_n.$$

In our case we may write the monodromy transformation as
$$\Psi = (-1) \otimes (-1) \otimes (-1) \otimes \omega.$$

Thus $V$ has rank 2 and $\Psi$ has order six. The action of $\sigma : T \to T$ with respect to this basis of $V$ is evident from (5.3.2), namely
$$\sigma = 1 \otimes 1 \otimes 1 \otimes \omega\ .$$

Thus $\Psi$ is a complex $(-\omega)$-reflection. Since this module is a direct summand of $\Lambda(T)$, any generator $\delta$ (there are six of them) must satisfy $h(\delta,\delta) = \pm 1$. At the end of the next paragraph we will show that $h(\delta,\delta) = 1$, thereby concluding the proof of the lemma.

(5.6) Suppose now that $F_0 \in \Delta_s^k$, where $k = 2, 3, 4$, so that $F_0$ has $k$ nodes and no other singularities. This situation is a straightforward generalization of the case $k = 1$ which has just been considered. A transversal to $\Delta_s^k$ at $F_0$ is of the form $D_\epsilon^k$, and it parametrizes a family $F_{t_1,\cdots,t_k}$ which in $k$ disjoint balls is defined by equations of the form (5.2.1). The family $\mathcal{T}|_{D_\epsilon^k}$ decomposes as
$$\mathcal{T}|_{D_\epsilon^k} = \mathcal{T}' \cup \mathcal{T}'', \tag{5.6.1}$$

where $\mathcal{T}''$ is topologically trivial, and the interior of $\mathcal{T}'$ is a disjoint union of $k$ open sets, each biholomorphic to a neighborhood of the origin in the family (5.3.2). Moreover, the family $\mathcal{T}|_{D_\epsilon^k}$ is homotopy equivalent to the family over the torus $(S^1)^k$, where the monodromy $\psi_i : T \longrightarrow T$ corresponding to the $i$-th factor is supported in the $i$-th component of $T'$. Reasoning as in (5.3) and (5.5), the space $V \subset H^3(T, \mathbb{Z})$ of vanishing cycles is of the form $V = V_1 \oplus \cdots \oplus V_k$, where the sum is orthogonal with respect to $h$, where each $V_i$ is a primitive $\mathcal{E}$-submodule of $H^3(T, \mathbb{Z})$ of rank one, and where the $i$-th monodromy tranformation $\Psi_i$ acts as multiplication by $-\omega$ on $V_i$. If we choose a generator $\delta_i$ of $V_i$, it is clear that $h(\delta_i, \delta_i) = \pm 1$. This value is independent of the choice and also independent of $i$. Since the signature of $h$ is $(4, 1)$ and there can be up to 4 independent vanishing cycles $\delta_i$ we see that $h(\delta_i, \delta_i) = -1$ is impossible, and so $h(\delta_i, \delta_i) = 1$ for all $i$. This completes the proof of Lemma 5.4 and also of the following lemma:

**(5.7) Lemma.** *Let $F_0 \in \Delta_s^k$. The vanishing cohomology $V$ of the family (5.6.1) splits as an $h$-orthogonal direct sum $V = V_1 \oplus \cdots \oplus V_k$ where each $V_i$ is a free $\mathcal{E}$-submodule of $\Lambda(T)$ of rank one generated by an element $\delta_i$ with $h(\delta_i, \delta_i) = 1$. The monodromy transformation $\Psi_i$ corresponding to the $i$-th factor of $D_\epsilon^k$ acts on $\Lambda(T)$ by the $(-\omega)$-reflection in $\delta_i$, and the local monodromy group at $F_0$ is isomorphic to $(\mathbb{Z}/6\mathbb{Z})^k$, generated by $(-\omega)$-reflections in the $\delta_i$.*

(5.8) It is easier to give a concrete description of a framed stable form than of a marked stable form. Suppose $Pf \in \mathcal{F}_s$ lies over $F \in \Delta_s^k$. By mimicking the proof of Lemma 4.11, with the local system $\Lambda(\mathcal{T})$ in place of $L(\mathcal{S})$, one shows that $Pf$ determines (i) a bijection $b_{Pf}$ between the nodes $p_i$ of $S$ and a collection of mutually orthogonal sublattices $R_i$ of $\Lambda$, each generated by a short root, and (ii) an isometric embedding $i_{Pf}$ of $\Lambda(T)$ into the orthogonal complement in $\Lambda$ of $\oplus R_i$, determined up to scalar multiplication. Here, by $\Lambda(T)$ we mean $H^3(T,\mathbb{Z})$ equipped with the hermitian $\mathcal{E}$-module structure defined by the same formulas as for smooth cubic forms; the arguments of (5.5) and (5.6) show that $\Lambda(T) \cong \mathcal{E}^{4-k,1}$. This situation is simpler than the one for marked forms because $\Lambda(T)$ and all the $R_i$ are summands of $\Lambda$, and so $P\operatorname{Aut}\Lambda$ acts transitively on pairs $(b,i)$ satisfying conditions (i) and (ii). Therefore the various framings of $F$ are in one-to-one correspondence with such pairs. When $k=1$, the situation is even simpler: the framings of $\Lambda$ are just the isometric embeddings $f:\Lambda(T) \to \Lambda$ modulo scalar multiplication. Consequently we can explicitly describe a neighborhood in $\mathcal{F}_s$ of a framed 1-nodal form. For $F_0 \in \Delta_s^1$ let $U = U_0 \times D$ be a neighborhood of $F_0$ as in (5.2), let $\mathcal{V} \subset \Lambda(\mathcal{T}|_{U-\Delta})$ be the local system of vanishing cohomology, and let $\mathcal{I} = \mathcal{V}^\perp$ be the local system of invariant cohomology. The deformation-retraction of $\mathcal{T}|_U$ to $T_0$ induces an isomorphism $\mathcal{I} \cong (U-\Delta) \times \Lambda(T_0)$. Then we have the following result, whose proof is easy.

**(5.9) Lemma.** *Let $F_0$ and $U$ be as above, and suppose $Pf_0 \in \mathcal{F}_s$ lies over $F_0$ and is represented by the embedding $f_0 : \Lambda(T_0) \to \Lambda$. Then the component $W$ of $\mathcal{F}_s|_U$ containing $Pf_0$ is $W_0 \cup (W - W_0)$, where $W_0$ is the component of $\mathcal{F}_s|_{U_0}$ containing $Pf_0$ and $W - W_0$ consists of those $Pf \in \mathcal{F}_0$ lying over forms $F \in U - U_0$ such that $f|_{\mathcal{I}_F}$ agrees with $f_0$ up to a scalar factor.*

## 6. Fractional differentials and extension of the period map

The purpose of this section is to concretely identify the extension of $g$ to the set of nodal surfaces. This requires reformulating our description of $H^3(T,\mathbb{C})$ and the various structures on it in terms of suitable differentials on $\mathbb{P}^3$, and an explicit understanding of the period map in terms of these differentials. Fractional differentials for double covers of $\mathbb{P}^3$ were studied by Clemens in [14]. Their theory, including the relation of the filtration by order of pole to the Hodge filtration, parallels that of rational integrals on projective introduced by Griffiths in [24]. The same is true of the theory for general cyclic covers, despite some additional structure (eigenvalue decompositions) and subtleties ($A_k$ singularities occur for $k$-fold covers when the branch locus acquires a node).

(6.1) For $F \in \mathcal{C}_0$, let $\mathcal{L}_F$, or simply $\mathcal{L}$, denote the sheaf over $\mathbb{P}^3$ defined by the following exact sequence of sheaves with $\sigma$-action:

$$0 \longrightarrow \mathbb{Z} \longrightarrow p_*(\mathbb{Z}) \longrightarrow \mathcal{L} \longrightarrow 0. \tag{6.1.1}$$

The subsheaf $\mathbb{Z}$ is the one on which $\sigma$ acts as the identity. It is easy to see that $\mathcal{L}|_{\mathbb{P}^3-S}$ is a local system with fiber $\mathbb{Z}^2$ and that $\sigma$ acts on it by an automorphism of order 3. Defining the action of $\omega \in \mathcal{E}$ on $\mathcal{L}|_{\mathbb{P}^3-S}$ by this action of $\sigma$, we see that the sheaf has the structure of a local system of rank one $\mathcal{E}$-modules. Furthermore, the positive generator of $\pi_1(\mathbb{P}^3-S)$ acts as multiplication by $\omega$. It is also easy to see that the sheaf $\mathcal{L}$ is concentrated on $\mathbb{P}^3-S$, namely $\mathcal{L} = i_!(\mathcal{L}|_{\mathbb{P}^3-S})$, where $i:\mathbb{P}^3 - S \longrightarrow \mathbb{P}^3$ is the inclusion and $i_!$ denotes, as usual, extension by zero. See, for example, lecture 8 of [16]. Since $p: T \longrightarrow \mathbb{P}^3$ is a finite map, $H^3(T,\mathbb{Z}) \cong H^3(\mathbb{P}^3,p_*(\mathbb{Z}))$. From the long exact sequence associated to (6.1.1), the vanishing of $H^3(\mathbb{P}^3,\mathbb{Z})$, and the injectivity of $H^4(\mathbb{P}^3,\mathbb{Z}) \longrightarrow H^4(T,\mathbb{Z})$, one obtains an isomorphism $H^3(\mathbb{P}^3,p_*(\mathbb{Z})) \cong H^3(\mathbb{P}^3,\mathcal{L})$. Combining

these two isomorphisms we obtain the first natural isomorphism of $\mathcal{E}$-modules displayed below. The second isomorphism results from the concentration of $\mathcal{L}$ on $\mathbb{P}^3 - S$:

$$H^3(T, \mathbb{Z}) \cong H^3(\mathbb{P}^3, \mathcal{L}) \cong H^3(\mathbb{P}^3 - S, \mathcal{L}). \tag{6.1.2}$$

Of course we also have $H^3(T, \mathbb{C}) \cong H^3(\mathbb{P}^3, \mathcal{L} \otimes \mathbb{C}) \cong H^3(\mathbb{P}^3 - S, \mathcal{L} \otimes \mathbb{C})$. The action of $\sigma$ decomposes $\mathcal{L} \otimes \mathbb{C}$ into its eigensheaves $\mathcal{L}_\omega$ and $\mathcal{L}_{\bar\omega}$, consequently giving us alternate descriptions of $H^3_\omega(T)$ and $H^3_{\bar\omega}(T)$ by local systems on $\mathbb{P}^3$. In particular $H^3_{\bar\omega}(T, \mathbb{C}) \cong H^3(\mathbb{P}^3 - S, \mathcal{L}_{\bar\omega})$.

(6.2) Similar considerations apply to sheaf cohomology. For the sheaves of holomorphic differential forms $\Omega^j_T$, $j = 0, 1, 2, 3$ we have that $H^i(T, \Omega^j)_{\bar\omega} \cong H^i(\mathbb{P}^3, p_*(\Omega^j_T)_{\bar\omega})$. Similarly, for the sheaves $\Omega^j_T(\infty S)$ of meromorphic differential $j$-forms on $T$ with poles (of arbitrary order) only on $S$, we have $H^i(T, \Omega^j_T(\infty S))_{\bar\omega} \cong H^i(\mathbb{P}^3, p_*(\Omega^j_T(\infty S))_{\bar\omega})$. The sheaves $p_*(\Omega^j_T)_{\bar\omega}$ and $p_*(\Omega^j_T(\infty S))_{\bar\omega}$ are no longer concentrated on $\mathbb{P}^3 - S$. In order to describe these sheaves, let us introduce the temporary notation $\mathcal{B} = i_*(\mathcal{L}_{\bar\omega}|_{\mathbb{P}^3 - S} \otimes \Omega^j_{\mathbb{P}^3}|_{\mathbb{P}^3 - S})$, and define the following subsheaves of $\mathcal{B}$. First, let $\Omega^j_{\mathbb{P}^3}(F^{-1/3})$ denote the subsheaf of $\mathcal{B}$ defined by the conditions (a) $\Omega^j_{\mathbb{P}^3}(F^{-1/3})|_{\mathbb{P}^3 - S} = \mathcal{B}|_{\mathbb{P}^3 - S}$, (b) the stalk of $\Omega^j_{\mathbb{P}^3}(F^{-1/3})$ at a point $P$ of $S$ consists of the elements of the stalk of $\mathcal{B}$ at $P$ which are of the form

$$\alpha w^{2/3} + \beta \wedge dw/w^{1/3}, \tag{6.2.1}$$

where $w$ is a local defining equation of $S$ at $P$, and (c) $\alpha$, $\beta$ are in the stalk of $\Omega^j_{\mathbb{P}^3}$, $\Omega^{j-1}_{\mathbb{P}^3}$ at $P$. Observe that the multivalued functions $w^{2/3}$ and $1/w^{1/3}$ define in a natural way single-valued sections of $\mathcal{L}_{\bar\omega} \otimes \mathcal{O}_{\mathbb{P}^3}|_{\mathbb{P}^3 - S}$ in $U - S$, where $U$ is a neighborhood of $P$ in $\mathbb{P}^3$, and therefore define sections of $i_*((\mathcal{L}_{\bar\omega} \otimes \mathcal{O}_{\mathbb{P}^3})|_{\mathbb{P}^3 - S})$ on $U$. See for instance the discussion in §2 of [20]. Thus the expression (6.2.1) makes sense. We also define $\Omega^j_{\mathbb{P}^3}(F^{-(\infty+1/3)})$ to be the subsheaf of $\mathcal{B}$ which equals $\mathcal{B}$ on $\mathbb{P}^3 - S$, and whose stalk at $P \in S$ consists of the elements of the form

$$w^{-l}\alpha/w^{1/3}, \tag{6.2.2}$$

where $l$ is a non-negative integer and $w$ and $\alpha$ are as above. Observe that the exterior derivative $d$ makes $\Omega^\bullet_{\mathbb{P}^3}(F^{-(\infty+1/3)})$ into a complex of sheaves on $\mathbb{P}^3$. This complex has a filtration which we call the *pole filtration*: the smallest filtration for which the expressions (6.2.2) have filtration level $\leq l$. The complex $\Omega^\bullet_{\mathbb{P}^3}(F^{-1/3})$ can be characterized as the *subcomplex* of $\Omega^\bullet_{\mathbb{P}^3}(F^{-(\infty+1/3)})$ consisting of all elements $\eta$ of filtration level zero such that $d\eta$ also has filtration level zero.

**(6.3)Lemma.** *There are natural isomorphisms $(p_*(\Omega^\bullet_T))_{\bar\omega} \cong \Omega^\bullet_{\mathbb{P}^3}(F^{-1/3})$ and $(p_*(\Omega^\bullet_T(\infty S)))_{\bar\omega} \cong \Omega^\bullet_{\mathbb{P}^3}(F^{-(\infty+1/3)})$ of complexes of sheaves over $\mathbb{P}^3$.*

**Proof:** Since the cyclic transformation $\sigma$ acts freely on $T - S$, the isomorphisms over $\mathbb{P}^3 - S$ are clear from the definitions of the objects involved. At a point $P \in S$, choose local coordinates $x, y, z$ centered at $P$ on $T$ and local coordinates $x, y, w$ centered at $P$ on $\mathbb{P}^3$ so that $\sigma(x, y, z) = (x, y, \omega z)$ and $p(x, y, z) = (x, y, z^3)$. A local holomorphic $j$-form on $T$ at $P$ has eigenvalue $\bar\omega$ if and only if it is of the form $az^2 + bz\, dz$, where $a$ and $b$ are $\sigma$-invariant forms. Therefore $a = p^*\alpha$ and $b = p^*\beta_1$ for suitable local holomorphic forms $\alpha$ and $\beta_1$ on $\mathbb{P}^3$. Setting $z = w^{1/3}$ and $\beta = \beta_1/3$ we obtain the expression (6.2.1) and the first assertion of the lemma follows. The second assertion follows by a similar (but slightly longer) calculation.

(6.4) Since $T - S$ is an affine variety, $H^*(T - S, \mathbb{C})$ is the cohomology of the complex of meromorphic forms on $T$ with poles only on $S$ (complex of global sections on $T$ of $\Omega_T^\bullet(\infty S)$) and consequently its summand $H^*(T - S, \mathbb{C})_{\bar{\omega}} \cong H^*(T, \mathbb{C})_{\bar{\omega}}$ is the cohomology of global sections on $\mathbb{P}^3$ of $\Omega_{\mathbb{P}^3}^\bullet(F^{-(\infty+1/3)})$. We write $\mathcal{A}_{1/3}^\bullet(\mathbb{P}^3)$ for this complex of global sections. Elements of $\mathcal{A}_{1/3}^\bullet(\mathbb{P}^3)$ are those forms which can be written as

$$\phi = \sum \frac{A_J \Omega_J}{F^{l+1/3}}, \qquad (6.4.1)$$

where $l$ is a non-negative integer, and where $\Omega_J = i(\partial/\partial X_{j_1}) \cdots i(\partial/\partial X_{j_q})\Omega$ for $J = (j_1, \ldots, j_q)$, where

$$\Omega = \sum_k (-1)^k X_k dX_0 \wedge \cdots \wedge \widehat{dX_k} \wedge \cdots \wedge dX_3.$$

Each $A_J$ is a homogeneous polynomial such that the expression (6.4.1) has degree zero. Thus the degree of the denominator equals the degree of the numerator: $3l + 1 = \deg A_J + 4 - |J|$. The pole filtration on $\Omega_{\mathbb{P}^3}^\bullet(F^{-(\infty+1/3)})$ induces a filtration of the same name on $\mathcal{A}_{1/3}^\bullet$. It the smallest filtration for which all expressions (6.4.1) have filtration level $\leq l$.

**(6.5) Theorem.** *Let $F \in \mathcal{C}_0$. Then the following statements hold.*

(1) *The cohomology of the complex $\mathcal{A}_{1/3}^\bullet$ is naturally isomorphic to $H_{\bar{\omega}}^\bullet(T)$.*

(2) *The filtration on $H_{\bar{\omega}}^3(T)$ induced by the pole filtration on $\mathcal{A}_{1/3}^\bullet$ agrees with the Hodge filtration in the sense that a class has pole filtration $l$ if and only if it has Hodge filtration $3 - l$.*

(3) $H_{\bar{\omega}}^{3,0}(T) = 0$ *and $H_{\bar{\omega}}^{2,1}(T)$ is the one-dimensional subspace of the five-dimensional space $H_{\bar{\omega}}^3(T)$ spanned by the cohomology class of $\Omega/F^{4/3}$.*

**Proof:** The proof of the first statement has just been given. The argument for the second part is by now standard, and we sketch the version that best suits our purposes, following the version in [12] of the original version of such theorems in [24]. For $i = 0, \ldots, 3$ let $U_i \subset \mathbb{P}^3$ be the set where $\partial F/\partial X_i \neq 0$. Since $F$ is not singular, $\mathcal{U} = \{U_i\}$ forms an open cover of $\mathbb{P}^3$ by affine open sets. The Čech-deRham bicomplex $C^\bullet(\mathcal{U}, \Omega_{\mathbb{P}^3}^\bullet(F^{-(\infty+1/3)}))$ contains both $\mathcal{A}_{1/3}^\bullet(\mathbb{P}^3)$ and $C^\bullet(\mathcal{U}, \Omega_{\mathbb{P}^3}^\bullet(F^{-1/3}))$ as quasi-isomorphic subcomplexes. The first subcomplex defines the pole filtration of $H^3(T)_{\bar{\omega}}$, the second subcomplex defines the Hodge filtration of $H^3(T)_{\bar{\omega}}$, and the two subcomplexes can be compared inside the larger complex. To do this, one uses the argument of §3b of [12]. The only changes needed are to replace the factor $1/(1 - l)$ in the definition of the partial homotopy $H$ at the top of p. 59 by the factor $1/(\frac{3}{2} - l)$, and to replace the complex of logarithmic forms by our complex $\Omega_{\mathbb{P}^3}^\bullet(F^{-1/3})$.

For the proof of the third part, simply observe that there is no non-zero form in $\mathcal{A}_{1/3}^3$ of pole filtration level 0, that is, with denominator $F^{1/3}$ and $J = \emptyset$ in the expression (6.4.1). Thus $H^{3,0}(T)_{\bar{\omega}} = 0$. Likewise, any form in $\mathcal{A}_{1/3}^3$ of pole filtration level 1, that is, with denominator $F^{4/3}$ in (6.4.1) and $J = \emptyset$ is a constant multiple of $\Omega/F^{4/3}$. This completes the proof of the theorem.

(6.6) **Remark.** The reader who would like to verify all the details of the proof we have just sketched may find it instructive to verify the following fact. Begin with the element $\Omega/F^{4/3}$ apply the partial homotopy to find an element of $C^\bullet(\mathcal{U}, \Omega^\bullet_{\mathbb{P}^3}(F^{-1/3}))$ cohomologous to the given one in the large bicomplex. The result has components of Čech-deRham bidegree $(0,3)$ and $(1,2)$. Its component of bidegree $(0,3)$ is, up to a constant factor, the Čech zero-cochain that to the open set $U_i$ assigns the differential form

$$\frac{F_i d\Omega_i \pm F_{ii}\Omega}{F_i^2 F^{1/3}}$$

on $U_i$, where $F_i = \partial F/\partial X_i$, $F_{ii} = \partial^2 F/\partial X_i^2$, and where $\Omega_i$ is as defined in (6.4). Its component of bidegree $(1,2)$ is a constant multiple of the Čech one-cochain that to $U_i \cap U_j$ assigns the differential form

$$\pm \frac{\Omega_{ij}}{F_i F_j F^{1/3}}$$

on $U_i \cap U_j$. The last two statements give, in explicit form, the assertion made about the Hodge filtration of the class of $\Omega/F^{4/3}$.

(6.7) We are now in a position to give explicit formulas for the period map. Let $Pf \in \mathcal{F}_0$ lie over $F \in \mathcal{C}_0$. Let $\Lambda' = \text{Hom}_{\mathbb{Z}}(\Lambda, \mathbb{Z})$ be the lattice dual to $\Lambda$, and let $f' : \Lambda' \longrightarrow \Lambda(T)' = H_3(T, \mathbb{Z}) \cong H_3(\mathbb{P}^3 - S, \mathcal{L}')$ be the map dual to $f$, where $\mathcal{L}' = Hom_{\mathbb{Z}}(\mathcal{L}, \mathbb{Z})$ is the dual local system. Observe that $\Lambda'$ and $H_3(T, \mathbb{Z})$ are $\mathcal{E}$-modules, that $\mathcal{L}'$ is a local system of $\mathcal{E}$-modules, and that $f'$ and the displayed isomorphisms are $\mathcal{E}$-linear. Moreover, $H^3_{\bar{\omega}}(T)$ is naturally isomorphic to the subspace $\text{Hom}(H_3(\mathbb{P}^3 - S, \mathcal{L}'), \mathbb{C})_{\bar{\omega}}$ of $\text{Hom}_{\mathbb{Z}}(H_3(\mathbb{P}^3 - S, \mathcal{L}'), \mathbb{C})$ consisting of all homorphisms $\alpha$ so that $\alpha(\sigma_*(x)) = \bar{\omega}\alpha(x)$. The subspace $\text{Hom}(\Lambda', \mathbb{C})_{\bar{\omega}}$ defined in the same way is isomorphic to $\Lambda \otimes_{\mathcal{E}} \mathbb{C} \cong \mathbb{C}^{4,1}$. The one-dimensional subspace $H^{2,1}_{\bar{\omega}}(T) \subset H^3_{\bar{\omega}}(T) = \text{Hom}(H_3(\mathbb{P}^3 - S, \mathcal{L}'), \mathbb{C})_{\bar{\omega}}$ is generated by the following function of $x \in H_3(\mathbb{P}^3 - S, \mathcal{L}')$:

$$x \longrightarrow \int_x \frac{\Omega}{F^{4/3}}.$$

Consequently, $g(Pf) \in \mathbb{C}H^4$ is the negative line in $\mathbb{C}^{4,1} = \text{Hom}(\Lambda', \mathbb{C})_{\bar{\omega}}$ generated by the following function of $x \in \Lambda'$:

$$x \longrightarrow \int_{f'(x)} \frac{\Omega}{F^{4/3}}.$$

Even more concretely, for a given $Pf_0 \in \mathcal{F}_0$ lying over $F_0 \in \mathcal{C}_0$ with zero set $S_0$, choose a tubular neighborhood $N$ of $S_0$ in $\mathbb{P}^3$, and choose a contractible neighborhood $U \subset \mathcal{C}_0$ of $F_0$ so that $\{F = 0\} \subset N$ for all $F \in U$. Since $\mathbb{P}^3 - \bar{N}$ is a deformation retract of $\mathbb{P}^3 - S_0$, there is a canonical isomorphism $H_3(\mathbb{P}^3 - \bar{N}, \mathcal{L}') \cong H_3(\mathbb{P}^3 - S_0, \mathcal{L}')$. Let $c_0, \cdots, c_4$ be the $\mathcal{E}$-basis of $H_3(\mathbb{P}^3 - \bar{N}, \mathcal{L}')$ which corresponds, under this isomorphism, to $f'_0(e_0), \cdots, f'_0(e_4)$, where $e_0, \cdots, e_4$ consists of the first five elements of the $\mathbb{Z}$-basis $e_0, \cdots, e_4, \omega e_0, \cdots, \omega e_4$ of $\Lambda'$ dual to the standard $\mathbb{Z}$-basis $(1, 0, \cdots, 0), \cdots, (0, \cdots, \omega)$ of $\Lambda$. Let $\tilde{U} \subset \mathcal{F}_0$ be the unique neighborhood of $Pf_0$ that lies isomorphically over $U \subset \mathcal{C}_0$. Then for $Pf \in \tilde{U}$ the period map is

$$g(Pf) = \left( \int_{c_0} \frac{\Omega}{F^{4/3}} : \cdots : \int_{c_4} \frac{\Omega}{F^{4/3}} \right). \tag{6.7.1}$$

From this we see immediately that $g$ is holomorphic, as asserted in (2.16).

(6.8) Now we study the case where the surface $S$ has nodes. Let $F \in \Delta_s$ and let $\Sigma \subset S$ denote the set of nodes of $S$ (the singular set of $S$), and let $\Sigma$ also denote the singular set of $T$. We can repeat all the definitions of (6.1) to (6.4) and the proof of Theorem 6.5 in this context, with only two exceptions: $T$ now has singularities, and the sets $U_i = \{\partial F/\partial X_i \neq 0\}$, defined and used in the proof of Theorem 6.5 no longer cover $\mathbb{P}^3$, but cover $\mathbb{P}^3 - \Sigma$. The proof of Theorem 6.5 will then work for the non-compact varieties $T - \Sigma$ and $\mathbb{P}^3 - \Sigma$, but we prefer to reformulate it in terms of compact non-singular varieties as follows. Blow up $\mathbb{P}^4$ at the singular set $\Sigma$ of $T$ and let $\tilde{T}$ denote the proper transform of $T$. It is easy to check from the local equation $x^2 + y^2 + z^2 + w^3 = 0$ of $T$ at each point of $\Sigma$ that $\tilde{T}$ is smooth, it is a resolution of singularities of $T$, and that each component of the exceptional divisor $E$ (the preimage of $\Sigma$ in $\tilde{T}$) is a singular quadric (the cone on a non-singular conic). Moreover the cyclic transformation $\sigma$ lifts to $\tilde{T}$. Consider the following maps

$$\begin{array}{ccccccc} T - S & \subset & T - \Sigma & = & \tilde{T} - E & \subset & \tilde{T} \\ \downarrow & & \downarrow & & & & \\ \mathbb{P}^3 - S & \subset & \mathbb{P}^3 - \Sigma & & & & \end{array}$$

It is easy to check (using the fact that the links of the singularities of $T$ are spheres) that the top right hand inclusion induces an isomorphism on $H^3$ with $\mathbb{Z}$ or $\mathbb{C}$ coefficients. It is clear that the cohomological formalism used in (6.1)–(6.4) applies to the two vertical arrows and the two left-hand inclusions. Putting all this together, we obtain a natural isomorphism of $H^3(\tilde{T})_{\bar{\omega}}$ with the cohomology of the complex $\mathcal{A}_{1/3}^\bullet$.

**(6.9) Theorem.** *Let $F \in \Delta_s^k$. Then the following statements hold.*

(1) $H^3(\mathcal{A}_{1/3}^\bullet)$ *is naturally isomorphic to $H^3_{\bar{\omega}}(\tilde{T})$.*

(2) *The filtration on $H^3_{\bar{\omega}}(\tilde{T})$ induced by pole filtration on $\mathcal{A}_{1/3}^\bullet$ agrees with the Hodge filtration on $H^3_{\bar{\omega}}(\tilde{T})$, in the sense that a class has pole filtration $l$ if and only if it has Hodge filtration $3 - l$.*

(3) $H^{3,0}_{\bar{\omega}}(\tilde{T}) = 0$ *and $H^{2,1}_{\bar{\omega}}(\tilde{T})$ is the one-dimensional subspace of the $5 - k$-dimensional space $H^3_{\bar{\omega}}(\tilde{T})$ spanned by the cohomology class of $\Omega/F^{4/3}$.*

**Proof:** The proof of the first statement was just explained in (6.8). For the second statement, it is clear that the obvious modification of the proof of Theorem 6.5 shows that the pole filtration on $H^3_{\bar{\omega}}(T - \Sigma)$ agrees with its Hodge filtration. Since the inclusion $T - \Sigma \subset \tilde{T}$ is holomorphic, the isomorphism $H^3_{\bar{\omega}}(\tilde{T}) \longrightarrow H^3_{\bar{\omega}}(\tilde{T} - E) = H^3_{\bar{\omega}}(T - \Sigma)$ carries the Hodge filtration of the first to that of the last, and therefore must also be an isomorphism of Hodge filtrations.

(6.10) We now combine these results with those of the previous section to explicitly identify the extension of the period map obtained in (3.16) to $\mathcal{F}_s|_{\Delta_s^1}$. Let $F_0 \in \Delta_s^1$, and recall the notation of (5.2). Let $F \in \mathcal{C}_0$ be associated to some point of $D_\epsilon^*$, so that $S$ is the corresponding fiber of $\mathcal{S}$. Note that the decomposition $H^3(T, \mathbb{Z}) = V \oplus V^\perp$ of (5.5) is dual to the decomposition

$$H_3(T, \mathbb{Z}) \cong H_3(\mathbb{P}^3 - S, \mathcal{L}') = H_3(B - S, \mathcal{L}') \oplus H_3(\overline{\mathbb{P}^3 - B} - S, \mathcal{L}'), \qquad (6.10.1)$$

where $\mathcal{L}'$ is the local system on $\mathbb{P}^3 - S$ dual to $\mathcal{L}$ as in (6.7) and $B$ is the ball in (5.2). In the notation of (5.2), let $N''$ be the union of the fibers of $\mathcal{S}''$. Then $\bar{N}''$ is a closed tubular neighborhood of $S_0''$ in $\mathbb{P}^3 - B$ and $\mathbb{P}^3 - B - \bar{N}''$ is a deformation retract of $\overline{\mathbb{P}^3 - B} - S$. Choose a basis $c_0, c_1, c_2, c_3$ for the $\mathcal{E}$-module $H_3(\mathbb{P}^3 - B - \bar{N}'', \mathcal{L}')$. Then their images in the second summand of (6.10.1) give a basis that is independent of $S$ (in other words, constant with respect to the Gauss-Manin connection) which is also a basis for $t = 0$. Moreover, in analogy with the choice of basis in (6.7), we choose this basis so that the first four elements $c^0, c^1, c^2, c^3$ of the $\mathbb{Z}$-basis of $H^3(T'', \mathbb{Z})$ dual to the $\mathbb{Z}$-basis $c_0, \cdots, c_3, \omega c_0, \cdots, \omega c_3$ of $H_3(T'', \mathbb{Z})$ satisfy $h(c^0, c^0) = -1$ and $h(c^i, c^i) = 1$ for $i > 0$. Next, choose a sector $I$ in the punctured disk $D_\epsilon^*$, choose a proper annulus $A \subset D_\epsilon^*$, and let $C$ denote the closure of $A \cap I$. Let $N'$ be the union of the fibers of $\mathcal{S}'$ lying over $C$. Then $B - S$ deformation retracts to $B - N'$. Let $c_4$ denote a generator of the rank one $\mathcal{E}$-module $H_3(B - N', \mathcal{L}')$, and let $c_4$ also denote its image in the first summand of (6.10.1). Finally, let

$$\phi_i(t) = \int_{c_i} \frac{\Omega}{F_t^{4/3}}. \tag{6.10.2}$$

Then we have the following result.

**(6.11) Lemma.** *For $0 \leq i \leq 3$, $\phi_i(t)$ is holomorphic in $D_\epsilon$, and $\phi_0(0) \neq 0$. There exists a holomorphic function $\eta(t)$ in $D_\epsilon$ so that $\phi_4(t) = t^{1/6}\eta(t)$ for all $t$ in the sector $I$.*

**Proof:** It is clear from the formula (6.10.2) and the choice of the cycles $c_0, \cdots, c_3$ that the $\phi_i$ are holomorphic in $D_\epsilon$ for $0 \leq i \leq 3$. Since the vector $(\phi_0(0), \cdots, \phi_3(0))$ spans the one-dimensional space $H_{\bar\omega}^{2,1}$, it is not zero. Therefore $-|\phi_0(0)|^2 + |\phi_1(0)|^2 + |\phi_2(0)|^2 + |\phi_3|^2 < 0$, and so $\phi_0(0) \neq 0$. From Lemma 5.4 it is clear that $\phi_4(t^6)$ can be analytically continued to a single-valued function on $D_\epsilon$. Since $(\phi_0(t), \cdots, \phi_4(t)) \in \mathbb{C}^{4,1}$ has negative norm, it follows from the first part that $\phi_4(t^6)$ is bounded, hence holomorphic at 0. It follows easily from Lemma 5.4 that $\phi_4(t) = t^{1/6}\eta(t)$ for some holomorphic function $\eta$ on $D_\epsilon$, and the proof is complete.

(6.12)  Now let $U = U_0 \times D$ be a neighborhood of $F_0 \in \Delta_s^1$ as in (5.2), and let $Pf_0 \in \mathcal{F}_s$ lie over $F_0$. Recall from (5.8) that $f_0 : \Lambda(T_0) \longrightarrow \Lambda$ is an isometric embedding with image the orthogonal complement of a short root. Assume that this short root is $e^4 = (0, 0, 0, 0, 1) \in \Lambda$, and let $f_0' : \Lambda' \longrightarrow \Lambda(T_0)' \cong H_3(\mathbb{P}^3 - \Sigma, \mathcal{L}') \cong H_3(\overline{\mathbb{P}^3 - B} - S_0, \mathcal{L}') \cong H_3(\mathbb{P}^3 - B - \bar{N}'', \mathcal{L}')$, where $N''$ is as in (6.10). Let $e_0, \cdots, e_4$ be the basis of $\Lambda'$ as in (6.7). Then $f_0'(e_0), \cdots, f_0'(e_3)$ is a basis of $H_3(\overline{\mathbb{P}^3 - B} - \bar{N}'', \mathcal{L}')$ as in (6.10). Now let $W$ be the connected component of $Pf_0$ in $\mathcal{F}_s|_U$ as described in Lemma 5.9. In the notation of Lemma 5.9, if $Pf \in W - W_0$, then $f(e^4)$ is orthogonal to $\Lambda(T_0)$ and $f'(e_4)$ is represented, up to a scalar multiple, by $c_4 \in H_3(B - N', \mathcal{L}')$ as in (6.10). Lemma 6.11 asserts that $g$ is single-valued in $W - W_0$ and that for $Pf_0 \in W_0$

$$g(Pf_0) = \left(\int_{c_0} \frac{\Omega}{F_0^{4/3}}, \cdots, \int_{c_3} \frac{\Omega}{F_0^{4/3}}, 0\right). \tag{6.12.1}$$

By Theorem 6.9, the integrand in this expression is the representative, in $\mathcal{A}_{1/3}^\bullet$, of $H_{\bar\omega}^{2,1}(\tilde{T}_0) \subset H_{\bar\omega}^3(\tilde{T}_0) = H_{\bar\omega}^3(T_0)$. Using the isomorphism in $H^3$ induced by the inclusions $T_0 \leftarrow T_0 - \Sigma = \tilde{T}_0 - E \to \tilde{T}_0$, we can regard $f_0$ as defined on $H^3(\tilde{T}_0)$. Thus Lemma 6.11 asserts that

$$g(Pf_0) = (f_0)_*(H_{\bar\omega}^{2,1}(\tilde{T}_0)) \in \mathbb{C}H(\Lambda \otimes_\mathcal{E} \mathbb{C}) = \mathbb{C}H^4. \tag{6.12.2}$$

Since any $Pf_0$ is equivalent, under an automorphism of $\Lambda$, to one with image the orthogonal complement of $e^4$, this formula for the period map is valid for any $Pf_0$ lying over $F_0 \in \Delta_s^1$. Moreover, the maps just used to define the needed isomorphisms on $H^3$ extend to maps of families over $\Delta_s^1$. This is because the threefolds $\tilde{T}$ for $F$ in $\Delta_s^1$ can be coherently assembled into a family $\pi : \tilde{\mathcal{T}}|_{\Delta_s^1} \to \Delta_s^1$. To construct it, observe that the singular loci of the fibers of $\mathcal{T}$ over $\Delta_s^1$ constitute a smooth subvariety $\Sigma \subset \mathcal{T}|_{\Delta_s^1} \subset \Delta_s^1 \times \mathbb{P}^4$. Therefore one can form $\tilde{\mathcal{T}} \subset B_\Sigma(\Delta_s^1 \times \mathcal{T})$, the proper transform of $\mathcal{T}$ in the blow-up of $\Delta_s^1 \times \mathcal{T}$ along $\Sigma$, and the resulting variety is smooth. Write $\Lambda(\tilde{\mathcal{T}})$ for the sheaf $R^3\pi_*(\mathbb{Z})$, which is a local system over $\Delta_s^1$. Using this notation we can state the following lemma, whose proof is now clear:

**(6.13) Lemma.** $\mathcal{F}_s|_{\Delta_s^1}$ *is isomorphic to the subsheaf of* $PHom(\Lambda(\tilde{\mathcal{T}}), \Delta_s^1 \times \Lambda)$ *consisting of projective equivalence classes of isometric embeddings onto complements of short roots. Under this isomorphism the restriction of the period map $g$ to $\mathcal{F}_s|_{\Delta_s^1}$ is given by formula (6.12.2).*

## 7. The monodromy group and the hyperplane configuration

(7.1) The purpose of this section is to prove Theorem 2.14, which identifies the monodromy group $\Gamma$ of the universal family $\mathcal{T}_0$ of cubic threefolds. We begin by listing some of the known properties of $\Gamma$.

1. $\Gamma$ acts by isometries of the 5-dimensional unimodular $\mathcal{E}$-lattice $\Lambda$ defined in (2.7).

2. By (5.4), $\Gamma$ is the image of $\pi_1(\mathcal{C}_0)$ under a homomorphism carrying meridians of $\Delta$ to $(-\omega)$-reflections in norm one vectors in $\Lambda$. This fact was originally observed in [13].

3. By (5.7) and the existence of a 4-nodal surface, $P\Gamma$ contains a group $(\mathbb{Z}/6)^4$.

4. There is a surjection from the Artin group $\mathcal{A}(E_6)$ to $\pi_1(\mathcal{C} - \Delta)$ such that each of the standard generators of $\mathcal{A}(E_6)$ maps to a meridian of $\Delta$. This result is due to Libgober [29], who states his theorem as a surjection $\mathcal{A}(E_6) \to \pi_1(P\mathcal{C} - P\Delta)$, where as usual $P\mathcal{C}$ denotes the projective space of the space of cubic forms and $P\Delta$ denotes the image in $P\mathcal{C}$ of $\Delta$. The surjection is constructed from a map from the complement of the discriminant in the base space of a universal unfolding of an $E_6$ singularity to the space $\mathcal{C}_0$ of smooth cubic forms. Since it maps a space with fundamental group $\mathcal{A}(E_6)$ to $\mathcal{C}_0$, our statement follows from his argument.

(7.2) For our purposes the Artin group is given by the following generators and relations. Take one generator for each node of the $E_6$ Dynkin diagram and declare that the generators $x$ and $y$ associated to two nodes braid ($xyx = yxy$) if the nodes are joined and commute ($xy = yx$) otherwise. Note that two group elements which braid are conjugate in the group they generate. Each Artin generator acts on $\Lambda$ as a 6-fold reflection, and we name these reflections $R_1, \ldots, R_6$ according to the following scheme:

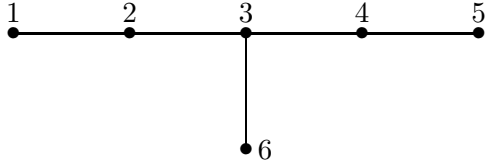

Now choose one of the short roots in which $R_i$ is a $(-\omega)$-reflection and call it $r_i$. Our strategy is to argue that without loss of generality we may take the $r_i$ to be a specific set of short roots, and then determine the group $\Gamma$ generated by the $R_i$.

**(7.3) Lemma.** *The six reflections $R_i$ are distinct.*

**Proof:** Adjoining any relation $R_i = R_j$ to $\Gamma$ (with $i \neq j$) reduces $\Gamma$ to $\mathbb{Z}/6$. To see this, observe that the braid relations imply that if $R_k$ braids with $R_i$ and commutes with $R_j$, then $R_i = R_j$ implies that $R_i = R_j = R_k$. Chasing around the Dynkin diagram, one shows that if any two of the $R_i$ coincide then all of them do. Since $P\Gamma$ contains a group $(\mathbb{Z}/6)^4$, this is impossible. Therefore the $R_i$ are distinct. (One could also use the Zariski-density of $\Gamma$ in $U(4,1)$, proved in [13], instead of the existence of a subgroup $(\mathbb{Z}/6)^4$.)

**(7.4) Lemma.** *The absolute value of $h(r_i, r_j)$ is 0 or 1 according to whether $R_i$ and $R_j$ commute or braid.*

**Proof:** If $R_i$ and $R_j$ commute, then $r_i$ and $r_j$ are either orthogonal or proportional. The latter is impossible by Lemma 7.3. So suppose that $R_i$ and $R_j$ braid. Writing $\alpha$ for $h(r_i, r_j)$, we can compute matrices for the action of the $(-\omega)$-reflections $R_i$ and $R_j$ on the span of $r_i$ and $r_j$:

$$R_i : r_i \mapsto -\omega r_i$$
$$r_j \mapsto r_j - (1+\omega)h(r_j, r_i)\, r_i = r_j + \bar{\omega}\bar{\alpha} r_i$$
$$R_j : r_i \mapsto r_i - (1+\omega)h(r_i, r_j)\, r_j = r_i + \bar{\omega}\alpha r_j$$
$$r_j \mapsto -\omega r_j.$$

Multiplying matrices together, we see that the braid relation requires $|\alpha| = 1$.

(7.5) Since the $(-\omega)$-reflection in $r_i$ is the same as the $(-\omega)$-reflection in any scalar multiple of $r_i$, we may without loss of generality replace the $r_i$ by scalar multiples of themselves. Since the $E_6$ diagram is a tree and pairs of the $r_i$ corresponding to non-adjacent nodes are orthogonal, we may do this in such a way that the inner product of $r_i$ and $r_j$ is 1 or 0, according to whether $R_i$ and $R_j$ braid or not.

(7.6) We have described all the inner products among the $r_i$. One can now compute the determinant of the span of $r_2, \ldots, r_6$, which turns out to be $-1$. This proves that the $r_i$ generate $\Lambda$. It also implies that $\operatorname{Aut} \Lambda$ acts transitively on configurations of vectors whose inner products are those of $r_1, \ldots, r_6$. Indeed, an isometry between two such configurations extends to an isometry between the lattices they span, which is to say an isometry of $\Lambda$. One can also check that the vector $r_1 + r_2 - r_4 - r_5$ is orthogonal to all the $r_i$, and so it vanishes. This identifies the expected linear dependence between 6 vectors spanning a 5-dimensional vector space.

(7.7) We are now ready to realize the $r_i$ as explicit vectors. So far we have described the inner product on $\Lambda$ by (2.7.1), but for the arguments below it is better to use the isometric lattice $I^{\mathcal{E}}_{4,1}$ with inner product matrix

$$A = \begin{pmatrix} 1 & & & & \\ & 1 & & & \\ & & 1 & & \\ & & & 1 & \\ & & & & 1 \end{pmatrix}. \tag{7.7.1}$$

(Blanks indicate zeros.) We regard elements of $I^{\mathcal{E}}_{4,1}$ as column vectors, where $\langle v|w\rangle = w^*Av$ and $w^*$ denotes the conjugate transpose of $w$. The transitivity of $\operatorname{Aut}\Lambda$ on '$E_6$-configurations' implies that we may introduce a basis for $\Lambda$ with inner product matrix $A$, with respect to which $r_1, \ldots, r_6$ have the coordinates below. (The bottom root, $r_7$, will be introduced in (7.17).)

$$\begin{array}{lll}
r_1 = (1, 0, 0, 0, 0) & r_3 = (0, 0, 0, 1, -\omega) & r_5 = (0, 1, 0, 0, 0) \\
r_2 = (1, 0, 0, 0, 1) & r_4 = (0, 1, 0, 0, 1) & \\
& r_6 = (0, 0, 1, 0, 1) & \\
& r_7 = (0, 0, 0, 0, 1) &
\end{array} \tag{7.7.2}$$

Now we will determine the group generated by $R_1, \ldots, R_6$, first considereing $\langle R_1, R_2, R_4, R_5 \rangle$ and then $\langle R_1, \ldots, R_5 \rangle$. Study of the latter group leads us to consider a seventh reflection $R_7$ which plays a key role. We determine the group generated by all the $R_i$ except $R_3$ and then the group $\Gamma = \langle R_1, \ldots, R_7 \rangle$ itself.

(7.8) We now present some background material from [1], specialized to the current situation. (We also adapt the material to work with our convention that a Hermitian form be linear in its first argument rather than in its second.) To begin, write vectors of $I^{\mathcal{E}}_{4,1}$ in the form $(\lambda; \mu, \nu)$ with $\lambda \in \mathcal{E}^3$ and $\mu, \nu \in \mathcal{E}$. Two important sublattices are $I^{\mathcal{E}}_{3,1}$, which consists of the vectors $(\lambda_1, \lambda_2, 0; \mu, \nu)$, and $I^{\mathcal{E}}_{1,1}$, which consists of the vectors $(0, 0, 0; \mu, \nu)$. We distinguish the null vector $\rho = (0; 0, 1)$ and define the height of $v \in I^{\mathcal{E}}_{n+1,1}$ to be $h(v, \rho)$, which is just the second-to-last coordinate of $v$. The matrices

$$T_{\lambda,z} = \begin{pmatrix} I_3 & \lambda & 0 \\ 0 & 1 & 0 \\ -\lambda^* & z - \lambda^2/2 & 1 \end{pmatrix}$$

define isometries of $I^{\mathcal{E}}_{4,1}$ that preserve $\rho$. Here $I_3$ is the $3 \times 3$ identity matrix, $\lambda$ is a column vector $(\lambda_1, \lambda_2, \lambda_3) \in \mathcal{E}^3$ and $z \in \operatorname{Im}\mathbb{C}$ satisfies $z - \lambda^2/2 \in \mathcal{E}$. We call these maps translations, and they satisfy the relations

$$T_{\lambda,z} T_{\lambda',z'} = T_{\lambda+\lambda', z+z'+\operatorname{Im}\langle\lambda'|\lambda\rangle}$$
$$T^{-1}_{\lambda,z} = T_{-\lambda,-z} \tag{7.8.1}$$
$$[T_{\lambda,z}, T_{\lambda',z'}] = T_{0, 2\operatorname{Im}\langle\lambda|\lambda'\rangle}.$$

Here $[T, T'] = TT'T^{-1}T'^{-1}$. Thus the translations form a Heisenberg group: its center and commutator subgroup coincide and consist of the central translations (the translations with $\lambda = 0$). Central translations are also called unitary transvections. Symmetries of $I_{4,1}^{\mathcal{E}}$ that fix $I_{1,1}^{\mathcal{E}}$ pointwise act by conjugation on the translations in the natural way: if $U$ is such a symmetry then

$$UT_{\lambda,z}U^{-1} = T_{U(\lambda),z}. \tag{7.8.2}$$

(7.9) The sublattice $I_{3,1}^{\mathcal{E}}$ is important because it contains $r_1, \ldots, r_5$. This implies that $R_1, \ldots, R_5$ preserve $I_{3,1}^{\mathcal{E}}$. Indeed we will see that $\langle R_1, \ldots, R_5 \rangle$ is almost all of $\operatorname{Aut} I_{3,1}^{\mathcal{E}}$. By a translation of $I_{3,1}^{\mathcal{E}}$ we mean a translation that preserves the sublattice — these are the $T_{\lambda,z}$ with $\lambda = (\lambda_1, \lambda_2, 0)$. By $\operatorname{Aut}^+ I_{3,1}^{\mathcal{E}}$ we mean the subgroup of $\operatorname{Aut}^+ I_{4,1}^{\mathcal{E}}$ preserving $I_{3,1}^{\mathcal{E}}$. It is easy to see that this has index 2 in $\operatorname{Aut} I_{3,1}^{\mathcal{E}}$.

**(7.10) Lemma.** *The group $\langle R_1, R_2, R_4, R_5 \rangle$ contains all the translations of $I_{3,1}^{\mathcal{E}}$.*

(7.11) **Proof:** Computation reveals that $R_2^{-1} R_1 = T_{\omega,0,0;\theta/2}$. By (7.8.2), we also have $T_{-\bar{\omega},0,0;\theta/2} = R_1 T_{\omega,0,0;\theta/2} R_1^{-1}$. By using $R_5$ and $R_4$ in place of $R_1$ and $R_2$ we obtain the translations $T_{0,\omega,0;\theta/2}$ and $T_{0,-\bar{\omega},0;\theta/2}$. One may avoid the second computation by using (7.8.2) where $U$ exchanges the first two coordinates. By (7.8.1) we have

$$[T_{\omega,0,0;\theta/2}, T_{-\bar{\omega},0,0;\theta/2}] = T_{0,0,0;\theta}.$$

Use of (7.8.1) shows that these five translations generate the group of all translations of $I_{3,1}^{\mathcal{E}}$.

**(7.12) Lemma.** *The group $\langle R_1, \ldots, R_5 \rangle$ coincides with $\operatorname{Aut}^+ I_{3,1}^{\mathcal{E}}$ and acts transitively on the primitive null vectors of $I_{3,1}^{\mathcal{E}}$.*

(7.13) **Proof:** Write $\mathcal{R}$ for the group generated by $R_1, \ldots, R_5$. Because $\mathcal{R}$ contains all the translations of $I_{3,1}^{\mathcal{E}}$ and also the reflections in the short root $r_3$, it contains the reflections in all the short roots of height 1 in $I_{3,1}^{\mathcal{E}}$. According to the proof of Theorem 6.6 of [1], the group generated by the translations and these reflections is $\operatorname{Aut}^+ I_{3,1}^{\mathcal{E}}$. The arguments of [1] simplify in this case, so we give a fairly complete sketch the proof here.

(7.14) First we argue that $\mathcal{R}$ acts transitively on the primitive null lattices in $I_{3,1}^{\mathcal{E}}$. Given a null vector $v = (v_1, v_2, 0; H, ?)$ of height $H$, we may change $v_1$ and $v_2$ by any elements of $\mathcal{E}H$ by applying a translation. Since the covering radius of $\mathcal{E}$ is $1/\sqrt{3}$, we may suppose further that each of $|v_1|^2$ and $|v_2|^2$ is bounded by $|H|^2/3$. With $v$ prepared in this way, the $(-\omega)$-reflection in one of the roots $(0, 0, 0; 1, n\theta - \omega)$, where $n \in \mathbb{Z}$, carries $v$ to a vector of height smaller (in absolute value) than $H$. Repeating this procedure as needed, we may reduce $H$ until it is 0, so that $v$ is a multiple of $\rho$.

(7.15) Now let $J$ be the group generated by $\mathcal{R}$ and the biflection $B$ in $b = (1, -1, 0; 0, 0)$. Observe that $B$ normalizes $\mathcal{R}$, since it exchanges the left and right arms of (7.7.2). This implies that $\mathcal{R}$ has index at most two in $J$. We will show that $J = \operatorname{Aut} I_{3,1}^{\mathcal{E}}$. First we show that $J$ acts transitively on the primitive null vectors of $I^\mathcal{E}_{3,1}$. Computation reveals that $R_3$ acts on $I^\mathcal{E}_{1,1}$ by the matrix
$$\begin{pmatrix} \theta\bar\omega & \bar\omega \\ \bar\omega & 0 \end{pmatrix},$$
and one can check that $R_3 T_{0,0,0;-\theta}$ acts on $I^\mathcal{E}_{1,1}$ by $\begin{pmatrix} 0 & \bar\omega \\ \bar\omega & 0 \end{pmatrix}$. Call this transformation $F$ and observe that $F^2$ acts on $I^\mathcal{E}_{1,1}$ by the scalar $\omega$. Defining $B'$ as the biflection in $b' = (0,0,0;1,1)$, computation reveals that
$$T_{\omega,-\omega,0;0}\, F\, T_{1,0,0;\theta/2}(b) = -\bar\omega b'.$$

This proves that $B' \in J$. Since $B'$ acts on $I^\mathcal{E}_{1,1}$ by the matrix $\begin{pmatrix} 0 & -1 \\ -1 & 0 \end{pmatrix}$, the transformation $B'F$ acts on $I^\mathcal{E}_{1,1}$ by the scalar $-\bar\omega$. Consequently, $J$ acts transitively on the six unit scalar multiples of $\rho$ and hence on the primitive null vectors of $I^\mathcal{E}_{3,1}$.

(7.16) Since the stabilizer of $\rho$ in $\operatorname{Aut} I^\mathcal{E}_{3,1}$ is generated by the translations, the 6-fold reflections in $r_1$ and $r_5$, and the biflection $B$, we see that $J = \operatorname{Aut} I^\mathcal{E}_{3,1}$. Since $\mathcal{R}$ is generated by reflections in short roots, we have $\mathcal{R} \subset \operatorname{Aut}^+ I^\mathcal{E}_{3,1}$. Above we observed that $\mathcal{R}$ has index $\leq 2$ in $J$; therefore $\mathcal{R} = \operatorname{Aut}^+ I^\mathcal{E}_{3,1}$. This proves the first claim. To prove the second claim it suffices (because $\mathcal{R}$ contains $F^2$) to show that $\rho$ and $-\rho$ are equivalent under $\mathcal{R}$. They are equivalent since the product of $B$ and the central involution of $I^\mathcal{E}_{1,1}$ exchanges them and has spinor norm $+1$, and so lies in $\mathcal{R}$.

**(7.17) Lemma.** $\Gamma$ contains the $(-\omega)$-reflection $R_7$ in the vector $r_7 = (0,0,1;0,0)$; this extends the $E_6$ diagram consisting of the top 6 nodes in (7.7.2) to the full affine $\tilde E_6$ diagram pictured there.

(7.18) **Proof:** It suffices to show that $r_7$ is equivalent under $\Gamma$ to a multiple of $r_6$. We have $R_6(r_7) = (0,0,-\omega;0,\bar\omega)$. By applying an element of $\operatorname{Aut}^+ I^\mathcal{E}_{3,1} \subset \Gamma$ fixing $r_7$ and carrying the null vector $(0,0,0;0,\bar\omega)$ to $(0,0,0;0,-\omega)$, we carry $R_6(r_7)$ to $-\omega r_6$. This completes the proof.

(7.19) **Remark.** Looijenga [31] has obtained a remarkable presentation of the fundamental group of the space of smooth cubic surfaces as a quotient of the affine Artin group $\mathcal{A}(\tilde E_6)$. If one used this presentation, then one could deduce the lemma immediately, because the $(-\omega)$-reflection in $r_7$ is the only possibility for image of the seventh Artin generator. However, we will need both the result of Lemma 7.12 and its proof in Theorem 7.21.

**(7.20) Lemma.** $\langle R_1, R_2, R_4, \ldots, R_7 \rangle$ contains all the translations of $I^\mathcal{E}_{4,1}$.

**Proof:** By the argument of Lemma 7.10, with $R_7$ and $R_6$ in place of $R_1$ and $R_2$ one sees that $\langle R_6, R_7 \rangle$ contains $T_{0,0,\omega;\theta/2}$ and $T_{0,0,-\bar\omega;\theta/2}$. Together with the result of that lemma, this proves the claim.

**(7.21) Theorem.** $\Gamma$ coincides with $\operatorname{Aut}^+ I^\mathcal{E}_{4,1}$ and acts transitively on the primitive null vectors and the short roots of $I^\mathcal{E}_{4,1}$.

**Proof:** The argument is very similar to the proof of Lemma 7.12. As before, it follows from the proof of Theorem 6.6 of [1]; we sketch the argument here.

(7.22) First we show that $\Gamma$ acts transitively on the primitive null vectors of $I^{\mathcal{E}}_{4,1}$. Given a null vector $v = (v_1, v_2, v_3; H, ?)$ of $I^{\mathcal{E}}_{4,1}$ with $|H| > 0$, then as before we may suppose that each of $|v_1|^2$, $|v_2|^2$ and $|v_3|^2$ is bounded by $|H|^2/3$. With $v$ prepared in this manner there is a short root of the form $r = (0, 0, 0; 1, n\theta - \omega)$ with $n \in \mathbb{Z}$ such that either $(-\omega)$-reflection in $r$ reduces the height of $v$ or else $r$ is orthogonal to $v$. (The latter case occurs only if $|v_1|^2 = |v_2|^2 = |v_3|^2 = |H|^2/3$.) Thus, after applying some sequence of reflections in $\Gamma$, we may suppose that either $H = 0$ or else $v$ is orthogonal to some short root $r$ of height 1. In the former case, $v$ is proportional to $\rho$, hence $\Gamma$-equivalent to $\rho$ by Lemma 7.12. In the latter case, after applying a translation, we may take $r = r_3$. Because $R_3$ braids with $R_6$ and $R_6$ with $R_7$, we see that $v$ is $\Gamma$-equivalent to an element of $r_7^\perp = I^{\mathcal{E}}_{3,1}$. Thus, by Lemma 7.12, $v$ is $\Gamma$-equivalent to $\rho$.

(7.23) The proof that the group generated by $\Gamma$ and $B$ is all of $\operatorname{Aut} I^{\mathcal{E}}_{4,1}$ is essentially the same as the corresponding part of the proof of Lemma 7.12, and so is the argument for the conclusion $\Gamma = \operatorname{Aut}^+ I^{\mathcal{E}}_{4,1}$. The transitivity of $\operatorname{Aut} I^{\mathcal{E}}_{4,1}$ on short roots follows from the fact that the orthogonal complement of any short root is unimodular of signature $(3, 1)$, and hence by Lemma 8.1 of [1] a copy of $I^{\mathcal{E}}_{3,1}$. This implies transitivity of $\operatorname{Aut}^+ I^{\mathcal{E}}_{4,1}$ on pairs $\pm v$ of short roots. Since $v$ and $-v$ are exchanged by biflection in either we actually have transitivity on short roots.

(7.24) The monodromy group is now identified completely. There is no further need for $I^{\mathcal{E}}_{4,1}$, and we revert to the previous notation, in which we have $\Gamma = \operatorname{Aut}^+ \Lambda$. We now study the congruence subgroup $\Gamma_\theta$ — the kernel of the natural action of $\Gamma$ on $\Lambda/\theta\Lambda$. To this end, recall from (3.12) that the quotient of $\Gamma$ by $\Gamma_\theta$ is isomorphic to the Weyl group $W(E_6)$.

**(7.25) Theorem.** $\Gamma_\theta$ *is generated by the triflections in the short roots of $\Lambda$.*

**Proof:** Declaring that the standard generators have order 2 reduces $\mathcal{A}(E_6)$ to $W(E_6)$. Therefore the normal subgroup of $\Gamma$ generated by squares of the $R_i$ has index $\leq |W(E_6)|$. Since the squares lie in $\Gamma_\theta$ and $\Gamma/\Gamma_\theta \cong W(E_6)$, they must normally generate all of $\Gamma_\theta$. Therefore their conjugates, which are triflections in the short roots of $I^{\mathcal{E}}_{4,1}$, generate $\Gamma_\theta$.

**(7.26) Lemma.** *Suppose that $L$ is a free $\mathcal{E}$-module and that $g$ is a nontrivial element of $GL(L)$ which acts trivially on $L/\theta L$. If $g$ has finite order then its order is 3 and $L$ is the direct sum of its eigenlattices.*

**Proof:** We adapt an argument from §39 of [45]. Let $g'$ be a power of $g$ that has prime order $p$ and matrix expression $g' = I + \theta^k A$ with respect to some basis for $L$, where $k > 0$ and $A \not\equiv 0 \bmod \theta$. If $p \neq 3$, then $\theta$ does not divide $p$, and the relation

$$I = (g')^p = I + p\theta^k A + \binom{p}{2}\theta^{2k} A^2 + \cdots$$

yields $A \equiv 0$ modulo $\theta$, contrary to the choice of $A$. Therefore the order of $g$ is a power of 3. If $g$ has order $> 3$, then there is a power $g'$ of $g$ with order 9. In this case we write $L_1$, $L_\omega$ and $L_{\bar\omega}$ for the eigenlattices of $(g')^3$. Since $g'$ has order 9, we know that at least one of $L_\omega$ and $L_{\bar\omega}$ in nontrivial — say $L_\omega$. Writing $g''$ for the restriction of $g'$ to $L_\omega$, we have $g'' = I + \theta^k B$ with respect to some basis of $L_\omega$, where $k > 0$ and $B \not\equiv 0(\theta)$. Reducing the relation

$$\omega I = (g'')^3 = I + 3\theta^k B + 3\theta^{2k} B^2 + \theta^{3k} B^3$$

modulo $3 = -\theta^2$ yields the absurdity $\omega \equiv 1$ modulo 3. Therefore the order of $g$ must be 3.

It remains to show that the eigenlattices of $g$, say $L_1$, $L_\omega$ and $L_{\bar\omega}$, are summands of $L$. Write $\lambda = \lambda_1 + \lambda_\omega + \lambda_{\bar\omega}$ with $\lambda_1$, $\lambda_\omega$ and $\lambda_{\bar\omega}$ being the eigencomponents of $\lambda$ in $L \otimes_{\mathcal{E}} \mathbb{Q}(\sqrt{-3})$. We must show that $\lambda_1 \in L_1$, $\lambda_\omega \in L_\omega$ and $\lambda_{\bar\omega} \in L_{\bar\omega}$. This follows easily from the fact that the three vectors

$$\lambda = \lambda_1 + \lambda_\omega + \lambda_{\bar\omega}, \qquad g\lambda = \lambda_1 + \omega\lambda_\omega + \bar\omega\lambda_{\bar\omega}, \quad \text{and} \quad g^2\lambda = \lambda_1 + \bar\omega\lambda_\omega + \omega\lambda_{\bar\omega}$$

lie in $L$ and the hypothesis that the difference of any two of them lies in $\theta L$.

**(7.27) Lemma.** *Any torsion element of $\Gamma_\theta$ is a product of one of the scalars $1$, $\omega$, $\bar\omega$ with triflections in some number of orthogonal short roots. In particular, $P\Gamma_\theta$ acts freely on $\mathbb{C}H^4 - \mathcal{H}$, and $P\Gamma_\theta \backslash \mathbb{C}H^4$ is a complex manifold.*

**Proof:** By Theorem 7.26, if $g \in \Gamma_\theta$ has finite order then its order is 3 and $\Lambda$ is the direct sum of its eigenlattices. Since $\Lambda$ is unimodular, each summand must be unimodular, so each summand is isomorphic to $\mathcal{E}^n$ ($1 \le n \le 4$) or $\mathcal{E}^{n,1}$ ($0 \le n \le 4$). This uses the classification of unimodular $\mathcal{E}$-lattices in low dimensions — see [1] or [21]. Exactly one of the summands contains negative-norm vectors, and by multiplying $g$ by one of the scalars $1$, $\omega$ or $\bar\omega$ we may suppose that $g$ acts trivially on this summand. Then each of the $\omega$- and $\bar\omega$-eigenlattices has a basis of mutually orthogonal short roots, so that $g$ is the product of the $\omega$- and $\bar\omega$-reflections in these roots. Thus the stabilizer of a point in $\mathcal{H}^k$ is $(\mathbb{Z}/3)^k$ acting in the standard way, and the lemma follows.

The next sequence of results shows that the configuration of geodesic hyperplanes $\mathcal{H}$ is quite special. We begin with the following.

**(7.28) Lemma.** *If $H$ and $H'$ are distinct components of $\mathcal{H}$ and $H \cap H' \ne \emptyset$, then $H$ and $H'$ are orthogonal.*

**Proof:** Let $H$ and $H'$ be defined by the equations $h(x,v) = 0$ and $h(x,v') = 0$, respectively, where $h(v,v) = h(v',v') = 1$. Since they meet in $\mathbb{C}H^4$, there is a vector of negative norm orthogonal to both $v$ and $v'$. Since $h$ has signature $(4,1)$, $v$ and $v'$ span a positive-definite space, which implies that their inner product matrix

$$\begin{pmatrix} 1 & h(v,v') \\ h(v',v) & 1 \end{pmatrix}$$

is positive definite. That is, $|h(v,v')|^2 < 1$. But this is an integer, so $h(v,v') = 0$ as claimed.

(7.29) The set $\mathcal{R}$ of short roots given in (2.19) projects, under reduction modulo $\theta$, to the set $\bar{\mathcal{R}} \subset V$ of vectors of norm one, where $V$ is the finite vector space of (2.12), which in turn projects onto $P\bar{\mathcal{R}} \subset PV$. It is easy to check that $\bar{\mathcal{R}}$ has 72 elements, or equivalently, that $P\bar{\mathcal{R}}$ has 36 elements (these are the 36 "minus points" of [17], p. 26). Their preimages in $\mathcal{R}$ gives a partition $\mathcal{R} = \mathcal{R}^1 \cup \cdots \cup \mathcal{R}^{36}$ into 36 sub-collections which are permuted by $\Gamma$. Each sub-collection is preserved and acted on transitively by $\Gamma_\theta$. There is a corresponding partition,

$$\mathcal{H} = \mathcal{H}_1 \cup \cdots \cup \mathcal{H}_{36}, \tag{7.29.1}$$

into 36 disjoint sub-collections permuted by $\Gamma$ and each acted on transitively by $\Gamma_\theta$.

**(7.30) Lemma.** *Let $H$ and $H'$ be components of $\mathcal{H}$ which belong to the same sub-collection of the decomposition (7.29.1), and suppose that $H \neq H'$. Then $H \cap H' = \emptyset$.*

**Proof:** We write $H = v^\perp$ and $H' = v'^\perp$ where $v$ and $v'$ are congruent modulo $\theta$. Then, modulo $\theta$, $h(v, v') \equiv h(v, v) \equiv 1$. Thus $h(v, v') \neq 0$, so by Lemma 7.28, $H \cap H' = \emptyset$.

**(7.31) Lemma.** *Let $H \in \mathcal{H}$ and let $\Gamma_H \subset \Gamma_\theta$ denote its stabilizer in $\Gamma_\theta$. Then the natural map $P\Gamma_H \backslash H \longrightarrow P\Gamma_\theta \backslash \mathbb{C}H^4$ (induced by the injection $H \subset \mathbb{C}H^4$) is injective.*

**Proof:** Let $x, y \in H$ and suppose that there exists $\gamma \in P\Gamma_\theta$ such that $x = \gamma y$. Then $H \cap \gamma H \neq \emptyset$. Since $H$ and $\gamma H$ belong to the same sub-collection in the decomposition (7.29.1), Lemma 7.30 implies that $H = \gamma H$, so that $\gamma \in P\Gamma_H$. The lemma follows.

**(7.32) Theorem.** *The space $P\Gamma_\theta \backslash \mathbb{C}H^4$ is smooth. The projection of $\mathcal{H}$ to $P\Gamma_\theta \backslash \mathbb{C}H^4$ consists of 36 smooth divisors intersecting transversally.*

**Proof:** It is clear from Lemmas 7.27 and 7.28 that the isotropy groups of the action of $P\Gamma_\theta$ on $\mathbb{C}H^4$ are isomorphic to $(\mathbb{Z}/3\mathbb{Z})^k$, $k \leq 4$, acting in the standard way. Thus $P\Gamma_\theta \backslash \mathbb{C}H^4$ is a smooth manifold. By the same token, for each component $H$ of $\mathcal{H}$, $\Gamma_H \backslash H$ is smooth. By Lemma 7.31 each component of $P\Gamma_\theta \backslash \mathcal{H}$ is embedded in $\Gamma_\theta \backslash \mathbb{C}H^4$, and by (7.28), any two components intersect at right angles, hence transversally.

Finally, we record the injectivity lemma analogous to Lemma 7.31 for $P\Gamma$. This will be used in §9

**(7.33) Lemma.** *Let $H$ be a component of $\mathcal{H}$ and let $\Gamma_H \subset \Gamma$ denote its stabilizer in $\Gamma$. Then the natural map $P\Gamma_H \backslash H \longrightarrow P\Gamma \backslash \mathbb{C}H^4$ (induced by the injection $H \subset \mathbb{C}H^4$) is injective.*

**Proof:** . Suppose $h_1, h_2 \in H$, $\gamma \in P\Gamma$, and $\gamma h_1 = h_2$. We produce $\gamma' \in P\Gamma_H$ such that $\gamma' h_2 = h_2$. Let $v$ be a short root orthogonal to $H$. Then $v$ and $\gamma v$ are orthogonal short roots, and so $v + \gamma v$ is a long root. Let $A$ be the biflection in this long root and let $\gamma_0 = PA\ (= P(-A))$. Then $\gamma_0 \in P\Gamma$, $\gamma_0 H = \gamma H$ and $\gamma_0$ is the identity on $H \cap \gamma H$, thus $\gamma_0^{-1} \gamma \in \Gamma_H$ gives the desired element $\gamma'$, and so proves the lemma.

## 8. Cuspidal degenerations

Our purpose here is to extend the period map $g : \mathcal{C}_s \to P\Gamma \backslash \mathbb{C}H^4$ to a map $\mathcal{C}_{ss} \to \overline{P\Gamma \backslash \mathbb{C}H^4}$, where $\overline{P\Gamma \backslash \mathbb{C}H^4}$ denotes the Satake compactification of $P\Gamma \backslash \mathbb{C}H^4$. For this we need some information on the monodromy of $\Lambda(\mathcal{T})$ near points $F \in \Delta_{ss} - \Delta_s$. Instead of doing a detailed analysis of topology of $\mathcal{T}$ near a point of $\Delta_{ss} - \Delta_s$ analogous to the discussion of §5 for points of $\Delta_s$, we find it expedient to use the results of that section and §7 to deduce the needed information on the local monodromy. We first prove that the monodromy transformation corresponding to the loop $\gamma_{\ell,r}$ defined below is a central translation (as defined in (7.8)), hence it has a unique fixed point in $\mathbb{C}H^4 \cup \partial\mathbb{C}H^4$, which lies in $\partial\mathbb{C}H^4$ and represents a boundary point of the Satake compactification. A simple classical argument based on the Schwarz lemma and the fact that the hyperbolic length of the loops $\gamma_{\ell,r}$ tend to 0 as $r$ tends to 0 shows that the period map takes neighborhoods of $F$ to neighborhoods of this fixed point. This observation easily yields the desired extension.

Let $F \in \Delta_{ss}^{a;b}$, with $b \geq 1$ cusps. By (3.8), there is a neighborhood $W$ of $F$ in $\mathcal{C}_{ss}$ and a Galois covering space of $W - \Delta$ such that the Fox completion $\tilde{W}$ of this space over $W$ is isomorphic to the unit ball in $\mathbb{C}^{20}$. Moreover, under that isomorphism the preimage $\tilde{\Delta}$ of $\Delta$ maps to the union of the hyperplanes for a root system of type $A_1^a \oplus A_2^b$, and the Galois group of $\tilde{W} \to W$ to the group $(\mathbb{Z}/2)^a \times S_3^b$ generated by the reflections in these hyperplanes. Furthermore, there is a unique point of $\tilde{W}$ lying over $F$, namely the origin.

Let $\ell$ be a line in $\mathbb{C}^{20}$ that intersects $\tilde{\Delta}$ only at the origin. For $0 < r < 1$ we write $D_{\ell,r}$ for the intersection of $\ell$ with the closed $r$-ball in $\mathbb{C}^{20}$, and $\gamma_{\ell,r}$ for the (positively oriented) boundary curve of $D_{\ell,r}$. These loops are all freely homotopic. Choose a form $F' \in \partial W$ and a point $m'$ of $\tilde{W}$ that lies over $F'$. Consider the loop $\gamma_{\ell,r}$ passing through $m'$, viewed as based at $m'$, and write $\bar{\gamma}$ for the element of $\pi_1(\bar{W} - \Delta, F') \cong \pi_1(W - \Delta, F')$ represented by the projection of $\gamma_{\ell,r}$.

**(8.1) Lemma.** *The monodromy action of $\bar{\gamma}^3$ on $\Lambda(T')$ is a central translation (as defined in (7.8)).*

**Proof:** Recall from (3.7) that $\pi_1(W - \Delta) \cong \mathbb{Z}^a \times B_3^b$. We will identify the action of $\pi_1(W - \Delta)$ on $\Lambda(T')$ in terms the actions of the standard generators for the various $\mathbb{Z}$ and $B_3$ factors, and then compute the action of $\bar{\gamma}$ by expressing $\bar{\gamma}$ in terms of these generators. By (5.4), each of the standard generators of $\mathbb{Z}^a \times B_3^b$ maps to the $(-\omega)$-reflection in a short root of $\Lambda(T')$. For each $B_3$ factor, its standard generators $g_1$ and $g_2$ act by the $(-\omega)$-reflections in short roots $s_1$ and $s_2$. Since $g_1$ and $g_2$ braid, we may assume that $h(s_1, s_2) = 1$ (see (7.4) and (7.5).) Observe the following two facts:

(i.) For each $B_3$ factor, the corresponding short roots $s_1$ and $s_2$ are linearly independent.

(ii.) The roots corresponding to distinct factors of $\mathbb{Z} \times \cdots \times \mathbb{Z} \times B_3 \times \cdots \times B_3$ are orthogonal.

The first fact follows from the fact that the corresponding reflections, reduced modulo $\theta$, generate a subgroup $S_3$ of the Weyl group (see (3.12)). If $s_1$ and $s_2$ were dependent then the reflections in them would generate an abelian group. The second fact follows from the fact that the geometric vanishing cycles in $T'$ corresponding to distinct singularities in $T$ are disjoint.

These facts determine the local monodromy group. If $s_1$ and $s_2$ are the roots for a $B_3$ factor, it follows from fact (i) that $s_1 - s_2$ is a primitive null vector. From fact (ii), the $b$ null vectors obtained in this way are mutually orthogonal. Since the signature of $\Lambda(T')$ is $(4,1)$, a maximal isotropic sublattice has dimension 1, these $b$ vectors must all be proportional, hence may be assumed to coincide (after multiplying the $s_i$ by suitable units). We call this common vector $\rho$. By the transitivity of $\mathrm{Aut}(\Lambda)$ on the primitive null vectors of $\Lambda$ (Theorem 7.21), we may choose a basis for $\Lambda(T')$ with inner product matrix as in (7.7.1), such that $\rho = (0,0,0,0,1)$. This is consistent with the use of the symbol $\rho$ in (7.8). Taking the root for each $\mathbb{Z}$ factor and the root $s_1$ for each $B_3$ factor, we obtain up to three mutually orthogonal short roots, which are also orthogonal to $\rho$. By the choice of the basis of $\Lambda(T')$ we may take these roots to be the first $a+b$ of $(1,0,0,0,0)$, $(0,1,0,0,0)$ and $(0,0,1,0,0)$. The remaining root for each $B_3$ factor is determined by the relation $s_1 - s_2 = \rho$. This implies that the hexflections to which the standard generators of $\mathbb{Z}^a \times B_3^b$ map are among the the $(-\omega)$-reflections in the roots $r_1, r_2, r_4, r_5, r_6$ and $r_7$ that form the 3 "limbs" of the affine $\tilde{E}_6$ diagram (7.7.2). More specifically, if $b=1$ then the roots may be taken to be $r_1$ and $r_2$, together with $r_5$ (resp. $r_5$ and $r_6$) if $F$ has one node (resp. two nodes). Similarly, if $b=2$ then the roots may be taken to be $r_1, r_2, r_4, r_5$, together with $r_7$ if $F$ has a node, and if $b=3$ then the roots are $r_1, r_2, r_4, r_5, r_6$ and $r_7$.

To express $\bar{\gamma}$ in terms of the generators used above, one recognizes $\bar{\gamma}$ as the square of the "fundamental element" of the Artin group

$$\mathcal{A}(A_1^a A_2^b) = \mathcal{A}(A_1)^a \times \mathcal{A}(A_2)^b = \mathbb{Z}^a \times B_3^b \ .$$

To explain this, we follow a suggestion of Ruth Charney. Since scalar multiplication on $\mathbb{C}^{20}$ commutes with the action of the Weyl group, and since all the points of $\gamma_{\ell,r}$ are scalar multiplies of each other, it is clear that $\bar{\gamma}$ lies in the center of the Artin group. The fundamenal element is a certain product of $l$ standard Artin generators, and has the property that either it or its square generates the center. (See [9], secs. 5 and 7.) It follows that $\bar{\gamma}$ is some power $n$ of the fundamental element. Furthermore, the number $l$ coincides with the number of mirrors. On the other hand, it is easy to see that $\gamma_{\ell,r}$ is homologous to a sum of positively-oriented circles, one for each mirror. (Perturb $\ell$ so that it meets the mirrors at $l$ distinct points.) Since the map $\tilde{W} \to W$ has branching of order 2 along the generic points of each mirror, we see that $\bar{\gamma}$ is homologous to the sum of $2l$ Artin generators. By considering the map to $\mathbb{Z}$ that carries each generator to 1, we see that the only possibility for $n$ is 2.

The fundamental element is the product of the fundamental elements of the various $\mathbb{Z}$ and $B_3$ factors. The fundamental element of a $\mathbb{Z}$ factor is the standard generator. The fundamental element of a $B_3$ factor with standard generators $g_1$ and $g_2$ is $g_1 g_2 g_1 = g_2 g_1 g_2$. Now we can compute the monodromy of $\bar{\gamma}$. The square of the fundamental element of one $B_3$ factor, say the one where $g_1$ and $g_2$ act as the $(-\omega)$-reflections $R_1$ and $R_2$ in $r_1$ and $r_2$, acts as $(R_1 R_2 R_1)^2 = (R_1 R_2 R_1)(R_2 R_1 R_2) = (R_1 R_2)^3 = T_{0,-\theta}$. Here the last term is the central translation defined in (7.8), and the last equality is obtained by computing matrices for $R_1$ and $R_2$ and multiplying them together. If there is another $B_3$ factor, then similar calculations (or the $S_3$ symmetry of the $\tilde{E}_6$ diagram) show that the square of its fundamental element also acts as $T_{0,-\theta}$. ) The squares of the fundamental elements of the $\mathbb{Z}$ factors act as triflections in the relevant roots. If no nodes are present then by multiplying together the fundamental elements of the $B_3$ factors, we see that $\bar{\gamma}$ acts as $T_{0,-b\theta}$. If at least one node is present then $\bar{\gamma}$ acts by the product of $T_{0,-b\theta}$ and either or both of the $\bar{\omega}$-reflections in $r_5$ and $r_7$. In any case, $\bar{\gamma}^3$ acts as $T_{0,-3b\theta}$, which is a central translation as claimed. The proof of the lemma is complete.

**(8.2) Theorem.** *The map $g : \mathcal{C}_s \to P\Gamma\backslash\mathbb{C}H^4$ extends to a map from $\mathcal{C}_{ss}$ to the Satake compactification $\overline{P\Gamma\backslash\mathbb{C}H^4}$. This carries $\mathcal{C}_{ss} - \mathcal{C}_s$ to the (unique) boundary point of $\overline{P\Gamma\backslash\mathbb{C}H^4}$ and factors as a composition of holomorphic maps*

$$\mathcal{C}_{ss} \to M \to \overline{P\Gamma\backslash\mathbb{C}H^4}$$

*which takes the unique point of $M - M_s$ to the unique boundary point of $\overline{P\Gamma\backslash\mathbb{C}H^4}$.*

**Proof:** For the sake of brevity we write $Q$ for $P\Gamma\backslash\mathbb{C}H^4$ and $\bar{Q}$ for $\overline{P\Gamma\backslash\mathbb{C}H^4}$. By the transitivity of $\Gamma$ on the one-dimensional isotropic sublattices of $\Lambda$ (see (7.21)), the Satake compactification $\bar{Q} = \overline{P\Gamma\backslash\mathbb{C}H^4}$ has a unique boundary point, which we will call $B$. The theorem is a consequence of the following claim: if $F \in \Delta_{ss}^{a,b}$, $b \geq 1$, and $\tilde{W}$ is as above, then for any neighborhood $V$ of $B$ in $\bar{Q}$, there is a neighborhood $\tilde{U}$ of the origin in $\tilde{W}$ so that, if $U \subset W$ denotes its image under $\tilde{W} \to W$, then $g(U) \subset V$. Assuming this for the moment, since $\bar{Q}$ is an analytic space, the Riemann extension theorem implies that the period map $g$ extends holomorphically to this neighborhood of $U$ of $F$. By shrinking $V$ we see that this extension takes $F$ to $B$. Thus $g : \mathcal{C}_s \to Q$ extends to holomorphic map $\mathcal{C}_{ss} \to \bar{Q}$. Since it carries $\mathcal{C}_{ss} - \mathcal{C}_s$ to $B$, and since the period map is constant on $G$-orbits in $\mathcal{C}_0$, it is also constant on orbits in $\mathcal{C}_{ss}$. This establishes the theorem.

Now we prove the claim. Let $\ell$ be any line through the origin in $\mathbb{C}^{20}$ that meets $\tilde{\Delta}$ only at the origin as before, and note that $\ell \cap \tilde{W}$ is isomorphic to the unit disk in $\mathbb{C}$. Thus for each such $\ell$ we may choose a biholomorphism $\alpha_\ell$ from the usual upper half-plane $H$, modulo the horizontal translation $\tau : z \mapsto z + 1$, to the punctured disk $(\ell \cap \tilde{W}) - \{0\}$. Note that the hyperbolic length of $\gamma_{\ell,r}$ (with respect to this uniformization) is *independent of $\ell$*, and decreases monotonically to 0 as $r \to 0$. Following $\alpha_\ell$ by the map $\tilde{W} - \tilde{\Delta} \to W - \Delta$, we obtain a map $\beta_\ell : H/\langle\tau\rangle \to \mathcal{C}_0$. Fix a deck transformation $\gamma : \mathcal{F}_0 \to \mathcal{F}_0$ in the conjugacy class in $\pi_1(\mathcal{C}_0)$ of the projection of $\gamma_{\ell,r}$ to $W - \Delta \subset \mathcal{C}_0$. Then there is a lift $\tilde{\beta}_\ell : H \to \mathcal{F}_0$ of $\beta_\ell$ under which the action of $\gamma$ on $\mathcal{F}_0$ corresponds to that of $\tau$ on $H$. Applying the period map, the action of $\gamma$ on $\mathcal{F}_0$ descends to an action on $\mathbb{C}H^4$, and $g \circ \tilde{\beta}_\ell : H \to \mathbb{C}H^4$ is equivariant with respect to these actions.

Since $\gamma^3$ acts on $\mathbb{C}H^4$ by a unitary transvection (Lemma 8.1), it fixes a unique one-dimensional isotropic sublattice of $\Lambda$, corresponding to an ideal point $\tilde{B}$ of $\mathbb{C}H^4$. Furthermore, the sets

$$N_\varepsilon = \{\, p \in \mathbb{C}H^4 \mid d_{\mathbb{C}H^4}(p, \gamma^3(p)) < \varepsilon \,\}$$

are the family of open horoballs centered at $\tilde{B}$, where we write $d$ for the hyperbolic metric. By the definition of the Satake compactification, the images of the $N_\varepsilon$ in $Q$ form a basis of deleted neighborhoods of $B$ in $\bar{Q}$. Given $\varepsilon > 0$, choose $r_0$ small enough so that the hyperbolic length of each $\gamma_{\ell,r_0}$ is less than $\varepsilon/3$, and let $\tilde{U}$ be the open $r_0$-ball in $\mathbb{C}^{20}$. If $x \in \tilde{U} - \tilde{\Delta}$ then we find the loop $\gamma_{\ell,r}$ containing $x$ and choose a point $y$ in $H$ that is carried to $x$ under the composition of the natural map $H \to H/\langle\tau\rangle$ and $\alpha_\ell$. Since the hyperbolic length of $\gamma_{\ell,r}$ is less than $\varepsilon/3$, we see that $d_H(y, \tau^3(y)) < \varepsilon$. By the Schwarz lemma $g \circ \tilde{\beta}_\ell$ does not increase hyperbolic distance. Therefore $g \circ \tilde{\beta}_\ell(y) \in N_\varepsilon$. Consequently the image of $x$ in $Q$ lies in the image of $N_\varepsilon$ in $Q$, and $\tilde{U}$ is the desired neighborhood. This completes the proof of the theorem.

## 9. Proof of the main theorem

The aim of this section is to prove our main results, Theorems 3.17 and 3.20. We will follow the strategy described in (3.22).

**Isomorphism on $M_0$**

**(9.1) Lemma.** *The differential of the period map $g : M_0^f \longrightarrow \mathbb{C}H^4$ is injective. Therefore $g : M_0^f \to \mathbb{C}H^4$ is locally biholomorphic and $g : M_0 \to P\Gamma\backslash\mathbb{C}H^4$ is an open map.*

**Proof:** Consider $Pf \in \mathcal{F}_0$ lying over $F \in \mathcal{C}_0$ and a family $F_t = F + tG$ of cubic forms. Recall that the differential form
$$\phi_t = \frac{\Omega}{(F+tG)^{4/3}}$$
spans $H_{\bar{\omega}}^{2,1}(T_t)$. Note that the corresponding expression (6.7.1) for $g$ can be differentiated under the integral sign so that the derivative of the period map is given by the cohomology class of
$$\left.\frac{d}{dt}\phi_t\right|_{t=0} = \text{const.}\frac{G\Omega}{F^{7/3}}. \tag{9.1.1}$$

This vanishes as an element of $T_{g(x)}\mathbb{C}H^4$ if and only if it lies in $H_{\bar{\omega}}^{2,1}(T)$. From Theorem 6.5 this is the case if and only if this expression is cohomologous to a fractional differential with filtration level one, i.e., with denominator $F^{4/3}$. By standard arguments as in [24], this is the case if and only if the numerator $G$ lies in the Jacobian ideal of $F$, that is, if and only if it is a tangent vector at $F$ to the orbit of the action of the general linear group. Since the tangent space to $M_0^f$ at $Pf$ is canonically isomorphic to the tangent space to $\mathcal{C}_0$ at $F$ modulo the tangent space to the orbit of $F$, this statement is equivalent to the injectivity of the differential of $g$ on $M_0^f$, thus proving the first assertion of the lemma. The second assertion is immediate, and the last assertion is clear from the $P\Gamma$-equivariance of $g$ and the definition of the quotient topology.

**(9.2) Lemma.** *The period map $g : M_0 \to P\Gamma\backslash\mathbb{C}H^4$ is injective, and its image is contained in $P\Gamma\backslash(\mathbb{C}H^4 - \mathcal{H})$.*

**Proof:** Let $N_0$ denote the moduli space of smooth cubic threefolds (a 10-dimensional space) and let $\mathfrak{S}^5$ denote the Siegel upper half space of genus five (a 15-dimensional space). The integral symplectic group $\text{Sp}(10, \mathbb{Z})$ acts properly discontinuously on $\mathfrak{S}^5$ and the Hodge structures of cubic threefolds are classified by $P\text{Sp}(10,\mathbb{Z})\backslash\mathfrak{S}^5$. In other words, the period map for cubic threefolds defines a holomorphic map $N_0 \longrightarrow P\text{Sp}(10,\mathbb{Z})\backslash\mathfrak{S}^5$. It fits, with the period map already defined, into a commutative diagram

$$\begin{array}{ccc} M_0 & \xrightarrow{g} & P\Gamma\backslash\mathbb{C}H^4 \\ \downarrow & & \downarrow \\ N_0 & \longrightarrow & P\text{Sp}(10,\mathbb{Z})\backslash\mathfrak{S}^5 \end{array}. \tag{9.2.1}$$

Indeed, there is a totally geodesic embedding $\mathbb{C}H^4 \longrightarrow \mathfrak{S}^5$ which can be described as follows: Let $\Omega$ be the standard unimodular symplectic form on $\mathbb{Z}^{10}$ and recall that $\mathfrak{S}^5$ is the space of Lagrangian subspaces in $\mathbb{C}^{10} = \mathbb{Z}^{10} \otimes \mathbb{C}$ on which the hermitian form $\Omega(X, \bar{Y})$ is positive definite. If $\mathbb{Z}^{10} = \mathcal{E}^5$ as in (2.2)–(2.4), decompose $\mathbb{C}^{10}$ into eigenspaces: $\mathbb{C}^{10} = \mathbb{C}^5_\omega \oplus \mathbb{C}^5_{\bar{\omega}}$. Then the embedding of $\mathbb{C}H^4$ into $\mathfrak{S}^5$ in question is the map that assigns to a negative line $l \in \mathbb{C}^{4,1}$ the Lagrangian subspace $\bar{l}^\perp \oplus l \subset \mathbb{C}^{10}$. From (2.4.1) we see that the diagram (9.2.1) is indeed commutative. Now the theorem of Clemens and Griffiths [15] asserts that the lower horizontal arrow of (9.2.1) is injective. The next lemma proves that the left vertical arrow is generically injective. Thus the top horizontal arrow is generically injective. By Lemma 9.1, the top horizontal arrow is an open map. Since a generically injective open map is injective, the first assertion follows.

The second assertion is equivalent to the statement that the image of $g : M^f_0 \to \mathbb{C}H^4$ is contained in $\mathbb{C}H^4 - \mathcal{H}$. To prove this equivalent assertion, suppose to the contrary that $g(x) \in \mathcal{H}$ for some $x \in M^f_0$. Then there is a short root $r \in \mathcal{R}$ so that the line $g(x)$ in $\mathbb{C}^{4,1}$ lies in $r^\perp$. This implies that the Lagrangian subspace $g(x)^\perp \oplus g(x) \subset \mathbb{C}^{10}$ contains the short root $r$. This in turn means that the intermediate Jacobian of the cubic threefold associated to $x$ has the elliptic curve $\mathbb{C}r/\mathcal{E}r$ as a factor. But the intermediate Jacobian of a smooth cubic threefold is an irreducible principally polarized abelian variety [15], so this is impossible, thus completing the proof of the lemma.

**(9.3) Lemma.** *Let $T$ and $T'$ be generic cyclic cubic threefolds with branch loci $S$ and $S'$. If $T$ and $T'$ are isomorphic via an automorphism of projective space, then so are $S$ and $S'$.*

**Proof:** If $T$ and $T'$ are generic then each admits just one cyclic structure, and a cyclic structure determines is branch set.

**(9.4) Lemma.** *The extended period map $g : M \to \overline{P\Gamma \backslash \mathbb{C}H^4}$ is surjective. Equivalently, the period map $g : M_s \to P\Gamma \backslash \mathbb{C}H^4$ is surjective and proper.*

**Proof:** In Theorem 8.2 we proved that the extended period map $g : M \longrightarrow \overline{P\Gamma \backslash \mathbb{C}H^4}$ takes the unique strictly semi-stable point of $M$ to the unique cusp of $\overline{P\Gamma \backslash \mathbb{C}H^4}$. This extension is a holomorphic map of analytic spaces of the same dimension which, by Lemma 9.1 is generically of maximal rank, hence of non-zero degree, thus proving the first statement of the lemma. The equivalence of the second statement is clear.

Combining Lemmas 9.2 and 9.4 with the fact that $g$ carries $\mathcal{F}_s - \mathcal{F}_0$ to $\mathcal{H}$, we obtain at once the following theorem, which is easily seen to be equivalent to the main theorem of [3]:

**(9.5) Theorem.** *The period map $g : M_0 \to P\Gamma \backslash (\mathbb{C}H^4 - \mathcal{H})$ is an isomorphism of analytic spaces.*

### Homeomorphism on the nodal divisor

Let $M_{nod} = M_s - M_0$, which is the same as $G \backslash \Delta_s$, the moduli space of cubic surfaces with nodes but no other singularities. Then $g|_{M_{nod}} : M_{nod} \longrightarrow P\Gamma \backslash \mathcal{H}$. Let $\mathcal{E}^{3,1} \subset \mathcal{E}^{4,1}$ be the standard embedding where the last coordinate is zero, let $\Gamma' = \mathrm{Aut}^+(\mathcal{E}^{3,1})$ be the group of spinor norm one automorphisms, and let $\Gamma' \subset \Gamma$ and $\mathbb{C}H^3 \subset \mathbb{C}H^4$ be the consequent embeddings. The natural map $P\Gamma'\backslash\mathbb{C}H^3 \to P\Gamma\backslash\mathcal{H}$ is bijective: surjectivity follows from transitivity of $P\Gamma$ on shorts roots (Theorem 7.21) and injectivity follows from Lemma 7.33. Since it is proper, this map is a homeomorphism. (It is also an isomorphism of analytic spaces, as are the other homeomorphisms in this subsection. We do not, however, need this for the proof. The analytic space isomorphism will also be a conseqence of the main theorem.) The aim of this subsection is to establish the following result:

**(9.6) Proposition.** *The map $g|_{M_{nod}} : M_{nod} \longrightarrow P\Gamma\backslash\mathcal{H} = P\Gamma'\backslash\mathbb{C}H^3$ is a homeomorphism.*

(9.7) For the proof we will use the fact that $M_{nod}$ is homeomorphic to the moduli space of stable sextuples of points in $\mathbb{P}^1$, and the fact, established in the work of Deligne and Mostow [20], that the latter space can be identified with $P\Gamma'\backslash\mathbb{C}H^3$. To make a more precise statement, let $(\mathbb{P}^1)_0^{(6)}$ denote the space of unordered sextuples of distinct points in $\mathbb{P}^1$, and let $(\mathbb{P}^1)_s^{(6)}$ denote the space of unordered sextuples in $\mathbb{P}^1$ with at most double points. The group $PGL(2)$ of automorphisms of $\mathbb{P}^1$ acts properly on both these spaces, and we let $M(6)_0 = PGL(2)\backslash(\mathbb{P}^1)_0^{(6)}$ and $M(6)_s = PGL(2)\backslash(\mathbb{P}^1)_s^{(6)}$ denote the corresponding quotients. Restating Lemma 2 of [11], we define a homeomorphism $\phi : M_{nod} \longrightarrow M(6)_s$ in the following way. Let $\bar\Delta_s$ denote the set of those $F \in \Delta_s$ such that $(0:0:0:1)$ is a node of $S$ and such that the tangent cone to $S$ at $(0:0:0:1)$ is $X_0X_2 - X_1^2 = 0$. Let $\bar G$ denote the subgroup of $G$ that leaving this point and cone invariant. It is easy to see that the natural map $\bar G\backslash\bar\Delta_s \to G\backslash\Delta$ is a bijection; this relies on the fact that the symmetry group of a cubic surface acts transitively on the nodes of the surface. (See the remark on p. 249 of [11].) Since it is proper, this map is a homeomorphism. Define $\bar\phi : \bar\Delta_s \to (\mathbb{P}^1)_s^{(6)}$ as follows. If $F \in \bar\Delta_s$, then (after possibly replacing $F$ by a constant multiple) we have

$$F(X_0, X_1, X_2, X_3) = X_3(X_0X_2 - X_1^2) + F_3(X_0, X_1, X_2) \tag{9.7.1}$$

for some cubic form $F_3$ in three variables. Then $\bar\phi(F)$ is the sextuple of zeros of the binary sextic form $F_3(U^2, UV, V^2)$. Geometrically, $\bar\phi(F)$ is the set of directions of the six lines (counting multiplicities appropriately) on $S$ through the node $(0:0:0:1)$. It follows from Lemma 2 of [11] that $\bar\phi(F)$ is stable, and that $\bar\phi$ descends to a homeomorphism between $\bar G\backslash\bar\Delta_s$ and $M(6)_s$. Defining $\phi : M_{nod} \to M(6)_s$ by precomposing with the identification $M_{nod} \cong \bar G\backslash\bar\Delta_s$, we have proved the following lemma:

**(9.8) Lemma.** *The map $\phi : M_{nod} \longrightarrow M(6)_s$ just defined is a homeomorphism.*

(9.9) In [20], Deligne and Mostow construct a number of period maps from suitable compactifications of spaces of distinct $N$-tuples of points in $\mathbb{P}^1$ to quotients of complex hyperbolic $(N-3)$-space by suitable discrete groups. Among their constructions is a period map $g' : M(6)_s \longrightarrow P\Gamma'\backslash\mathbb{C}H^3$, where $M(6)_s$ and $\Gamma'$ are as above. This is the first example in the list $N = 6$ of §14.4 of [20]. Some comments on the precise form in which we use their results are in order. First, we could start from sextic forms in 2 variables rather than the space $\mathcal{C}$ of (2.1) and perform all the constructions analogous to the ones of our §§2–3, using the first cohomology of the triple cover of $\mathbb{P}^1$ branched over a sextuple $x = \langle x_1, \cdots, x_6\rangle \in (\mathbb{P}^1)_0^{(6)}$. This is a non-singular curve $C_x$ of genus

4, with affine equation $w^3 = (z - x_1) \cdots (z - x_6)$. We obtain a period map

$$g' : M(6)_0 \longrightarrow P\Gamma' \backslash \mathbb{C}H^3$$

by assigning to $x$ the space $H_{\bar{\omega}}^{1,0}(C) \subset H_{\bar{\omega}}^1(C)$ spanned by the differential

$$\frac{dz}{w} = \frac{dz}{(z - x_1)^{1/3} \cdots (z - x_6)^{1/3}} \ .$$

Proceeding as before, this map extends to Fox completions and fits into a diagram analogous to (3.16.1):

$$\begin{array}{ccc} M(6)_s^f & \longrightarrow & \mathbb{C}H^3 \\ \downarrow & & \downarrow \\ M(6)_s^m & \longrightarrow & P\Gamma'_\theta \backslash \mathbb{C}H^3 \\ \downarrow & & \downarrow \\ M(6)_s & \longrightarrow & P\Gamma' \backslash \mathbb{C}H^3 \end{array} \qquad (9.9.1)$$

Here $M(6)_s^m$ is the moduli space of *ordered* stable sextuples in $\mathbb{P}^1$, and the top and middle horizontal arrows are shown in [20] to be isomorphisms. The lower horizontal arrow is obtained by passing to unordered sextuples, and it follows that it must also be an isomorphism. Since the actual groups are not identified in [20], we make the following remark. It is shown in [20] that the discrete group of hyperbolic motions in the right middle position of (9.9.1) is a homomorphic image of the pure braid group on six strands, with the standard generators acting by triflections. From this, or directly from considerations as in our §5, it follows that the discrete group occurring in the lower right hand corner of (9.9.1) is an image of the braid group on six strands and is generated by hexflections. From Lemma 7.12 it follows that this group is $P\Gamma'$ as asserted (and, as in Theorem 7.25, that the group in the middle right hand position is the congruence subgroup $P\Gamma'_\theta$). We thus have the following theorem:

**(9.10) Theorem (Deligne-Mostow [20]).** *The period map $g'$ just described gives an isomorphism $g' : M(6)_s \longrightarrow P\Gamma' \backslash \mathbb{C}H^3$ of analytic spaces.*

The proof of Proposition 9.6 reduces to the following proposition:

**(9.11) Lemma.** $g' \circ \phi = g|_{M_{nod}} : M_{nod} \longrightarrow P\Gamma' \backslash \mathbb{C}H^3$.

**Proof:** It is enough to check that $g' \circ \phi = g|_{M_{nod}}$ on the open subset $G \backslash \Delta_s^1$ of $M_{nod} = G \backslash \Delta_s$. If $F \in \Delta_s^1$, then from Lemma 6.13 we see that $g(F) \in P\Gamma' \backslash \mathbb{C}H^3$ is the equivalence class of the line $f_*(H_{\bar{\omega}}^{2,1}(\tilde{T})) \in \mathbb{C}H^3$, where $f : H_{\bar{\omega}}^3(\tilde{T}) \longrightarrow \mathbb{C}^{3,1}$ is any isometry that maps the projection of the integral lattice $H^3(\tilde{T}, \mathbb{Z})$ in $H_{\bar{\omega}}^3(\tilde{T})$ isomorphically to the lattice $\mathcal{E}^{3,1}$ in $\mathbb{C}^{3,1}$. Recall from (6.8) that $\tilde{T}$ is the blow-up of the singular cubic threefold $T$ associated to $F$. From (9.9) we see that $g' \circ \phi(F) \in P\Gamma' \backslash \mathbb{C}H^3$ is the equivalence class of $f'_*(H_{\bar{\omega}}^{1,0}(C_{\phi(F)})) \in \mathbb{C}H^3$, where $f' : H_{\bar{\omega}}^1(C_{\phi(F)}) \longrightarrow \mathbb{C}^{3,1}$ is any isometry which maps of the integral lattice $H^1(C_{\phi(F)}, \mathbb{Z})$ isomorphically onto $\mathcal{E}^{3,1}$. Recall that $\phi(F) \in M(6)_s$ is represented by the sextuple of zeros of the binary sextic $F_3(U^2, UV, V^2)$, where $F_3$ is as in (9.7.1), and that $C_{\phi(F)}$ is the triple cover of $\mathbb{P}^1$ branched along the sextuple $\phi(F)$. Thus to prove that $g' \circ \phi(F) = g(F)$ it is enough to produce an isomorphism $H^3_{\tilde{\omega}}(\tilde{T}) \longrightarrow H^1_{\tilde{\omega}}(C_{\phi(F)})$ that takes the projection of the integral lattice $H^3(\tilde{T}, \mathbb{Z})$ to the projection of the integral lattice $H^1(C_{\phi(F)}, \mathbb{Z})$ and which preserves the Hodge structures. For this it is enough to produce an isomorphism $H^3(\tilde{T}, \mathbb{Z}) \longrightarrow H^1(C_{\phi(F)}, \mathbb{Z})$ which commutes with the actions of the cyclic automorphisms $\sigma$ on $\tilde{T}$ and $C_{\phi(F)}$, and whose complexification is an isomorphism of Hodge structures (of bidegree $(-1, -1)$).

The latter is easy to do using the well known structure of cubic hypersurfaces with a double point. First observe (see (9.7.1) that the singular threefold $T \subset \mathbb{P}^4$ has equation

$$X_3(X_0 X_2 - X_1^2) + F_3(X_0, X_1, X_2) - Y^3 = 0.$$

If $F \in \Delta^1_s$ then the only singular point of $T$ is $(0 : 0 : 0 : 1 : 0)$. Projecting $T$ from this point to the hyperplane $X_3 = 0$, we see that $\tilde{T}$ is isomorphic to $B_{C'}\mathbb{P}^3$, the blow-up of $\mathbb{P}^3$ along the curve $C' \subset \mathbb{P}^3$ which is the complete intersection of the quadric $X_0 X_2 - X_1^2 = 0$ and the cubic $F_3(X_0, X_1, X_2) + Y^3 = 0$. See Lemma 6.5 of [15] for more details. Setting $X_0 = U^2, X_1 = UV, X_2 = V^2$, we see that $C'$ is isomorphic to the sextic curve

$$F_3(U^2, UV, V^2) - Y^3 = 0$$

in the weighted projective space $\mathbb{P}(1, 1, 2)$ with homogeneous coordinates $(U, V, Y)$. In other words, $C'$ is isomorphic to the curve $C_{\phi(F)}$, the cyclic triple cover of $\mathbb{P}^1$ branched along $\phi(F)$. Let us write simply $C$ for $C_{\phi(F)} = C'$. Thus $\tilde{T} = B_C \mathbb{P}^3$. Note that the cyclic automorphism $\sigma$ of $\mathbb{P}^4$ defined in (2.1) induces cyclic group actions on $\tilde{T}$, $\mathbb{P}^3$ and $C$.

Now let $D \subset B_C \mathbb{P}^3$ denote the preimage of $C$ under $B_C \mathbb{P}^3 \longrightarrow \mathbb{P}^3$. Then $D$ is isomorphic to the projectivized normal bundle of $C$ in $\mathbb{P}^3$, and the the cyclic automorphism $\sigma$ acts on it in a natural way. We also have natural maps

$$H^3(\tilde{T}, \mathbb{Z}) = H^3(B_C \mathbb{P}^3, \mathbb{Z}) \longrightarrow H^3(D, \mathbb{Z}) \longrightarrow H^1(C, \mathbb{Z}) \tag{9.11.1}$$

where the first arrow is induced by the inclusion of $D$ in $\tilde{T} = B_C \mathbb{P}^3$ and the second is the Gysin map (integration over the fiber). It is easy to see that both maps are isomorphisms, that they commute with the action of $\sigma$, that the complexification of the first arrow is an isomorphism of Hodge structures (of bidegree $(0, 0)$), and that the complexification of the second arrow is an isomorphism of Hodge structures (of bidegree $(-1, -1)$). (See Lemma 3.11 of [15] for a related discussion.) The composition of the two arrows gives the desired isomorphism $H^3(\tilde{T}, \mathbb{Z}) \longrightarrow H^1(C, \mathbb{Z})$, thereby completing the proof of the lemma, and of Proposition 9.6.

**Proof of Theorem 3.17**

(9.12) First observe that that Theorem 9.5 and Proposition 9.6 together imply that the bottom horizontal arrow of (3.16.1) is a homeomorphism. Since a holomorphic homeomorphism with image a normal analytic space is biholomorphic, Theorem 3.17 is proved for the bottom horizontal arrow of (3.16.1).

Now consider the middle horizontal arrow of (3.16.1). Its restriction to $M_0^m$ maps $M_0^m$ to $P\Gamma_\theta \backslash (\mathbb{C}H^4 - \mathcal{H})$ and is a proper open map of degree one; thus it is an isomorphism. Using the fact that $\mathcal{M}_s$ is the Fox completion of $\mathcal{M}_0$ over $\mathcal{C}_s$ it is easy to see that $M_s^m = G\backslash\mathcal{M}_s$ is the Fox completion of $M_0^m = G\backslash\mathcal{M}_0$ over $M_s = G\backslash\mathcal{C}_s$. Consequently the required isomorphism follows from the uniqueness of the Fox completion. In more detail, following the terminology of [22], one checks that $M_s^m$ is a spread over $M_s$ (because $\mathcal{M}_s$ is a spread over $\mathcal{C}_s$) and that it is a complete spread (because of the completeness of $\mathcal{M}_s$ over $\mathcal{C}_s$). It is obvious that $P\Gamma_\theta\backslash\mathbb{C}H^4$ is the Fox completion of $P\Gamma_\theta\backslash(\mathbb{C}H^4 - \mathcal{H})$ over $P\Gamma\backslash\mathbb{C}H^4$. From the extension theorem of §3 of [22] it follows that there is a map $P\Gamma_\theta\backslash\mathbb{C}H^4 \longrightarrow M_s^m$ which restricts on $P\Gamma_\theta\backslash(\mathbb{C}H^4 - \mathcal{H})$ to the inverse of the restriction of the middle horizontal arrow of (3.16.1) to $M_0^m$ and which covers the inverse of the bottom horizontal arrow of (3.16.1). This extension must be a continuous inverse to the middle horizontal arrow of (3.16.1). It follows that the middle horizontal arrow is an isomorphism of complex manifolds.

The top horizontal arrow can be handled in the same way: its restriction to $M_0^f$ is an isomorphism onto $\mathbb{C}H^4 - \mathcal{H}$ because the restrictions of the two top vertical arrows to these spaces are unbranched covers and the middle arrow is an isomorphism. The isomorphism on $M_s^f$ is handled by a similar extension argument for Fox completions.

The fact that $g : M_s^f \to \mathbb{C}H^4$ carries framed surfaces with exactly $k$ nodes onto $\mathcal{H}^k$ follows from the equivariance of the map and the fact that it is biholomorphic. Finally, the isomorphisms of (3.17.1) are clear because the horizontal arrows are holomorphic bijections and the Satake compactifications are normal analytic spaces [7]. This completes the proof of Theorem 3.17.

**Proof of Theorem 3.20:**

(9.13) We write $\hat{Q}$ and $\hat{Q}_\theta$ for the spaces $P\Gamma\backslash\mathbb{C}H^4$ and $P\Gamma_\theta\backslash\mathbb{C}H^4$ equipped with the orbifold structures defined in (3.19). The orbifold isomorphism $M_s^m \cong \hat{Q}_\theta$ follows because, as orbifolds, they are manifolds ((3.15) and (7.27)), and $g$ is an analytic-space isomorphism. The isomorphism $M_s \cong \hat{Q}$ requires a different technique, which could also be applied to $\hat{Q}_\theta$. Let $x \in M_s$ be represented by a form $F \in \Delta_s^k \subset \mathcal{C}_s$, and let $f \in \mathcal{F}_s$ be a framing of $F$. We may choose a small neighborhood $W$ of $F$ that is invariant under the stabilizer $J$ of $F$ in $G$. Let $\tilde{W}$ be the component of the preimage of $W$ that contains $f$. Because $J$ preserves $\Delta \cap W$, every element of $J$ lifts to an automorphism of $\tilde{W}$ preserving the preimage $\tilde{\Delta}$ of $\Delta$. These lifts normalize the group $N \cong (\mathbb{Z}/6)^k$ of covering transformations of $\tilde{W}$ over $W$, yielding a group $H = N.J$ acting on $\tilde{W}$. Choose a small 4-dimensional disk $D$ that is transverse at $f$ to the $G$-orbit of $f$ and is invariant under $H$. We will compare the two maps $D \to \mathbb{C}H^4 \to \hat{Q}$ and $D \to \mathcal{C}_s \to M_s$ and observe that they define the same orbifold structure. The quotient $D/N$ is a complex manifold, its image in $\mathcal{C}_s$ is a disk transverse at $F$ to the $G$-orbit of $F$, and by construction it is $J$-invariant. The orbifold structure of $M_s$ is given by the quotients of such invariant transverse disks in $\mathcal{C}_s$ by the local isotropy groups, so the map to $(D/N)/(H/N) = (D/N)/J$ is an orbifold chart for $M_s$. On the other hand, $D$ projects biholomorphically to a neighborhood of $g(f)$ in $\mathbb{C}H^4$ because $\mathcal{F}_s$ is a principal $G$-bundle over $\mathbb{C}H^4$. By $P\Gamma$-equivariance, we know that the subgroup $N$ of $H$ acts on $\mathbb{C}H^4$ by the 6-fold reflections in the components of $\mathcal{H}$ that pass through $g(f)$. The definition of $\hat{Q}$ as an orbifold shows that an orbifold chart for a neighborhood of $g(x)$ is the map from the complex manifold $D/N$ to $(D/N)/(H/N) = (D/N)/J$. The orbifold isomorphism is now clear.

## 10. The universal cubic surface

(10.1) In this section we prove a simple but remarkable theorem which relates two vector bundles over the moduli space of marked smooth cubic surfacs. To state it, observe that since $G$ acts freely on $\mathcal{M}_s$ (Lemma 3.14), the universal surface (2.8.1) descends to a universal marked cubic surface $\mathcal{S}^m$ over $M_s^m$. Its total space is, however, no longer smooth. In addition, the trivial $\mathbb{C}^4$-bundle in whose projective bundle $\mathcal{S}$ lies (cf. (2.8)) descends to a vector bundle $\mathbb{V}$ over $M_s^m$, and $\mathcal{S}^m \subset P\mathbb{V}$. Then we have the following.

(10.2) **Theorem.** *The bundles $TM_0^m$ and $\mathbb{V}|_{M_0^m}$ are $W(E_6)$-equivariantly isomorphic.*

As a consequence, used in the next section, symmetries of cubic surfaces are in one-to-one correspondence with symmetries of the moduli space of marked forms which fix a point:

(10.3) **Corollary.** *For each $F \in \mathcal{C}_0$, the linear action of its isotropy group on $\mathbb{C}^4$ is isomorphic to the linear action of the isotropy group of the class of $m$ (resp. $f$) on the tangent space to $M_0^m$ (resp. $M_0^f$), where $m \in M_s^m$ (resp. $f \in M_s^f$) lies over $F$.*

The theorem follows from an easily proved lemma. To state it, let $\mathcal{G} \subset T\mathcal{C}_0$ be the sub-bundle of vectors tangent to the $G$-orbits, and let $\mathcal{Q}$ denote the quotient bundle $T\mathcal{C}_0/\mathcal{G}$. Let $GL(4,\mathbb{C})$ act on $\mathcal{C}_0 \times \mathbb{C}^4$ as in the formulas of (2.17). Then we have the following.

(10.4) **Lemma.** *The bundles $\mathcal{Q}$ and $\mathcal{C}_0 \times \mathbb{C}^4$ are equivariantly isomorphic $GL(4,\mathbb{C})$-bundles.*

(10.5) Given the lemma, observe that after lifting the $GL(4,\mathbb{C})$-action as in lifting the $G$-action defined in (3.2), there is an equivariant isomorphism of the bundles lifted to $\mathcal{M}_0$. Over $\mathcal{M}_0$ there is an isomorphism $\mathcal{Q} \cong p^*TM_0^m$ of $GL(4,\mathbb{C})$-bundles, where $p: \mathcal{M}_0 \longrightarrow M_0^m$ is the orbit map, as well as an isomorphism $\mathcal{M}_0 \times \mathbb{C}^4 \cong p^*\mathbb{V}$ of $GL(4,\mathbb{C})$-bundles. Since these $GL(4,\mathbb{C})$-equivariant isomorphisms over $\mathcal{M}_0$ are equivalent to the asserted $W(E_6)$-equivariant isomorphisms over $M_0^m$, the lemma does indeed prove Theorem 10.2.

(10.6) **Proof of Lemma 10.4**: For $F \in \mathcal{C}_0$, let $J_F$, the Jacobian ideal of $F$, denote the ideal in the polynomial ring $\mathbb{C}[X_0, \cdots, X_3]$ generated by the partial derivatives of $F$, and let $R_F$ denote the quotient ring $\mathbb{C}[X_0, \cdots, X_3]/J_F$. In the notation of (10.3), $\mathcal{Q}_F = R_F^3$, and $\mathbb{C}^{4*} = R_F^1$; these are the graded pieces of $R_F$ of degrees one and three respectively. By a theorem of Macaulay, the graded ring $R_F$ satisfies Poincaré duality with a fundamental class in $R_F^4$. In particular, the pairing
$$R_F^1 \otimes R_F^3 \longrightarrow \mathbb{C}$$
defined by
$$P \otimes Q \longrightarrow Res_0 \frac{PQ\, dX_0 \wedge \cdots \wedge dX_3}{\partial F/\partial X_0 \cdots \partial F/\partial X_3} \tag{10.6.1}$$
is a perfect $GL(4,\mathbb{C})$-invariant pairing. See §3d of [12] for more details. The $R_F^i$ define $GL(4,\mathbb{C})$-bundles $R^i$ over $\mathcal{C}_0$, where $R^1 \cong \mathcal{C}_0 \times \mathbb{C}^{4*}$ and $R^3 \cong \mathcal{Q}$. The pairing (10.6.1) provides a $GL(4,\mathbb{C})$-equivariant isomorphism between $R^3$ and $R^{1*}$, hence the required $GL(4,\mathbb{C})$-equivariant isomorphism between $\mathcal{Q}$ and $\mathcal{C}_0 \times \mathbb{C}^4$.

## 11. Automorphisms of Cubic Surfaces

(11.1) Theorem 2.20 identifies the moduli space $M_0$ of smooth cubic surfaces with $P\Gamma\backslash(\mathbb{C}H^4 - \mathcal{H})$. However, the isomorphism is far from explicit, and so it is natural to ask which surfaces correspond to which points of the ball quotient. In this section we solve this problem for the two most symmetric cubic surfaces, the Fermat cubic defined in $\mathbb{P}^3$ by $\sum_{i=0}^3 X_i^3 = 0$, and the diagonal surface defined in $\mathbb{P}^4$ by $\sum_{i=0}^4 Y_i = \sum_{i=0}^4 Y_i^3 = 0$. The idea is that smooth symmetric surfaces correspond to points in $\mathbb{C}H^n - \mathcal{H}$ with nontrivial stabilizers in $P\Gamma$. This assertion is equivalent to the isomorphism of orbifolds in Theorem 2.20, and is in analogy with the case of cubic curves, where the curves with symmetry correspond to orbifold points of the moduli space.

(11.2) To make this last point explicit, suppose that $f \in \mathcal{F}_0$ lies over $F \in \mathcal{C}_0$ and that $\alpha \in G$ leaves $F$ invariant. There is a unique element $\hat{\alpha} \in P\Gamma$ such that $\hat{\alpha}f = \alpha f$, namely, $\hat{\alpha} = f \circ \alpha^* \circ f^{-1}$, where $\alpha^*$ indicates the action of $\alpha$ on $\Lambda(T)$. (The uniqueness relies on the fact that both $G$ and $P\Gamma$ act freely on $\mathcal{F}_0$.) The map $\alpha \mapsto \hat{\alpha}$ defines an isomorphism from the stabilizer in $G$ of $F$ to the stabilizer in $P\Gamma$ of the image of $f$ in $\mathbb{C}H^4 - \mathcal{H}$ under the period map. We call $\hat{\alpha}$ the transform of $\alpha$ (with respect to $f$).

(11.3) The most important kind of symmetry of a cubic surface is a biflection of $\mathbb{P}^3$. Biflections play a large role in Segre's investigation, §§98-100 of [44], where they are called harmonic homologies. Segre showed that if a cubic surface admits a nontrivial symmetry then it admits a biflection, and that in almost all cases the full group of the surface is generated by biflections. We now compute the transform of a biflection and then give explicit coordinates for the images of the diagonal and Fermat surfaces under the period map.

**(11.4) Lemma.** *If $\alpha \in GL(4, \mathbb{C})$ is a biflection then any transform $\hat{\alpha}$ of $\alpha$ acts on $\mathbb{C}H^4$ as the biflection in some long root $r$.*

(11.5) **Proof:** Recall the identification in Theorem 10.2 of the tangent bundle to $M_0^m$ with the $\mathbb{C}^4$-bundle $\mathbb{V}$. From that identification, and specifically from (10.3), one sees that $\alpha$ acts on $\mathbb{C}^4$ in the same way it acts on the tangent space to $M_0^m$, hence in the same way as it acts on the tangent space to $\mathbb{C}H^4$ at $g(f)$, where $f$ is a framing of some smooth form preserved by $\alpha$. Thus $\hat{\alpha}$ acts on $\mathbb{C}H^4$ as a biflection. From $\hat{\alpha} = f \circ \alpha^* \circ f^{-1}$, it follows that the action of $\hat{\alpha}$ on $\mathbb{C}H^4$ may be represented by an isometry of $\Lambda$ of order 2. Such an isometry can only be a reflection of $\Lambda$, and a simple arithmetic argument (Lemma 8.1 of [1]) shows that any reflection of a unimodular Eisenstein lattice is a reflection in a lattice vector $r$ of norm $\pm 1$ or $\pm 2$. The case of negative norm is excluded because such reflections do not act on $\mathbb{C}H^4$ as biflections (each has an isolated fixed point). The case of norm 1 is excluded because the image of $f$ would, as a fixed point of $\hat{\alpha}$, lie in $\mathcal{H}$. This is impossible because $F$ is smooth. Thus $r$ must have norm 2.

**(11.6) Theorem** *There is a single $\Gamma$-orbit of vectors of $\Lambda$ of norm $-5$ that are orthogonal to no short roots. The corresponding point in $P\Gamma\backslash\mathbb{C}H^4$ is the image under the period map of the diagonal surface.*

(11.7) **Proof:** Suppose $v \in \Lambda$ has norm $-5$. Then $\Lambda$ contains the sublattice $v^\perp \oplus \langle v \rangle$, and the projections into the spans of the two summands define an $\mathcal{E}$-module isomorphism (the gluing map) of $(v^\perp)'/v^\perp$ with $\langle v\rangle'/\langle v\rangle$. Here the prime indicates the dual lattice. The norm of an element of one of these dual quotients is well-defined modulo 1, and its sum with the norm of its image under the gluing map is 0 modulo 1. This allows one to compute the norms (mod 1) of the elements of $(v^\perp)'/v^\perp$. This in turn allows one to verify that $v^\perp \oplus \langle w\rangle$ may be enlarged in an essentially unique way to a unimodular lattice $L$, where $\langle w\rangle$ denotes a 1-dimensional lattice over $\mathcal{E}$ with a generator $w$ of norm 5. Since $\mathcal{E}^5$ is the only positive-definite unimodular $\mathcal{E}$-lattice in five dimensions, $L$ is isomorphic to $\mathcal{E}^5$. One can reverse the procedure to recover $v$ from $w$. This establishes a bijection between $\text{Aut}(\Lambda)$-orbits of norm $-5$ vectors in $\Lambda$ with $\text{Aut}(\mathcal{E}^5)$-orbits of norm 5 vectors in $\mathcal{E}^5$, the bijection being given by the relation of having isometric orthogonal complements. A trivial calculation shows that there is a unique orbit of norm 5 vectors in $\mathcal{E}^5$ that are orthogonal to no short roots: a representative is $(1,1,1,1,1)$. Therefore, $\text{Aut}\,\Lambda$ acts with a single orbit on the norm $-5$ vectors of $\Lambda$ that are orthogonal to no short roots. One can check that $v = (3,1,1,1,1)$ is orthogonal to no short roots, so it is a representative for this orbit. Since $\text{Aut}\,\Lambda = \Gamma \times \{\pm I\}$, to show that $\Gamma$ acts transitively on this $(\text{Aut}\,\Lambda)$-orbit it suffices to show that $\Gamma$ contains an element exchanging $v$ and $-v$. The product of the central involution and the biflection in any long root of $v^\perp$, say $(0,0,0,1,-1)$, accomplishes this and has spinor norm $+1$. Therefore it lies in $\Gamma$. This proves the first part of the theorem.

(11.8) Now suppose that $f$ is a framed cubic form whose underlying form $F$ defines a copy of the diagonal surface, which has symmetry group $S_5$. Let $x = g(f) \in \mathbb{C}H^4$, and observe that it is stabilized by a group $\hat{S}_5$. We may choose coordinates for $x^\perp$ in $\mathbb{C}^{4,1}$ as five coordinates which sum to zero, with the ten biflections given as transpositions of coordinates. With respect to these coordinates there is a long root of $\Lambda$ proportional to $r_1 = (1,-1,0,0,0)$, and after scaling the coordinate system we may suppose that $\Lambda$ actually contains this vector. Then $\Lambda$ also contains a long root proportional to $r_2 = (0,1,-1,0,0)$; its inner product with $r_1$ must be a unit of $\mathcal{E}$, so that by replacing this root by a scalar multiple of itself we may suppose that $\Lambda$ contains $r_2$. Continuing in this fashion we may suppose that $\Lambda$ contains $r_3 = (0,0,1,-1,0)$ and that $r_4 = (0,0,0,1,-1)$. These four roots generate an $A_4$ root lattice. Since it has determinant 5, which is square-free (indeed prime) in $\mathcal{E}$, it is a maximal integral 4-dimensional $\mathcal{E}$-lattice. Since $\Lambda$ is unimodular, its orthogonal complement is generated by an element of norm $-5$. Since $x \notin \mathcal{H}$, we have established that $x \in \mathbb{C}H^4$ corresponds to one of the norm $-5$ vectors considered above.

**(11.9) Theorem.** *There is a unique $\Gamma$-orbit of vectors $v \in \Lambda$ of norm $-3$ which are orthogonal to no short roots. The corresponding point in $P\Gamma\backslash\mathbb{C}H^4$ is the image under the period map of the Fermat cubic.*

**Proof:** If $v$ is not primitive then it is a multiple of a norm $-1$ vector, so that $v^\perp \cong \mathcal{E}^{4,1}$ contains short roots. So suppose $v$ is primitive. As in the previous argument we may compute the norms modulo 1 of the elements of $(v^\perp)'/v^\perp$, and in this way we find that $v^\perp$, of determinant 3, has index 3 in a unimodular lattice. Since the only positive-definite unimodular $\mathcal{E}$-lattice in 4 dimensions is $\mathcal{E}^4$, $v^\perp$ is a sublattice of $\mathcal{E}^4$ of index 3 containing no short roots. This implies that $v^\perp$ is isometric

to
$$D_4(\theta) = \{(z_1, \ldots, z_4) \in \mathcal{E}^4 \mid z_1 + z_2 + z_3 + z_4 \equiv 0 \pmod{\theta}\}.$$

One can classify the unimodular lattices containing $v^\perp \oplus \langle v \rangle$ by considering the possible gluings of $D_4(\theta)'/D_4(\theta)$ to $\langle v \rangle'/\langle v \rangle$. Up to isometry of the summands $v^\perp$ and $\langle v \rangle$, there is a unique such unimodular lattice in which $v$ is orthogonal to no short roots. If $w$ is another norm $-3$ vector in $\Lambda$ orthogonal to no short roots, then the arguments above show that $w^\perp$ is isometric to $v^\perp$, and the essential uniqueness of the gluing map shows that the isometry $w^\perp \oplus \langle w \rangle \cong v^\perp \oplus \langle v \rangle$ extends to $\Lambda$. We have shown that $\operatorname{Aut}\Lambda$ has only one orbit of primitive norm $-3$ vectors that are orthogonal to no short roots. One can check that $v = (2 - \bar{\omega}, 1, 1, 1, 1)$ is orthogonal to no short roots, so it represents this orbit. The argument of (11.7) shows that $\Gamma$ also acts transitively on these vectors.

(11.10) Now suppose that $f$ is a framed cubic form with underlying form $F$ defining a copy of the Fermat surface, and let $x = g(f)$ as before. The stabilizer of $F$ in $G$ is the complex reflection group $(\mathbb{Z}/3)^3{:}S_4$, generated by the biflections in the long roots of $D_4(\theta)$. Mimicking the previous proof, one can find a copy of $D_4(\theta)$ in $x^\perp$. Since $x \notin \mathcal{H}$, and since any 4-dimensional enlargement of $D_4(\theta)$ that is integral contains short roots, $D_4(\theta)$ is all of $x^\perp$. Therefore $D_4(\theta)^\perp$ is spanned by a vector of norm $-3$, and indeed one of the sort we have just considered.

(11.11) **Remark.** One can use a known property of cubic surfaces to show that $P\Gamma$ provides a counterexample to a natural conjecture regarding complex reflection groups. It is known that there is a smooth cubic surface whose symmetry group is not generated by complex reflections of $\mathbb{P}^3$ (see [44], pp. 151–152). Therefore the stabilizer in $P\Gamma$ of a corresponding point $x$ in $\mathbb{C}H^4$ is not generated by complex reflections. This occurs despite the fact that $P\Gamma$ itself is generated by complex reflections. This phenomenon contrasts sharply with both real hyperbolic reflection groups and finite complex reflection groups. In both these settings, the stabilizer in the reflection group of any vector or point of projective space is itself generated by reflections.

(11.12) We observed above that any cubic surface admitting a nontrivial symmetry admits a biflection. By Lemma 11.4 this implies that the points of $\mathbb{C}H^4$ corresponding to symmetric cubic surfaces may be found along the orthogonal complements of long roots of $\Lambda$. The following theorem shows that all the symmetric surfaces occur along the orthogonal complement of any chosen long root.

**(11.13) Theorem.** *$\Gamma$ acts transitively on the long roots of $\Lambda$.*

**Proof:** The orthogonal complement of such a root $r$ contains a norm 0 vector: consider the real part of $h$ and use Meyer's Theorem [36]. Therefore by Theorem 7.21 we may take $r$ to have the form $(x, y, 0; 0, z)$ in the coordinates of (7.7), where $x$ and $y$ are units of $\mathcal{E}$. Applying elements of $\Gamma$, we may take $x = y = 1$. Then, applying a translation, we may take $z = 0$, completing the proof.

## 12. Notation

$\operatorname{Aut}^+ V$, $\operatorname{Aut}^+ \Lambda$, $\operatorname{Aut}^+ I_{4,1}^{\mathcal{E}}$, $\operatorname{Aut}^+ I_{3,1}^{\mathcal{E}}$
               automorphisms of spinor norm 1; (2.12)

| | |
|---|---|
| $\mathcal{C}, \mathcal{C}_0, \mathcal{C}_s, \mathcal{C}_{ss}$ | space of nonzero (smooth, stable, semistable) cubic forms; (2.1) and (3.1) |
| $\tilde{\mathcal{C}}_0$ | a cover of $\mathcal{C}_0$; (2.10) |
| $\mathbb{C}H(W), \mathbb{C}H^n$ | complex hyperbolic spaces; (2.15) |
| $D$ | central subgroup of $GL(4,\mathbb{C})$ of order 3; (2.17) |
| $\Delta, \Delta_s, \Delta_{ss}, \Delta_s^k, \Delta_{ss}^{a,b}$ | |
| | the discriminant locus, its stable and semistable parts, and |
| | strata therein; (2.1), (3.1), (3.5) and (3.7) |
| $\mathcal{E}$ | The Eisenstein integers $\mathbb{Z}[\omega]$; (2.2) |
| $\mathcal{E}^{n,1}$ | standard $\mathcal{E}$-lattice of signature $(n,1)$; (2.7) |
| $\eta(S), \eta$ | hyperplane class in $L(S)$, and the "standard copy" of it in $L$; (3.2) |
| $F$ | a cubic form; determines $S$ and $T$ implicitly; (2.1) |
| $\mathcal{F}_0, \mathcal{F}_s$ | space of framed smooth (stable) cubic forms; (3.9) and (3.10) |
| $G$ | $GL(4,\mathbb{C})/D$ |
| $g_0, g$ | period maps; (2.16), (2.18), (3.16) and §8 |
| $\Gamma, \Gamma_\theta$ | monodromy group and congruence subgroup; (2.11) and (3.12) |
| $h$ | hermitian form on $\Lambda(T)$; (2.3) |
| $h'$ | hermitian form on $H^3(T,\mathbb{C})$; (2.4) |
| $\mathcal{H}, \mathcal{H}^k$ | hyperplane arrangement in $\mathbb{C}H^4$, and strata therein; (2.19) and (3.16) |
| $I_{n,1}^\mathcal{E}$ | another coordinate system for $\mathcal{E}^{n,1}$; (7.7) |
| $id$ | the identity map |
| $L, L(S), L(\mathcal{S})$ | $\mathbb{Z}^{1,6}$, the lattice $H^2(S) \cong L$, and the local system of these lattices; (3.2) |
| $\Lambda, \Lambda(T), \Lambda(\mathcal{T})$ | $\mathcal{E}^{4,1}$, the $\mathcal{E}$-lattice associated to $T$, and the local system of |
| | these $\mathcal{E}$-lattices; (2.2), (2.7), (2.9), and for singular $T$, (5.8) |
| $M, M_s, M_0, M^m, M_s^m, M_0^m, M_s^f, M_0^f$ | |
| | moduli spaces; (2.18) and (3.15) |
| $\mathcal{M}_0, \mathcal{M}_s, \mathcal{M}_{ss}$ | space of marked smooth (stable, semistable) cubic forms; (3.2)–(3.4) |
| $p$ | projection $T \to \mathbb{P}^3$; (2.1) |
| $\pi$ | projection $\mathcal{S} \to \mathcal{C}$ or $\mathcal{T} \to \mathcal{C}$; (2.8) |
| $q$ | quadratic form on $V$; (2.12) |
| $\rho_0, \rho$ | monodromy representations; (2.10) and (2.11) |
| $S, S'$, etc. | cubic surface defined by $F, F'$, etc.; (2.1) |
| $\mathcal{S}, \mathcal{S}_0$ | universal family of (smooth) cubic surfaces; (2.8) |
| $\sigma$ | branched covering transformation $T \to T$; (2.1) |
| $T, T'$, etc. | cyclic cubic threefold defined by $F, F'$, etc.; (2.1) |
| $\mathcal{T}, \mathcal{T}_0$ | universal family of (smooth) cyclic cubic threefolds; (2.8) |
| $\theta$ | $\omega - \bar\omega = \sqrt{-3} \in \mathcal{E}$; (2.3) |
| $V, V(S), V(T)$ | $\mathbb{F}_3^5$, and finite vector spaces associated to $S$ and $T$; (2.12) and (4.8) |
| $Z$ | isometric embedding $\Lambda(T) \to H^3_{\bar\omega}(T,\mathbb{C})$; (2.4.2) |
| $\omega$ | a primitive cube root of unity; (2.1) |
| $\Omega$ | symplectic pairing on $H^3(T,\mathbb{Z})$; (2.3) |

D. Allcock, Department of Mathematics, Harvard University, Cambridge, MA, USA.
*allcock@math.harvard.edu*

J. Carlson, Department of Mathematics, University of Utah, Salt Lake City, UT, USA.
*carlson@math.utah.edu*

D. Toledo, Department of Mathematics, University of Utah, Salt Lake City, UT, USA.
*toledo@math.utah.edu*